\pdfoutput=1
\RequirePackage{fix-cm}
\documentclass[12pt,reqno]{amsart}
\usepackage{amsmath,amssymb,amsfonts,latexsym,amsthm,enumerate,mathabx}

\usepackage[T1]{fontenc}
\usepackage{lmodern}

\usepackage{xfrac}
\usepackage{mathrsfs}

\DeclareFontFamily{U}{rsfs}{\skewchar\font127 }
\DeclareFontShape{U}{rsfs}{m}{n}{%
   <-6> rsfs5
   <6-8> rsfs7
   <8-> rsfs10
}{}

\DeclareFontFamily{U}{matha}{\skewchar\font127 }
\DeclareFontShape{U}{matha}{m}{n}{%
   <-6> matha5
   <6-8> matha7
   <8-> matha10
}{}

\DeclareFontFamily{U}{mathb}{\skewchar\font127 }
\DeclareFontShape{U}{mathb}{m}{n}{%
   <-6> mathb5
   <6-8> mathb7
   <8-> mathb10
}{}

\DeclareFontFamily{U}{mathx}{\skewchar\font127 }
\DeclareFontShape{U}{mathx}{m}{n}{%
   <-6> mathx5
   <6-8> mathx7
   <8-> mathx10
}{}

\usepackage[margin=1in]{geometry}
\usepackage{stackrel}

\usepackage{pkg3}
\usepackage{kbordermatrix}

\newcommand{\lit}{L}
\newcommand{\ulit}{\vec L}

\newcommand{\qp}{\acute{q}}
\newcommand{\Qp}{\acute{Q}}

\newcommand{\cc}{7}

\newcommand{\kacirc}{\ka^\circ}
\newcommand{\mucirc}{\mu^\circ}

\newcommand{\tree}{T_{d,k}}
\newcommand{\etree}{\acute{T}_{d,k}}
\newcommand{\Gndk}{\mathcal{G}_{n,d,k}}

\newcommand{\nuaux}{\bm{\nu}}

\newcommand{\nubar}{\overline{\nu}}

\newcommand{\bemax}{\be_{\max}}

\newcommand{\dubd}{d_{\mathrm{ubd}}}
\newcommand{\dlbd}{d_{\mathrm{lbd}}}
\newcommand{\dfm}{d_\square} 

\newcommand{\qrig}{q}
\newcommand{\qf}{q_\free}
\newcommand{\vr}{v_\rig}

\newcommand{\bpdzg}{\DOT{z}_g}
\newcommand{\bpdzh}{\DOT{z}_h}
\newcommand{\bphzg}{\HAT{z}_g}
\newcommand{\bphzh}{\HAT{z}_h}

\newcommand{\hqf}{\hq_\free}
\newcommand{\hqo}{\hq_\one}
\newcommand{\hqz}{\hq_\zro}
\newcommand{\dqf}{\dq_\free}
\newcommand{\dqo}{\dq_\one}
\newcommand{\dqz}{\dq_\zro}

\newcommand{\Dq}{D_{\qrig}}

\newcommand{\grigid}[1][]{\mathbf{g}_\rig^{#1}}

\newcommand{\naesol}{\Om}

\newcommand{\upmsg}{\smash{\acute{\bm{u}}}}
\newcommand{\downmsg}{\smash{\grave{\bm{\upsilon}}}}

\newcommand{\mf}{m_\mathrm{f}}
\newcommand{\Rne}{\mathscr{R}}
\newcommand{\Avp}{\widetilde{\text{\textsc{a}}}_v}
\newcommand{\Av}{\text{\textsc{a}}_v}

\newcommand{\mA}{m_{\text{\sc a}}}

\newcommand{\dtq}{\bm{\DOT{q}}}
\newcommand{\htq}{\bm{\HAT{q}}}
\newcommand{\hqeq}{\htq_=}
\newcommand{\hqne}{\htq_{\ne}}
\newcommand{\dqeq}{\dtq_=}
\newcommand{\dqrr}{\dtq_\rr}
\newcommand{\dqne}{\dtq_{\ne}}
\newcommand{\hqrf}{\htq_{\rf}}
\newcommand{\hqfr}{\htq_{\fr}}
\newcommand{\hqzz}{\htq_{\zz}}
\newcommand{\hqoo}{\htq_{\oo}}
\newcommand{\hqzo}{\htq_{\zo}}
\newcommand{\hqzf}{\htq_{\zf}}
\newcommand{\hqfz}{\htq_{\fz}}
\newcommand{\hqof}{\htq_{\of}}
\newcommand{\hqfo}{\htq_{\fo}}
\newcommand{\hqff}{\htq_{\ff}}
\newcommand{\hqoz}{\htq_{\oz}}

\newcommand{\dqzz}{\dtq_\zz}
\newcommand{\dqzo}{\dtq_\zo}
\newcommand{\dqzf}{\dtq_\zf}
\newcommand{\dqoz}{\dtq_\oz}
\newcommand{\dqoo}{\dtq_\oo}
\newcommand{\dqof}{\dtq_\of}
\newcommand{\dqfz}{\dtq_\fz}
\newcommand{\dqfo}{\dtq_\fo}
\newcommand{\dqff}{\dtq_\ff}
\newcommand{\dqrd}{\dtq_\rd}
\newcommand{\dqrf}{\dtq_\rf}
\newcommand{\dqfr}{\dtq_\fr}
\newcommand{\dqdr}{\dtq_{\cdot\rig}}

\title{Satisfiability threshold for random regular \textsc{nae-sat}}

\author[J.~Ding]{${}^{*}$Jian Ding}
\author[A.~Sly]{${}^{\dagger}$Allan Sly}
\author[N.~Sun]{${}^{\ddagger}$Nike Sun}
\thanks{Research supported by 
${}^{*}$NSF grant DMS-1313596;
${}^{\dagger}$Sloan Research Fellowship; ${}^{\ddagger}$NDSEG and NSF GRF}

\date{\today}

\begin{document}

\begin{abstract}
We consider the random regular $k$-\textsc{nae-sat} problem with $n$ variables each appearing in exactly $d$ clauses. For all $k$ exceeding an absolute constant $k_0$, we establish explicitly the satisfiability threshold $d_\star\equiv d_\star(k)$. We prove that for $d<d_\star$ the problem is satisfiable with high probability while for $d>d_\star$ the problem is unsatisfiable with high probability. If the threshold $d_\star$ lands exactly on an integer, we show that the problem is satisfiable with probability bounded away from both zero and one. This is the first result to locate the exact satisfiability threshold in a random constraint satisfaction problem exhibiting the condensation phenomenon identified by Krzaka{\l}a et al.\ (2007). Our proof verifies the one-step replica symmetry breaking formalism for this model. We expect our methods to be applicable to a broad range of random constraint satisfaction problems and combinatorial problems on random graphs.
\end{abstract}

\maketitle

\section{Introduction}

Given a Boolean formula in conjunctive normal form (i.e., expressed as an \textsc{and} of \textsc{or}s), a \emph{not-all-equal-\textsc{sat}} (\textsc{nae-sat}) solution is an assignment $\ux$ of literals to variables such that both $\ux$ and its negation $\neg\ux$ evaluate to \textsc{true}. A $k$-\textsc{nae-sat} problem is one in which each clause involves exactly $k$ literals.

The $k$-\textsc{nae-sat} problem is a symmetrized version of $k$-\textsc{sat}. A major direction of research has concerned the
large-system limit of \emph{random problem instances}, seeking to establish typical behavior and phase transitions. In particular, much effort has been directed towards  locating the \emph{satisfiability transition}: the critical density $\al_\star$ where solutions cease to exist \cite{achlioptas2005rigorous}.

Random $k$-\textsc{nae-sat} is perhaps the simplest of a broad universality class of sparse random constraint satisfaction problems (\textsc{csp}s) --- including $k$-\textsc{sat}, colorings, and independent sets on random graphs --- which has been intensively studied in statistical physics, combinatorics and theoretical computer science. Statistical physicists~\cite{PNAS19062007, mezard1985replicas, mezard2002random} have described these problems by \emph{replica symmetry breaking} (\textsc{rsb}): above a certain \emph{condensation} threshold $\al_\mathrm{c}$ which is strictly below $\al_\star$, the solution space is dominated by large clusters. This deep but non-rigorous theory makes explicit predictions for the satisfiability thresholds of these models, the one-step replica symmetry breaking solution. However, no such prediction has been rigorously verified in a \textsc{csp} exhibiting condensation, with all previous satisfiability bounds leaving a constant gap.

In this paper we consider random $d$-regular $k$-\textsc{nae-sat}, in which each variable is involved in exactly $d$ clauses, and clause literals are chosen uniformly at random. We establish the following sharp satisfiability threshold, the first of its kind among
this class of \textsc{csp}s:

\bThm\label{t:main}
For $k\ge k_0$ there is a threshold $d_\star\equiv d_\star(k)$,
given by the largest zero of the explicit function \eqref{e:explicitIntro},
such that the probability for a random $d$-regular $k$-\textsc{nae-sat} instance to be solvable tends to one for $d<d_\star$, and tends to zero for $d>d_\star$.
\eThm

The threshold $d_\star$ is given by the largest zero of the function
\beq\label{e:explicitIntro}
\bPhistar(d)
\equiv \log2-\log(2-\qrig)
	-d(1-k^{-1}-d^{-1})\log[ 1- 2 (q/2)^k]
	+(d-1)\log[1-(q/2)^{k-1}]
\eeq
where $\qrig=\qrig(d)$ is the unique solution in the interval $[1-2^{-k},1]$ of
\[ d=1+\Big(\log \smf{2(1-\qrig)}{2-\qrig}\Big)
	/ \Big( \log \smf{1-2(\qrig/2)^{k-1}}{1-(\qrig/2)^{k-1}} \Big).
\]
We will find (see Propn.~\ref{p:explicit})
that $\bPhistar$ is decreasing with a unique zero on the interval $(2^{k-1}-2)k\log2 \le d\le 2^{k-1}k\log2$.

As the threshold is given by the root of an equation, it is possible for $d_\star$ to be integer-valued, though we have no reason to believe that this ever occurs. Nevertheless, we also address this hypothetical possibility by showing that if $d=d_\star$ then the probability for the \textsc{nae-sat} instance to be solvable is asymptotically bounded away from both zero and one. This completes the characterization of the satisfiability transition.

The methods developed in this paper offer a new approach to tackling other problems in the same class and establishing exact thresholds. Indeed, in a companion paper \cite{dss-is} we consider the maximum independent set problem on random regular graphs, where we determine the explicit threshold, and furthermore show tight concentration of the maximum independent set size about the threshold value.

Previous work on the satisfiability transition has identified sharp thresholds in models not exhibiting condensation, e.g. \textsc{xor-sat} \cite{MR1972120, pittelsorkin}. The $2$-\textsc{sat} satisfiability transition is also much simpler, and can be identified by a branching process argument \cite{2677892sat, MR1423858, fernandez2001random}. See also~\cite{MR2518205} for detailed discussions of these problems. Previous work on \textsc{nae-sat} has centered on the Erd\doubleacute{o}s--R\'enyi version in which variables are included in clauses independently at random, with a series of improving bounds on the satisfiability transition \cite{MR2263010, coja2012condensation, MR2961553}.

Shortly prior to the posting of this paper, A.~Coja-Oghlan posted a paper~\cite{cojaoghlanksat} on a different symmetrization of regular $k$-\textsc{sat} in which a $2$-clause joins each consecutive pair of variables, forcing them to take opposite literals. While not establishing a satisfiability threshold, his paper establishes a \textsc{1rsb}-type formula for the existence of solutions that satisfy all but $o(n)$ clauses. His approach of modeling clusters of configurations is similar to our own.

\subsection{Notation}

\begin{figure}[ht]
\includegraphics[height=1.2in,trim=.6in .6in .7in .6in,clip]{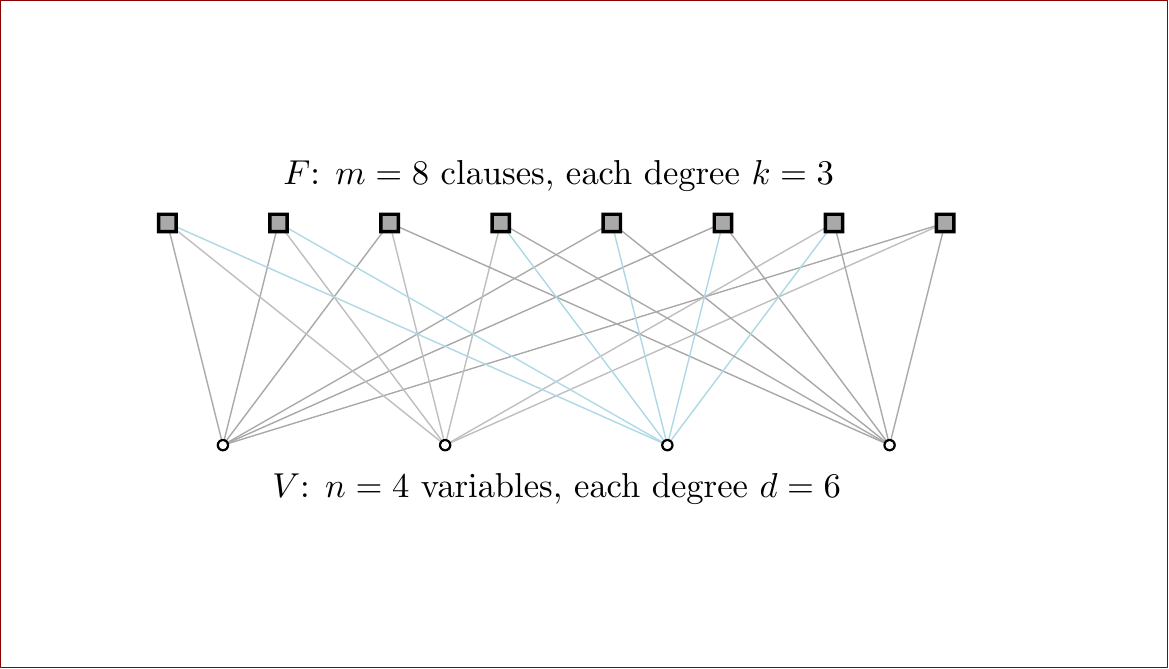}
\caption{$(d,k)$-regular bipartite factor graph}
\label{f:bfg}
\end{figure}

Throughout this paper $G$ denotes a $(d,k)$-regular bipartite graph with bipartition $(V,F)$, where $V\equiv\set{v_1,\ldots,v_n}$ is the set of degree-$d$ vertices (variables), $F\equiv\set{a_1,\ldots,a_m}$ is the set of degree-$k$ vertices (clauses), and every edge $e\in E$ is of form $e=(av)$ with $a\in F$, $v\in V$ (Fig.~\ref{f:bfg}). A \emph{variable assignment} is a configuration $\ux\in\set{\zro,\one}^V$, and a \emph{literal assignment} is a configuration $\ulit\in\set{\zro,\one}^E$. Let $\pd a=(v_1,\ldots,v_k)$ denote the $k$ vertices adjacent to $a$, with repetition if the graph has multi-edges. The evaluation of $\ux$ by clause $a$ is the vector $$(\ulit\ux)_a\equiv(\lit_{av}\oplus x_v)_{v\in\pd a}\in\set{\zro,\one}^k$$ where $\oplus$ indicates addition modulo two. We write $\neg x\equiv x \oplus \one$.

\bdfn
A variable assignment $\ux\in\set{\zro,\one}^V$ is a \emph{\textsc{sat} solution} for $(G,\ulit)$ if $(\ulit\ux)_a$ is not the identically-$\zro$ vector $(\zro^k)$ for any $a\in F$. The assignment $\ux$ is a \emph{not-all-equal-\textsc{sat}} (\textsc{nae-sat}) \emph{solution} for $(G,\ulit)$ if both $\ux$ and $\neg\ux$
 are \textsc{sat} solutions for $(G,\ulit)$.
\edfn

We hereafter write $G\equiv \Gndk$ to indicate that $G$ is
chosen according to the configuration model for uniformly random $(d,k)$-regular bipartite (multi-)graphs with $n$ degree-$d$ vertices, with $\ulit$ a uniformly random literal assignment.

\subsection{Outline of proof}
Our proof of Thm.~\ref{t:main} has two main parts which we now describe. The first is an application of the moment method: if $Z_n$ ($n\ge1$) are non-negative random variables then the Cauchy--Schwarz inequality implies $(\E Z_n)^2 / \E[Z_n^2] \le \P(Z_n>0)$, so if $(\E Z_n)^2\asymp\E[Z_n^2]$ then $Z_n>0$ \emph{with positive probability} in the limit $n\to\infty$.\footnote{The event $B_n$ is said to hold with positive probability if $\liminf_n\P(B_n)>0$.} On the other hand, if the $Z_n$ are integer-valued with $\E Z_n\to0$, then $\P(Z_n>0)\to0$ by Markov's inequality.

It is most natural to apply the moment method with the \textsc{nae-sat} partition function
\[
Z\equiv Z_{G,\ulit} \equiv |\naesol_{G,\ulit}|,\quad
	\naesol_{G,\ulit}\subseteq\set{\zro,\one}^V
	\text{ the set of {\sc nae-sat} solutions on $(G,\ulit)$.}
\]
We emphasize that $\naesol_{G,\ulit}$ is a random subset of $\set{\zro,\one}^V$ determined by $(G,\ulit)$. By symmetry, $\P(\ux\in\naesol)$ is constant over $\ux\in\set{\zro,\one}^V$,
so it suffices to consider the identically-$\zro$ vector $\ux=\vec\zro$:
\beq\label{e:ez}
\begin{array}{l}
\E Z
= 2^n\P(\vec\zro\in\naesol_{G,\ulit})
=2^n(1-2/2^k)^m = \exp\{ n\,\Phi_k(d) \}\\
{}\qquad\qquad\text{with }
\Phi_k(d)\equiv \log 2
	+(d/k)\log(1-2/2^k) ].
\end{array}
\eeq
For fixed $k$ the rate function $\Phi_k(d)$ is clearly decreasing in $d$, with unique zero at
\beq\label{e:first.moment.threshold}
\dfm\equiv \smf{k\log 2}{ -\log(1-2/2^k) }
=[2^{k-1}-\tf12
	- \tf{1}{6 \cdot 2^k}
	- O(4^{-k})]\,k\log2
\le 2^{k-1}k\log 2 \equiv\dubd
\eeq
In \S\ref{s:assignments} we will see that a rather straightforward application of the second moment method on $Z$ gives the following

\bppn\label{p:sat.below}
For $k\ge k_0$
and
$d\le\dlbd\equiv\dlbd(k)\equiv(2^{k-1}-2)k\log 2$,
$\E[Z^2]\asymp_k (\E Z)^2$, implying $(G,\ulit)$ has a \textsc{nae-sat} solution
with positive probability as $n\to\infty$.
\eppn

However the second moment method \emph{fails} for small $\rho$: there is a regime of $\rho$ in which both $\E Z$ and $\E[Z^2]/(\E Z)^2$ are exponentially large in $n$, giving no information on the limiting behavior of $\P(Z>0)$. In a sense, this issue characterizes this class of \textsc{csp}s.

To determine the exact threshold, we introduce a \emph{frozen model} with spins $\zro$ and $\one$ together with a third spin $\free$ (``free''), such that each configuration $\ueta\in\set{\zro,\one,\free}^V$ effectively encodes an entire \emph{cluster} of ($\zro/\one$-valued) \textsc{nae-sat} solutions. Roughly speaking, in a frozen model configuration, the $\zro/\one$ spins indicate variables which are rigid (cannot be flipped by local perturbations) while the $\free$ spin indicates variables which can be flipped by making local changes. In \S\ref{s:assignments} we show how to project \textsc{nae-sat} solutions to frozen configurations by a certain ``coarsening'' algorithm, and show further that this only produces 
frozen configurations with a very low density of frees.
Conversely, in \S\ref{s:sat} we explain how to recover an \textsc{nae-sat} assignment from such a frozen configuration.
This reduces the proof of Thm.~\ref{t:main}
to showing a sharp threshold for the existence of
frozen configurations with a low density of frees.

We establish this by the second moment method applied to the partition function $\ZZ$ of the frozen configurations. In \S\ref{s:first}-\ref{s:nd} we prove
\bThm\label{t:moments}
For $k\ge k_0$ and $\dlbd\le d\le \dubd$,
there is an explicit constant $\bPhistar\equiv\bPhistar_k(d)$
which is decreasing in $d$ such that
$$\EZZ \asymp_k \exp\{ n\,\bPhistar \},\quad
\E[\ZZ^2] \lesssim_k(\EZZ)^2+n^{O(1)}\,\EZZ.$$
\eThm

The proof of Thm.~\ref{t:moments} comprises a large portion of the present paper. The first moment is addressed in \S\ref{s:first}, where we identify the exact local neighborhood profile that gives the maximal contribution to the expectation. This is done by a Bethe variational principle which relates stationary points of the rate function to fixed points of certain tree recursions. A major technical difficulty is the high dimensionality of the maximization problem, and the possibility of multiple stationary points which must be ruled out. This is done by delicate \emph{a priori} estimates which allow us to reduce the dimensionality by certain symmetry conditions. 

The second moment can be understood in the same framework by regarding it as the first moment of the pair model, but clearly the dimensionality is substantially increased. We show in \S\ref{s:two} that the dominant contribution comes from two local maximizers: one corresponding to pairs whose overlap distribution looks like a product measure, and the other corresponding to pairs which are perfectly correlated --- in each case, with both marginals given by the first moment maximizer. The results of \S\ref{s:first}~and~\ref{s:two} control the moments up to polynomial prefactors, which are determined in \S\ref{s:nd} by establishing negative-definiteness of the Hessians for the first- and second-moment rate functions at their maximizers.

Thm.~\ref{t:moments} allows us to locate the exact threshold $d_\star\equiv d_\star(k)$ such that $\P(\ZZ>0)$ converges to zero for $d>d_\star$, and 
is bounded away from zero for $d>d_\star$. The following theorem improves this to a sharp threshold:
\bThm\label{t:whp}
For $k\ge k_0$,
\bnm[(a)]
\item\label{t:whp.frozen}
	$\lim_{n\to\infty}\P(\ZZ>0)=1$ for $\dlbd\le d < d_\star$; and
\item\label{t:whp.below}
	$\lim_{n\to\infty}\P(Z>0)=1$ for $d\le\dlbd$.
\enm
\eThm

Thm.~\ref{t:whp} is proved in \S\ref{s:whp} by a variance reduction argument. This issue occurs commonly in applications of the second moment method, and is often dealt with by a somewhat standard machinery known as the \emph{subgraph conditioning method} (see \cite{RW92, RW94}) which ``explains'' the variance in terms of the short cycles in the graph. Applying this method is technically demanding, and seems to us intractable in our models due to the large number of variables.

We develop instead a novel approach of taking a certain log-transform of the partition function, and bounding the incremental fluctuations of its Doob martingale with respect to the edge-revealing filtration; each increment amounts to the effect of adding a clause. We control the variance by discrete Fourier analysis applied on the spins at the boundary of a large local neighborhood of the added clause, and we show that the main contribution comes from the degree-two Fourier coefficients which correspond to the formation of short cycles in the graph.

\subsection{Notation}
For non-negative functions $f(k,d,n)$ and $g(k,d,n)$ we use any of the equivalent notations $f=O_k(g)$, $g=\Om_k(f)$, $f\lesssim_k g$, $g \gtrsim_k f$ to indicate $f\le C(k)\, g$ for a finite constant $C(k)$ depending on $k$ but not on $d,n$. (In this paper, if $f\le C(k,d)\,g$ then $f\le C(k)\, g$ simply by taking the maximum of $C(k,d)$ over the finitely many integers $d\le \dubd(k)$.) We drop the subscript $k$ to indicate when we can take the same constant $C(k)\equiv C$ for all $k\ge k_0$.

\subsection*{Acknowledgements}
We thank Amir Dembo, Elchanan Mossel, Andrea Montanari, and David Wilson for helpful conversations.

\section{Satisfying assignments}\label{s:assignments}

\subsection{Satisfiability below critical regime}

We now prove Propn.~\ref{p:sat.below} by applying the second moment method to the \acr{nae-sat} partition function.
Write $\bin_{n,p}(j)\equiv \binom{n}{j} p^j(1-p)^{n-j}$.

\bpf[Proof of Propn.~\ref{p:sat.below}]
Assume throughout that $d\le\dlbd$. By definition, $\E[Z^2]$ is the sum over pairs $\ux^1,\ux^2\in\set{\zro,\one}^V$ of the probability that both $\ux^i$ are valid \acr{nae-sat} solutions. By symmetry, the sum over $\ux^2$ is the same for all $\ux^1$; further, conditioned on $\ux^1$ being a valid solution, the probability that $\ux^2$ is also valid depends only on the number $n\al$ of vertices in which the $\ux^i$ agree. Therefore
\[
\E[Z^2]
=(\E Z)\,\sum_\al \smb{n}{n\al}\sum_\gm \mathbf{p}^\al_\gm (1-\vth)^{m\gm}
\]
where $\vth\equiv 2/(2^k-2)$, and $\mathbf{p}^\al_\gm$ is the probability, given vectors $\ux^1,\ux^2\in\set{\zro,\one}^V$ which agree in $n\al$ coordinates, that there are exactly $m\gm$ clauses $a\in F$ for which $(\ux^1\oplus\ux^2)_{\pd a}$ is not identically $\zro$ or identically $\one$. Let $D_1,\ldots,D_n$ be i.i.d.\ $\Bin(k,\al)$ random variables: then
\begin{align}
\mathbf{p}^\al_\gm \label{e:P.alpha.gamma}
&=\TS
\P(
\sum_{a=1}^m \Ind{D_a\notin\set{0,k}}=m\gm
\giv\sum_{a=1}^m D_a=m\al k
)\\ \nonumber
&\le \DS
\f{\P(\sum_{a=1}^m \Ind{D_a\notin\set{0,k}}=m\gm)}
	{\P(\sum_{a=1}^m D_a=m\al k)}
= n^{O(1)}\,\exp\{ -m\,\relent{\gm}{ \gm_0} \},
\end{align}
where $\gm_0\equiv\gm_0(\al)\equiv1-\al^k-(1-\al)^k$. Thus we conclude
\[\begin{array}{l}
\E[Z^2]
	\le n^{O(1)}\,(\E Z)\,\exp\{n\,\sup_{\al,\gm}
	\mathbf{a}(\al,\gm)\}\\
{}\qquad\qquad\text{with }\mathbf{a}(\al,\gm)
	\equiv H(\al) + (d/k)
		[ -\relent{\gm}{\gm_0}+\gm\log(1-\vth) ].
\end{array}\]
For fixed $\al$, $\mathbf{a}$ is strictly concave in $\gm$ with second derivative $-[\gm(1-\gm)]^{-1}\le -4$, and is uniquely maximized at
$\gm^\star(\al) = \gm_0(1-\vth)/(1-\vth\gm_0)$ with optimal value
\beq\label{e:Phi2}
\begin{array}{rl}
\overline{\mathbf{a}}(\al)
\hspace{-6pt}&= H(\al) +(d/k) \log(1-\gm_0\vth)\\
&= (d/k) \log(1-\vth)+ H(\al)
	+(d/k) \log( 1+ \tf{\vth}{1-\vth} [\al^k+(1-\al)^k] ).
\end{array}
\eeq
The function $\overline{\mathbf{a}}$ is symmetric in $\al$ with
$\overline{\mathbf{a}}(
\slf12)=\Phi$, 
the first-moment exponent of \eqref{e:ez}. We now show that $\al=\slf12$ is the global 
maximizer for $d\le\dlbd$. 
Since $\overline{\mathbf{a}}(\al)-\Phi$
is nondecreasing in $d$, it suffices to show this for $d=\dlbd$.
Since $(d/k) [\vth/(1-\vth)] = \log2$
for $d=\dlbd$ we find
\[ \overline{\mathbf{a}}(\al)-\Phi
\le H(\al)+ [-1+\al^k+(1-\al)^k-2/2^k]\log2
	 + O(8^{-k}).\]
It is straightforward to calculate that for $k^{-1} (\log k)^2 \le \al\le 1-k^{-1} (\log k)^2$,
\[(\overline{\mathbf{a}})''(\al)
=H''(\al) + O( k^{-(\log k)/2} ) < -3,\]
so clearly $\al=\slf12$ is the unique maximizer on this interval.
For $0\le \al\le\slf12$,
$H(\al)$ is increasing while $\al^k+(1-\al)^k$ is decreasing,
and we use this to bound
\[\TS
\sup_{k^{-3/2}\al\le k^{-4/5} } [\overline{\mathbf{a}}(\al)-\Phi]
	\le O(k^{-4/5}\log k)-k^{-1/2}\log2
	<0.
\]
For $\al\le k^{-3/2}$ we have $(1-\al)^k = 1-k\al +O(k^{1/2}\al)$, therefore
\[\TS
\sup_{2^{-3k/4}\le \al\le k^{-3/2}}
	[\overline{\mathbf{a}}(\al)-\Phi]
\le \al (3k/4)\log2
	-\al k\log2	+O(k^{1/2}\al)
	<0
\]
Lastly, recalling $H(x)+x\log c\le \log(1+c)\le c$ gives
\[\begin{array}{rl}
\sup_{0\le\al\le 2^{-3k/4}}
[\overline{\mathbf{a}}(\al)-\Phi]
\hspace{-6pt}&\le
	-(2/2^k)\log2
	+O(k^2/2^{3k/2})
	+\sup_{\al\le 2^{-3k/4}}[H(\al)-\al k\log2]\\
&\le 2^{-k}[ 1-2\log2 ]
	+O(k^2/2^{3k/2})<0.
\end{array}
\]
Therefore $\mathbf{a}$ is uniquely maximized at $(\al^\star,\gm^\star(\al^\star))=(\slf12,1-\vth)$ with maximal value $\Phi$, which proves
$\E[Z^2] \le n^{O(1)}\,(\E Z)^2$.

To remove the polynomial factor we now give a more precise calculation of the probabilities $\mathbf{p}^\al_\gm$ of \eqref{e:P.alpha.gamma}. Let $D_1,\ldots,D_m$ be i.i.d.\ $\Bin(k,\al)$ as before,
and for $0\le j\le k$ define
$p_j(\al)\equiv (1-\vth)^{\Ind{j\ne 0,k}} \bin_{k,\al}(j)$.
Then, since $\overline{\mathbf{a}}(\al)$ is uniquely maximized at $\al=\slf12$,
$$
\E[Z^2]/(\E Z)^2
= o(1)
+ 
\sum_{|\al-1/2|\le 1/3}
	\f{
	\tbinom{n}{n\al}
	\sum_{\nu}
	\Ind{\sum_j j\nu_j = k\al}
	\binom{m}{m\nu} \prod_j p_j(\al)^{m\nu_j} }
	{ 
	(\E Z) 
	e^{O( (mk)^{-1} )}
	/ \sqrt{ 2\pi mk \al(1-\al) }
	 },
$$
where the inner sum is taken over
probability measures $\nu$ on $\set{0,\ldots,k}$
such that $m\nu$ is integer-valued. By Stirling's approximation,
$$\E[Z^2]/(\E Z)^2
= o(1)
+ 
\sum_{|\al-1/2|\le 1/4}
\sum_\nu
\f{\Ind{\sum_j j\nu_j=k\al}}{\mathscr{P}(\al,\nu)}
\f{\exp\{ n\,\mathbf{b}(\al,\nu)\}}{\E Z}
$$
where $\mathbf{b}(\al,\nu)\equiv H(\al)-(d/k)\sum_j \nu_j \log [\nu_j/p_j(\al)]$ is strictly concave in $(\al,\nu)$, and the correction term $\mathscr{P}(\al,\nu)$ is $n^{O(1)}$ in general, and is $\asymp_k n^{(k+1)/2}$ for $\nu$ satisfying $\max_j 1/\nu_j\lesssim_k 1$. It is easily seen that this is indeed satisfied by $\argmax_\nu \mathbf{b}(\al,\nu)$ for $\slf14\le\al\le\slf34$, so it follows using the strict concavity of $\mathbf{b}$ that
$$\E[Z^2]/(\E Z)^2
= o(1)
+ O_k(1) 
\sum_{|\al-1/2|\le1/4}
\f{\exp\{n\,\sup_\nu\mathbf{b}(\al,\nu)\}}{n^{1/2}\,\E Z}.
$$
Of course $\sup_\nu\mathbf{b}(\al,\nu)$ need not be concave in $\al$, however, since we previously took an upper bound on $\mathbf{p}^\al_\gm$, 
$\sup_\nu\mathbf{b}(\al,\nu)\le \overline{\mathbf{a}}(\al)$
which is strictly concave near $\al=\slf12$ with global maximum $\overline{\mathbf{a}}(\slf12)=\Phi$. This proves $\E[Z^2]\lesssim_k(\E Z)^2$ for $d\le \dlbd$.
\epf

\subsection{Coarsening algorithm and frozen model}
\label{ss:coarsen.frozen}

In view of Propn.~\ref{p:sat.below} we hereafter assume unless indicated otherwise that $k\ge k_0$ large,
\beq\label{e:regime}
\dlbd\le d=(2^{k-1}-\rho)k\log 2\le \dubd
 \ (0\le\rho\le2),
\quad\text{so }\Phi
= 2^{-k} (2\rho-1)\log2+ O(4^{-k}).
\eeq
In this regime, we define the following algorithm to map a satisfying variable assignment $\ux\in\set{\zro,\one}^V$ to a \emph{coarsened} configuration $\ueta\equiv\ueta(\ux)\in\set{\zro,\one,\free}^V$. In the coarsened model, $\zro$ and $\one$ indicate variables which are ``rigid'' or ``forced'' while $\free$ indicates variables which are ``free,'' as follows:

\bdfn\label{d:forcing}
A clause--variable edge $(av)$ is said to be \emph{$\ueta$-forcing} if $\ueta_{\pd a}\in\set{\zro,\one}^k$ with $\lit_{av}\oplus\eta_v = \neg\lit_{aw}\oplus\eta_w$ for all $w\in\pd a\setminus v$. We also say that $a$ is $\ueta$-forcing and $v$ $\ueta$-forced. Recall that $\pd a=(v_1,\ldots,v_k)$ denotes the neighbors of $a$ \emph{with multiplicity}, so each clause can have at most one $\ueta$-forcing edge.\footnote{For example, if $\pd a=(v,v,v,w,\ldots,w)$ with $\eta_v\ne \eta_w$ and $\ulit_a=(\zro^k)$, the clause is not considered $\ueta$-forcing.} Given $\ueta_{\pd a}\in\set{\zro,\one}^k$, of the $2^k-2$ valid configurations of $\ulit_a$ there are exactly two which are $(av_i)$-forcing for each $1\le i\le k$, with the remaining $2^k-2-2k$ configurations not forcing to any $v_i$. A variable which is not $\ueta$-forced is said to be \emph{$\ueta$-free}.
\edfn

\noindent\textbf{Coarsening algorithm.}

{\addtolength{\leftskip}{.25in}
\noindent
Set $\ueta^0\equiv\ux$. For $t\ge0$, if there exists $v\in V$ which has $\eta^t_v\ne\free$ but which is not $\ueta^t$-forced, then take the first\footnote{First with respect to the ordering on $V=[n]$.} such $v$ and set $\eta^{t+1}_v=\free$. Set $\eta^{t+1}_w=\eta^t_w$ for all $w\ne v$.\\ Iterate until the first time $t_1$ that no such vertex $v$ remains.
\par}

\smallskip\noindent
Denote the terminal configuration $\ueta\equiv\ueta(\ux)\equiv\ueta^{t_1}$.\footnote{We could define a \emph{cluster} of \acr{nae-sat} solutions to be the pre-image of any $\ueta$ under the coarsening algorithm.} Let $Z_{\ge n\be}$ denote the contribution to $Z$ from assignments $\ux\in\set{\zro,\one}^V$ such that the coarsened configuration $\ueta(\ux)\in\{\zro,\one,\free\}^V$ has more than $n\be$ free variables.

\bppn\label{p:coarsen}
In regime \eqref{e:regime}, $\E Z_{\ge n\be}$ is exponentially small in $n$ for $\be=\slf{\cc}{2^k}$.

\bpf
By symmetry, $\E Z_{\ge n\be}= (\E Z)\,\mathbf{f}_{n\be}$ where $\mathbf{f}_{n\be}$ denotes the probability, conditioned on $\ux=\vec\zro$ being a valid \acr{nae-sat} solution, that its coarsening $\ueta$ has at least $n\be$ free variables. 

We simulate the coarsening algorithm as follows:
of the $nd$ half-edges incident to variables,
choose $e_1,\ldots,e_m$ uniformly at random (with random ordering)
to be \emph{potentially forcing}.
Edge $e_a$ corresponds to clause $a$, though here the clauses are not explicitly formed.
Conditioned on $\ux$ being a valid solution, each clause independently has probability $\vth\equiv \slf{2k}{(2^k-2)}$
 to be $\ux$-forcing (cf.~Defn.~\ref{d:forcing}):
therefore set each $e_a$ to be \emph{initially forcing} 
with probability $\vth$, independently over $a$.
Then, for each $t\ge0$,
if there exists $v\in V$ which
is incident to no (remaining) initially forcing half-edge, then
take the first such $v$ and
\bnm[(i)]
\item Delete all $d^t_v$ remaining potentially forcing half-edges incident to $v$; and
\item Delete the first $d-d^t_v$ potentially forcing half-edges among all those remaining.
\enm
The interpretation is that the coarsening algorithm sets $v$ to be a free variable at stage $t$. Thus the $d-d^t_v$ clauses incident to $v$ and potentially forcing to other variables can no longer be forcing, so we remove these clauses from consideration (step (ii)).\footnote{Step (i) does not delete any initially forcing half-edges, but step (ii) can.}

Say a variable $v$ is \emph{$t$-free} if it has no initially forcing half-edges remaining after $nt$ iterations of the above procedure. Since initially forcing edges are deleted in order, $v$ must avoid the set $E_t$ of initially forcing edges $e_a$ with index $a>ndt$. If there are $\ge nt$ free variables in the coarsened configuration $\ueta$, then the above process must survive at least $nt$ iterations. The law of $|E_t|$ is $\Bin(m-ndt,\vth)$, so (by a union bound)
\[\begin{array}{rl}
\mathbf{f}_{nt}
\hspace{-6pt}&\le\tbinom{n}{nt}
	\E[ \tbinom{ nd(1-t) }{|E_t|} / \tbinom{nd}{|E_t|} ]
	\le \tbinom{n}{nt}\,\E[(1-t)^{|E_t|}]\\
&= n^{O(1)}\,\exp\{ n[H(t) + d(\slf1k-t)\log(1-\vth t)] \}.
\end{array}
\]
If $t = \slf{C}{2^k}$
 with $C\asymp1$ then $\mathbf{f}_{nt} \le n^{O(1)}\,\exp\{ n
	(\slf{C}{2^k})[1-\log C + O(k^2 /2^k)] \}$. Then recalling \eqref{e:regime} we have $\E Z_{\ge nt} \le e^{n\Phi}\,\mathbf{f}_{nt}$ exponentially small in $n$ for $C=\cc$.
\epf
\eppn

\bdfn\label{d:frozen}
We say $\ueta\in\set{\zro,\one,\free}^V$ is a \emph{$\zof$ frozen model} configuration on $(G,\ulit)$ if
\bnm[(a)]
\item No clause $a\in F$ is \emph{unsatisfied} (meaning $\ueta_{\pd a}\in\set{\zro,\one}^k$ with $(\ulit\ueta)_a$ identically $\zro$ or $\one$);
\item Each variable $v\in V$ has $\eta_v\ne\free$ if and only if there is a clause $a\in\pd v$ with $\smash{\ueta_{\pd a}\in\set{\zro,\one}^k}$ and $\smash{\lit_{av}\oplus\eta_v=\neg\lit_{aw}\oplus\eta_w}$ for all $w\in\pd a\setminus v$ (cf.~Defn.~\ref{d:forcing}).
\enm
Some of our computations are simplified by working with the image of the $\zof$ frozen model under the projection $\set{\zro,\one}\mapsto\rig$, hereafter \emph{$\rig/\free$ frozen model}.
\edfn

Let $\ZZ_{n\be}$ denote the frozen model partition function on $(G,\ulit)$
restricted to configurations with exactly $n\be$ $\free$-vertices. In view of Propn.~\ref{p:coarsen}, in regime \eqref{e:regime}
we hereafter restrict all consideration to the \emph{truncated} $\zof$ frozen model partition function
\beq\label{e:truncated.frozen.model} \TS
\ZZ\equiv \sum_{t\le\bemax} \ZZ_{nt},\quad
\bemax\equiv \slf{\cc}{2^k}.
\eeq
We will show in \S\ref{s:sat} that restricted frozen model solutions indeed correspond to true \acr{nae}-solutions.\footnote{Some truncation is indeed necessary: the identically-$\free$ vector is a valid configuration of the unrestricted frozen model, and in fact it turns out that the dominant contribution to the partition function of the unrestricted $\zof$ frozen model comes from configurations with much higher density of free variables (roughly $\asymp (\log k)/k$) --- hence not corresponding to \acr{nae}-solutions.}

\section{First moment of frozen model}\label{s:first}

In this section we identify the leading exponential order $\bPhistar=\lim_{n\to\infty} n^{-1}\log\EZZ$ of the first moment of the (truncated) frozen model partition function \eqref{e:truncated.frozen.model}. The random $(d,k)$-regular bipartite factor graph $G\equiv\Gndk$ converges locally weakly (in the sense of \cite{MR1873300,MR2354165}) to the infinite $(d,k)$-regular tree $\tree$ --- the infinite tree with levels indexed by $\Z_{\ge0}$ such that all vertices at even integer levels are of degree $d$ (variables) and all vertices at odd integer levels are of degree $k$ (clauses). Our calculation is based on a variational principle which relates the exponent $\bPhistar$ to a certain class of Gibbs measures for the frozen model on $\tree$ which are characterized by fixed-point recursions. In fact the recursions can have multiple solutions, and much of the work goes into identifying (via \emph{a~priori} estimates) the unique fixed point which gives rise to $\bPhistar$. We begin by introducing the Gibbs measures which will be relevant for the variational principle.

\subsection{Frozen model tree recursions}\label{ss:frozen.tree}

We shall specify a Gibbs measure $\nu$ on $\tree$ by defining a consistent family of finite-dimensional distributions $\nu_t$ on the depth-$t$ subtrees $\tree(t)$. A typical manner of specifying $\nu_t$ is to specify a ``boundary law'' on the configuration on the depth-$t$ vertices, and then to define $\nu_t$ as an appropriate finite-volume Gibbs measure on $\tree(t)$ \emph{conditioned} on the boundary configuration.

In our setting some difficulty is imposed by the fact that the frozen model is not a factor model (or Markov random field) in the conventional sense that $\ueta|_A$ and $\ueta|_B$ are conditionally independent given the configuration $\ueta|_C$ on any subset $C$ separating $A$ from $B$ --- in particular, given the variable spins at level $2t$ of $\tree$, whether a variable at level $2(t-1)$ is permitted to take spin $\free$ depends on whether its neighboring $\zro$'s and $\one$'s in level $2t$ are forced by clauses in level $2t+1$.

We shall instead specify Gibbs measures for the frozen model via a message-passing system, as follows. First sample uniformly random literals $\ulit(t)$ on $\tree(t)$. Given the literals, each variable $v$ will send a message $\si_{v\to a}$ to each neighboring clause $a\in\pd v$ which represents the ``state of $v$ ignoring $a$'', and will receive in return a message $\si_{a\to v}$ representing the ``state of $a$ ignoring $v$.'' That is, $\si_{v\to a}$ will be a function $\dmp_{d-1}$ of $d-1$ incoming messages $(\si_{b\to v})_{b\in\pd v\setminus a}$, and likewise $\si_{a\to v}$ will be a function $\hmp_{av}$ (which will involve the literals at $a$) of $k-1$ incoming messages $(\si_{w\to a})_{w\in\pd a\setminus v}$. The actual state $\eta_v$ of $v$ is then a function $\dmp_d$ of all its incoming messages $\dsi_{\pd v\to v}\equiv (\si_{u\to v})_{u\in\pd v}$; the configuration may be invalidated if any variable receives conflicting incoming messages.

If on the boundary of $\tree(t)$ we are given a vector $\ueta^\uparrow\equiv(\si_{v\to w})_{v,w}$ (for $v$ at level $t$, $w$ the parent of $v$), then there is \emph{at most one} completion of $\ueta^\uparrow$ to a (bi-directional) message configuration on $\tree(t)$: iterating $\dmp_{d-1},\hmp_{av}$ gives all the messages upwards in the direction of the root, and once those are known we can recurse back down to determine the messages in the opposite direction. The measure $\nu_t$ can then be specified by giving the law of the boundary messages $\ueta^\uparrow$: our choice will be to take $\ueta^\uparrow$ i.i.d.\ according to a law $\dq$ ($t$ even) or $\hq$ ($t$ odd); consistency of the family $(\nu_t)_t$ will then amount to fixed-point relations on $\dq,\hq$.

The message-passing rules for our frozen model are as follows:
\bnm[1.]
\item \emph{Vertex message-passing rule} $\dmp_D:\set{\zro,\one,\free}^D \to\set{\zro,\one,\free}$: output
\[
\begin{array}{rl}
\free&\text{if all $D$ incoming messages are $\free$;}\\
\zro &\text{if at least one $\zro$ but no $\one$'s incoming;}\\
\one &\text{if at least one $\one$ but no $\zro$'s incoming;}\\
\text{{\sc unsat}} &\text{otherwise
	(i.e.\ both $\zro,\one$ incoming).}
\end{array}
\]

\item \emph{Clause message-passing rule} $\hmp_{av}: \usi_{\pd a\setminus v\to a}\mapsto \si_{a\to v}$: output
\[\begin{array}{rl}
	\zro & \text{if }
		\lit_{av}\oplus\zro = \neg\lit_{aw}\oplus \si_{w\to a}
		\text{ for all }w\in\pd a\setminus v; \\
	\one & \text{if }
		\lit_{av}\oplus\one = \neg\lit_{aw}\oplus \si_{w\to a}
		\text{ for all }w\in\pd a\setminus v; \\
	\free & \text{otherwise.}
	\end{array}\]
\enm
We then define
\beq\label{e:gibbs.frozen}
\TS
Z_t\,\nu_t(\ulit(t),\ueta(t),\ueta^\uparrow)
=\begin{cases}
\dbq(\ueta^\uparrow)\equiv\prod_i \dq(\eta^\uparrow_i),
	& \text{$t$ even},\\
\hbq(\ueta^\uparrow)\equiv\prod_i \hq(\eta^\uparrow_i),
	& \text{$t$ odd},\\
\end{cases}
\eeq
\emph{provided} $\ueta^\uparrow$ completes
leads to a valid message configuration
(no \acr{unsat} messages) on $\tree(t)$
with respect to literals $\ulit(t)$. The root marginal is then given by
$$
\nu_1(\eta_\rt=x)
= \f{(\hqf+\hq_x)^d-(\hqf)^d}
	{ (\hqf+\hqz)^d+(\hqf+\hqo)^d-(\hqf)^d }
\quad
\text{for }x=\zro\text{ or }\one,
$$
with the remaining probability going to $\eta_\rt=\free$.
The $\nu_t$ are consistent if and only if $q\equiv(\dq,\hq)$ satisfies the \emph{frozen model recursions}
\beq\label{e:frozen.recursions}
\hq_\zro=\hq_\one
	=(\slf{2}{2^k})
	( \dqz+\dqo )^{k-1},\quad
\dqz
=
\f{(\hqf+\hqz)^{d-1}-(\hqf)^{d-1}}
	{(\hqf+\hqz)^{d-1}+(\hqf+\hqo)^{d-1}-(\hqf)^{d-1}}
=\dqo
\eeq
with $\hqf=1-\hqz-\hqo$
and $\dqf=1-\dqz-\dqo$.

\blem\label{l:contract}
In the regime $\dqf\lesssim 2^{-k}$, 
the recursion \eqref{e:frozen.recursions} has a unique solution
$q^\star$, which furthermore satisfies
$2^k (\dqf)^\star = \slf12 + O(k^2/2^k)$.

\bpf
Writing $\qrig\equiv1-\dqf$ and $v\equiv \hqf/(\hqz+\hqf)$, we see that a solution of \eqref{e:frozen.recursions} corresponds to a solution of
the equations
\beq\label{e:q.rec}
\qrig=\qrig_{d-1}(v)\equiv \f{2-2v^{d-1}}{2-v^{d-1}},\quad
v=v_{k-1}(\qrig) \equiv \f{1-2(\qrig/2)^{k-1}} {1-(\qrig/2)^{k-1}}.
\eeq
If $1-\qrig\lesssim 2^{-k}$ then $v_{k-1}(\qrig)=1-\slf{2}{2^k} + O(\slf{k}{4^k})$, therefore $v_{k-1}(\qrig)^{d-1} = 2^{-k} + O(\slf{k^2}{4^k})$ and $\qrig_{d-1}\circ v_{k-1}(\qrig) = 1- 2^{-k-1} + O(\slf{k^2}{4^k})$. In this regime we also calculate
$$v_{k-1}'(\qrig)
=-\f{(k-1) (\qrig/2)^{k-1}}{\qrig[1-(\qrig/2)^{k-1}]^2}
\asymp k 2^{-k},\quad
\qrig_{d-1}'(v)
=
-\f{2(d-1)v^{d-1}}{v(2-v^{d-1})^2}
\asymp k,$$
thus $(\qrig_{d-1}\circ v_{k-1})'\asymp\slf{k^2}{2^k}$ so in this regime \eqref{e:q.rec} must have the unique solution as claimed.
\epf
\elem

\brmk
Note that if $\nu$ is the Gibbs measure on $\tree$ corresponding to a solution $q^\star$ of \eqref{e:frozen.recursions}, then $\nu(\si_\rt\ne\free)$ is a fixed point of $\qrig_d\circ v_{k-1}$. In the regime of Lem.~\ref{l:contract} the fixed points of $\qrig_d\circ v_{k-1}$ and $\qrig_{d-1}\circ v_{k-1}$ are nearly identical, so in view of Propn.~\ref{p:coarsen} we are justified in restricting attention to fixed points with $\dqf\lesssim 2^{-k}$.
\ermk

\subsection{Auxiliary model}

On the tree $\tree$, the frozen model configuration $\ueta$ can be uniquely recovered from the configuration $\usi$ of messages on all the directed edges: each vertex spin $\eta_v$ is determined by applying $\dmp_d$ to the incoming messages. We refer to $\usi$ as the \emph{auxiliary} configuration, and we now observe that we can define a model on auxiliary configurations on $(d,k)$-regular bipartite graphs which is in bijection with the frozen model but has the advantage of being a \emph{factor model} in a relatively simple sense.

The spins of the auxiliary model on the bipartite factor graph are the bidirectional messages $\si_{va}\equiv\si_{av}\equiv(\si_{v\to a},\si_{a\to v})$, taking values in the alphabet $\msg\equiv\set{\zro,\one,\free}^2\setminus\set{\zo,\oz}$. Write $\dsi_v$ for the $d$-tuple of spins on the edges incident to variable $v\in V$, and write $\ksi_a$ for the pair of spins on the edges incident to clause $a\in F$.

In the \emph{auxiliary model},
each configuration $\usi\in\msg^E$
receives the factor model weight
\beq\label{e:aux}
\Psi(\usi)
\equiv
\Psi_{G,\ulit}(\usi)
\equiv
\prod_{v\in V}
\dpsi(\dsi_v)
\prod_{a\in F}
\hpsi^a(\ksi_a)\eeq
where the variable factor weight $\dpsi(\dsi_v)$  is simply the indicator that each outgoing message $\si_{v\to a}$ is determined by the message-passing rule $\dmp_{d-1}$ from the incoming messages $\si_{b\to v}$, $b\in\pd v\setminus a$; and likewise the clause factor weight $\hpsi^a(\ksi_a)$ is the indicator that each outgoing message $\si_{a\to v}$ is determined by the message-passing rule $\hmp_{av}$ from the incoming messages $\si_{w\to a}$, $w\in\pd a\setminus v$. Then, with $\neg\free\equiv\free$,
we have $\hpsi^a(\ksi_a)=\hpsi^\circ(\ksi_a\oplus \ulit_a)$
where $\dpsi$ and $\hpsi^\circ$ are given explicitly by
\beq\label{e:zof}
\dpsi(\dsi)\equiv\text{\footnotesize$
\begin{cases}
1,& \dsi=(\ff^d),\\
1,& \dsi\in\perm[(\fz,\zf^{d-1})],\\
1,& \dsi\in\perm[(\fo,\of^{d-1})],\\
1,& \dsi\in\perm[(\zz^j,\zf^{d-j})_{j\ge2}],\\
1,& \dsi\in\perm[(\oo^j,\of^{d-j})_{j\ge2}],\\
0, & \text{else;}
\end{cases}$},\quad\hpsi^\circ(\ksi)
\equiv\text{\footnotesize$\begin{cases}
1, & \ksi\in\perm[ (\zz \text{ or } \fz,\of^{k-1}) ]\\
1, & \ksi\in\perm[ (\oo \text{ or } \fo,\zf^{k-1}) ]\\
1, &
\ksi\in\perm[(\zro\free^j,\one\free^{k-j})_{2\le j\le k-2}],\\
1, & \ksi\in\perm[(\ff,\zro\free^j,\one\free^{k-1-j})_{1\le j\le k-2}],\\
1, & \ksi\in\perm[(\ff^j,\zro\free^{s},\one\free^{k-1-s})_{j\ge 2}],\\
0, & \text{else;}
\end{cases}$}
\eeq
with $\mathrm{Per}(\usi)$ the set of permutations of $\usi$. 
We refer to this as the \emph{factor model} with \emph{specification}
$\phi\equiv(\dpsi,\hpsi)$.

\brmk\label{r:frozen.aux.bij}
The frozen model is in exact bijection with the auxiliary model. Given an auxiliary configuration $\usi\in\msg^E$, the corresponding frozen configuration $\ueta$ is given by coordinate-wise application of $\dmp_d$. The inverse mapping $\ueta\mapsto\usi$ can be defined as follows:
first determine the clause-to-variable messages by setting $\si_{a\to v}$ to be $\eta_v$ if $(av)$ is $\ueta$-forcing and $\free$ otherwise, equivalently $\si_{a\to v}=\hmp_{av}(\eta_{\pd a\setminus v})$. Then determine the variable-to-clause messages $\si_{v\to a}$ by applying $\dmp_{d-1}$ (since we assumed $\ueta$ is a valid frozen model configuration, $v$ cannot receive conflicting incoming messages $\si_{a\to v}=\zro$ and $\si_{b\to v}=\one$).
\ermk

\bdfn
The \emph{$\zof$ auxiliary model on $G$} is defined to be the average of the auxiliary model \eqref{e:aux} over all literal configurations $\ulit$. The \emph{$\rig/\free$ auxiliary model} is the image of the $\zof$ auxiliary model under the projection $\Pi:\set{\zro,\one}\mapsto\rig,\free\mapsto\free$.
\edfn

It is easily seen that the $\zof$ auxiliary model is again a factor model on $G$, with variable factor $\dpsi$ as before and 
clause factor $\hpsi(\ksi)\equiv 2^{-k}\sum_{\ulit}\hpsi^\circ(\ksi\oplus\ulit)$. Further, $\dpsi(\dsi)$ and $\hpsi(\ksi)$ depend on $\dsi$ and $\ksi$ only through their projections under $\Pi$, so we conclude that $\rig/\free$ auxiliary model on $G$ is a factor model with specification
\beq\label{e:rf}
\text{\footnotesize
$\dpsi(\dsi)
= \begin{cases}
1,& \dsi=(\ff^d),\\
2,& \dsi\in\perm(\fr,\rf^{d-1}),\\
2,& \dsi\in\perm[(\rr^j,\rf^{d-j})_{j\ge2}],\\
0, & \text{else;}
\end{cases}\quad
\hpsi(\ksi)
\equiv
2^{-k}\begin{cases}
2,
& \ksi\in\perm[(\rr\text{ or }\fr,\rf^{k-1})],\\
2^k-2-2k,
& \ksi = (\rf^k),\\
2^k-4,
& \ksi\in\perm(\ff,\rf^{k-1}),\\
2^k,
& \ksi\in\perm[(\ff^j,\rf^{k-j})_{j\ge2}],\\
0,
& \text{else.}
\end{cases}$}
\eeq

\subsection{Bethe variational principle}\label{ss:bethe}

The primary purpose of defining the auxiliary model is that it gives us the following approach for calculating $\EZZ$. Given an auxiliary configuration $\usi$, consider the normalized empirical measures
$$\begin{array}{llr}
\dbh(\dsi)
	\equiv
	n^{-1}\sum_{v\in V} \Ind{\dsi_v=\dsi}
	& (\dsi\in\msg^d)
	\quad
	& \text{variable empirical measure;}\\
\hbh(\ksi)
	\equiv
	m^{-1}
	\sum_{a\in F}\Ind{\ksi_a=\ksi}
	& (\ksi\in\msg^k)
	\quad
	& \text{clause empirical measure.}
\end{array}$$
We regard $\bh\equiv(\dbh,\hbh)$ as a vector indexed by
$\supp\vph\equiv(\supp\dpsi,\supp\hpsi)$. For $\si\in\msg$ and $\dsi\in\supp\dpsi$ let $\dH_{\si,\dsi}$ denote the number of appearances of $\si$ in $\dsi$, and similarly write $\hH_{\si,\ksi}$ for the number of appearances of $\si$ in $\ksi$. For $\bh$ to correspond to a valid configuration $\usi$,
the variable and clause empirical measures must give rise to the same edge marginals
$$\TS
\vh = d^{-1} \dH\dbh
	= k^{-1} \hH\hbh,\quad
\vh(\si)\equiv (nd)^{-1}\sum_{(va)\in H}\Ind{\si_{va}=\si}.
$$

\bdfn\label{d:simplex}
Given $\vph\equiv(\dpsi,\hpsi)$ let $\simplex$ denote the space of probability measures $\bh\equiv(\dbh,\hbh)$ on
$\supp\vph$
(that is, $\dbh$ is a probability measure on $\supp\dpsi$ while $\hbh$ is a probability measure on $\supp\hpsi$) such that
\bnm[(i)]
\item $(\dbh,\tf{d}{k}\hbh)$ lies in the kernel of matrix
	$H_{\simplex}\equiv \bpm \dH & -\hH \epm$, and
\item $\dbh(\dmp_d(\dsi)=\free)\le\bemax$ (cf.~\eqref{e:truncated.frozen.model}).
\enm
Let $\dbs\equiv|\supp\dpsi|$, $\hbs\equiv|\supp\hpsi|$, and $\bar{s}\equiv|\supp\hpsi|=|\msg|$: we shall show  (Lem.~\ref{l:lin.bij}) that $H_{\simplex}$ is surjective, therefore $\simplex$ is an $(\dbs+\hbs-\bar{s}-1)$-dimensional space.
\edfn

The expected number of auxiliary configurations on $\Gndk$ with empirical measure $\bh$ is
$$\EZZ(\bh)
=
\f{\binom{n}{n\dbh}
\binom{m}{m\hbh}}
	{ \binom{nd}{ nd\vh } }
\dpsi^{n\dbh}
\hpsi^{m\hbh}
\equiv
\f{n! m! \prod_\si ( nd\vh(\si) )!}{(nd)! }
\prod_{\dsi}
	\DSf{\dpsi(\dsi)^{n\dbh(\dsi)} }{(n\dbh(\dsi))!}
\prod_{\ksi}
	\DSf{\hpsi(\ksi)^{m\hbh(\ksi)} }{(m\hbh(\ksi))!}.
$$
Stirling's formula gives
$\EZZ(\bh)= n^{O(1)}\, \exp\{n\bPhi(\bh)\}$ where
\beq\label{e:bethe}
\bPhi(\bh)
\equiv\sum_{\dsi}\dbh(\dsi)\log\DSf{\dpsi(\dsi)}{\dbh(\dsi)}
+\tf{d}{k}\sum_{\ksi}\hbh(\ksi)\log\DSf{\hpsi(\ksi)}{\hbh(\ksi)}
-d\sum_\si\vh(\si)\log\DSf{1}{\vh(\si)}.\eeq
If further $\min\bh\gtrsim_k 1$ as $n\to\infty$, then
\beq\label{e:poly.correction}
\EZZ(\bh)
= \f{e^{O_k(n^{-1})}}{ (2\pi n)^{( \dbs+\hbs-\bar{s}-1 )/2} }
	\underbrace{
	\Big[\f{
\prod_\si d\vh(\si)}
	{k \prod_{\dsi} {\DS\dbh(\dsi)}
	\prod_{\ksi}
		{\DS 
		(\tf{d}{k}\hbh(\ksi)) }
		}
		\Big]^{1/2}
		}_{\mathscr{P}(\bh)}
	\exp\{n\,\bPhi(\bh)\}
\eeq
The first moment of frozen model configurations is $\EZZ=\sum_{\bh\in\simplex}\EZZ(\bh)$. The aim of this section is to compute the exponent
$\bPhistar=\lim_n n^{-1}\log\EZZ$
by determining the maximizer $\bhstar\equiv(\dbhstar,\hbhstar)$ of $\bPhi$ on $\simplex$. Observe it is clear from the functional form of $\bPhi$ that $\dbhstar$ and $\hbhstar$ must be \emph{symmetric} functions on $\msg^d$ and $\msg^k$ respectively.

If $\bhstar$ lies in the interior $\simplexint$ of $\simplex$ then it must be a stationary point for $\bPhi$. Such points correspond to a generalization of the tree Gibbs measures considered in \S\ref{ss:frozen.tree}, where the boundary conditions are specified by a law on incoming \emph{and outgoing} messages, as follows: first sample uniformly random literals $\ulit(t)$ on $\tree(t)$ as before. If $\usi(t)$ is a message configuration on the edges of $\tree(t)$ --- including the edges $E(t-1,t)$ joining levels $t-1$ and $t$ --- then let $\Psi_t(\ulit(t),\usi(t))$ denote the product of the factor weights $\dpsi(\dsi_v)$, $\hpsi^a(\ksi_a)$ over all $v,a\in \tree(t-1)$. For probability measures $\hd,\hh$ on $\msg$ we define the measures
\beq\label{e:gibbs.aux}
Z_t\,\nuaux_t(\ulit(t),\usi(t))
= \begin{cases}
\Psi_t(\ulit(t),\usi(t))
	\prod_{e\in E(t-1,t)} \hd_{\si_e},
	& \text{$t$ even,}\\
\Psi_t(\ulit(t),\usi(t))
	\prod_{e\in E(t-1,t)} \hh_{\si_e},
	& \text{$t$ odd,}
\end{cases}
\eeq
with $Z_t$ the normalizing constant which makes $\nuaux_t$ a probability measure. This generalizes the definition of $\nu_t$ in
\eqref{e:gibbs.frozen} by taking $\hd_{\eta\eta'}$ proportional to $\dq_\eta$ and $\hh_{\eta\eta'}$ proportional to $\hq_{\eta'}$, i.e.
\beq\label{e:bp.sym}
\hd_{\eta\eta'} = \dq_\eta / (2+\dqf),\quad
\hh_{\eta\eta'} = \hq_{\eta'}/(2+\hq_\free).
\eeq
The family $(\nuaux_t)_t$ is consistent if and only if $h\equiv(\hd,\hh)$ satisfies the \emph{Bethe recursions}
\beq\label{e:bp}\TS
\bpdzh\,\hd_\si
	= 
	\sum_{\dsi\,:\,\si_1=\si}
	\dpsi(\dsi) \prod_{i=2}^d \hh_{\si_i},\quad
\bphzh\,\hh_\si
=\sum_{\ksi\,:\,\si_1=\si}\hpsi(\ksi) \prod_{i=2}^k \hd_{\si_i}
\eeq
(with $\bpdzh,\bphzh$ the normalizing constants); these generalize the frozen model recursions \eqref{e:frozen.recursions}, as we shall see explicitly below. Thus a solution $h$ of \eqref{e:bp} specifies a Gibbs measure $\nuaux$ for the auxiliary model on $\tree$ which generalizes the measures $\nu$ described in \S\ref{ss:frozen.tree}.

It is clear from the $\zro/\one$ symmetries of $\hpsi$ that any solution $h$ of the $\zof$ Bethe recursions must also have the $\zro/\one$ symmetry, and as a consequence must correspond to a solution $g$ of the $\rig/\free$ Bethe recursions via
\beq\label{e:bp.zof.rf}
\hh_{\eta\eta'}
	= \gh_{\Pi\eta,\Pi\eta'}/(2-\gh_\ff),\quad
\hd_{\eta\eta'} =
	2^{\Ind{\eta=\eta'=\free}}
	\gd_{\Pi\eta,\Pi\eta'}/2.
\eeq
The $\rig/\free$ Bethe recursions read explicitly as follows:
\[\begin{array}{l}
\bphzg\,\gh_\rr
	=\bphzg\,\gh_\fr=(\slf{2}{2^k}) (\gd_\rf)^{k-1},\\
\bphzg\,\gh_\ff
	=(\gd_\rf+\gd_\ff)^{k-1}
		-(\slf{4}{2^k}) (\gd_\rf)^{k-1},\\
\bphzg\,\gh_\rf
=(\gd_\rf+\gd_\ff)^{k-1}
	-(\slf{2}{2^k}) (k+1) (\gd_\rf)^{k-1}\\
\qquad\qquad{}+(\slf{2}{2^k}) (k-1) (\gd_\rf)^{k-2}( \gd_\rr+\gd_\fr-2\gd_\ff ),\\
\bpdzg\,\gd_\ff= (\gh_\ff)^{d-1},
	\qquad
	\bpdzg\,\gd_\fr=2(\gh_\rf)^{d-1},\\
\bpdzg\,\gd_\rr
	=\bpdzg\,\gd_\rf
	=2[ (\gh_\rr+\gh_\rf)^{d-1}-(\gh_\rf)^{d-1} ],
\end{array}\]
where $\gd_\rf$ was simplified using $\gh_\rr=\gh_\fr$. The recursion for $\gh_\rf$ then simplifies to
$$\bphzg\,\gh_\rf=\bphzg\,\gh_\ff
	+(\slf{2}{2^k})(k-1)(\gd_\rf)^{k-2}(\gd_\rf-2\gd_\ff),$$
so we see that $\gd_\fr=2\gd_\ff$ if and only if $\gh_\rf=\gh_\ff$,
in which case the corresponding solution $h$ of the $\zof$ Bethe recursions satisfies the symmetries \eqref{e:bp.sym}. A fixed point of the recursion \eqref{e:q.rec} is given by $\qrig=\gd_\rf/(\gd_\rf+\gd_\ff)$ and $v=\gh_\rf/(\gh_\rr+\gh_\rf)
	=\hh_\ff/(\hh_\oo+\hh_\ff)
	=4\hh_\ff/(1+\hh_\ff)$,
using the relation $4\hh_\oo+3\hh_\ff=1$. In the reverse direction, any solution $\qrig,v$ of \eqref{e:q.rec} gives rise to a Bethe solution via
\beq\label{e:grf.q}
	\begin{array}{l}
	\gh_\rf=\gh_\ff
		=2\hh_\ff/(1+\hh_\ff)
		=v/2, \qquad
	\gh_\rr=\gh_\fr
		=2\hh_\oo/(1+\hh_\ff)
		=(1-v)/2, \\
	\gd_\fr
		=2\gd_\ff=2\hd_\ff
		=2 \dqf/(2+\dqf), \qquad
	\gd_\rr=\gd_\rf=2\hd_\of=
		\qrig/(2+\dqf).
	\end{array}
\eeq
This proves our claim that the measures $\nuaux$ generalize the measures $\nu$ of \S\ref{ss:frozen.tree}.

The connection between these Gibbs measures and the rate function $\bPhi$ is given by the following variational principle:

\blem\label{l:interior.bp}
If $\phi\equiv(\dpsi,\hpsi)$ is such that both $\dH$ and $\hH$ are surjective, then any stationary point $\bh$ of $\bPhi$ belonging to $\simplexint$ corresponds to a Bethe fixed point solving \eqref{e:bp} via
\beq\label{e:bij}
\TS
\dbz_h\,\dbh(\dsi)=\dpsi(\dsi)\prod_{i=1}^d \hh_{\si_i},\quad
\hbz_h\,\hbh(\ksi)=\hpsi(\ksi)\prod_{i=1}^k\hd_{\si_i},\quad
\barz_h\,\vh(\si)= \hd_\si \hh_\si
\eeq
with $\dbz_h,\hbz_h,\barz_h$ normalizing constants satisfying $\barz_h=\dbz_h/\bpdzh=\hbz_h/\bphzh$ for $\bpdzh,\bphzh$ as in \eqref{e:bp}.

\bpf
At an interior stationary point $\bh$, consider differentiating $\bPhi$ in direction $\bde\equiv(\dbde,0)$ with $\dH\dbde=0$, so that $\bh+s\bde\in\simplexint$ for $|s|$ small. Writing
$\dba\equiv\log[ \dpsi(\dsi)/\dbh(\dsi) ]$,
\[\begin{array}{l}
0=\pd_s\bPhi(\bh+s\bde)|_{s=0}
	=\sum_{\dsi}\dbde(\dsi)\dba(\dsi)
	=\sum_{\dsi}\dbde(\dsi)\dbep(\dsi),\\
\qquad\qquad
\dbep(\dsi)
	\equiv\dba(\dsi)+\sum_{i=1}^d\dlm(\si_i)
	\text{ with }
	\dlm:\msg\to\R\text{ arbitrary.}
\end{array}
\]
We claim it is possible to choose $\dlm$ such that $\dbep$ has marginals $\bar\ve\equiv0$: in vector notation $\dbep=\dba+\dH^t\dlm$, so this amounts to solving $\dH\dba+\dH\dH^t\dlm=0$,
which has a unique solution $\dlm$ by surjectivity of $\dH$.
Taking $\dbde=\dbep$ with this value of $\dlm$ in the above derivative gives
\beq\label{e:h.prod}
\TS
\dbh(\dsi)
=\dpsi(\dsi) \prod_{i=1}^d e^{\dlm(\si_i)},\quad
\text{likewise}\quad
\hbh(\ksi)=\hpsi(\ksi)\prod_{i=1}^k e^{\hlm(\si_i)}.
\eeq
Now differentiate in the direction of general $\bde$
with $\tf1d\dH\dbde=\vde=\tf1k\hH\hbde$, so $\bh+s\bde\in\simplex$ for small $|s|$. Applying \eqref{e:h.prod} and simplifying gives
$$\TS
0=\pd_s\bPhi(\bh+s\bde)|_{s=0}
=d\sum_\si \vde(\si)
\bar\rho(\si),\quad
\bar\rho(\si)\equiv \log\vh(\si)-\dlm(\si)-\hlm(\si) .$$
By surjectivity we may choose $\bde$ with
$\vde(\si)=\bar\rho(\si)-|\msg|^{-1}\sum_{\si'}\bar\rho(\si')$, and then substituting into the above we find that $\log\vh-\dlm-\hlm$ is a constant function of $\si$, that is,
\[
\vh(\si)\text{ equals }
	e^{\dlm(\si)}\,e^{\hlm(\si) }
	\text{ up to normalizing constant.}
\]
On the other hand, the marginal of \eqref{e:h.prod} reads
\[\TS
\vh(\si)
= e^{\dlm(\si)} \sum_{\dsi:\si_1=\si}
	\dpsi(\dsi) \prod_{i=2}^d e^{\dlm(\si_i)}
= e^{\hlm(\si)} \sum_{\ksi:\si_1=\si}
	\hpsi(\ksi) \prod_{i=2}^k e^{\hlm(\si_i)}.\]
Comparing the expressions for $\vh(\si)$ shows that the probability measures $\hd$ and $\hh$ on $\msg$ obtained by normalizing respectively $\smash{e^{\hlm(\si)}}$ and $\smash{e^{\dlm(\si)}}$ must solve the Bethe recursions \eqref{e:bp}. Lastly \eqref{e:h.prod} shows that $\bh$ corresponds to $h\equiv(\hd,\hh)$ via \eqref{e:bij}, concluding the proof.
\epf
\elem

\bthm\label{t:first.moment.exponent}
In the $\zof$ auxiliary model, let $\bhstar$ denote the unique stationary point of $\bPhi$ which corresponds --- via \eqref{e:bij} and \eqref{e:bp.sym} --- to the solution $q^\star$ of the frozen model recursions \eqref{e:frozen.recursions} which was identified in Lem.~\ref{l:contract}. The unique maximizer of $\bPhi$ on $\simplex$ is given by $\bhstar$.
\ethm

In view of Lem.~\ref{l:interior.bp} and our preceding discussion of Gibbs measures, Thm.~\ref{t:first.moment.exponent} will follow by showing
\bnm[1.]
\item Any global maximizer $\bh$ of $\bPhi$ on $\simplex$ must lie in the interior $\simplexint$, and so corresponds via \eqref{e:bij}
to a solution $h$ of the Bethe recursions \eqref{e:bp}. (For the required surjectivity of $\dH,\hH$ see Lem.~\ref{l:lin.bij}.)
\item Any such Bethe solution $h$ satisfies the symmetries \eqref{e:bp}, therefore reduces to a solution $q$ of the frozen model recursions \eqref{e:frozen.recursions}. Further $q$ is in the regime of Lem.~\ref{l:contract}, which uniquely identifies $\bh=\bhstar$.
\enm

\subsection{Boundary maximizers}

In this section we verify (by {\it a priori} estimates) that $\bPhi$ has no maximizers on the boundary of $\simplex$. By Rmk.~\ref{r:frozen.aux.bij} we may work interchangeably with the frozen and auxiliary models.

We begin with a preliminary calculation. For a vector $\vec\ell\in\Z^n$ let $\grigid[\vec\ell](n,E)$ denote the probability, with respect to a uniformly random assignment of $E$ forcing half-edges to $n$ degree-$d$ variables, that variable $i$ receives at least $\ell_i$ of the $E$ edges for each $1\le i\le n$. If $\vec\ell$ is the constant vector $(\minval,\ldots,\minval)$ we write $\grigid[\vec\ell](n,E)\equiv\grigid[\minval](n,E)$.

\blem\label{l:forcing}
For $\target=y(\log d)/d$ with $y\asymp1$ and $\vec\ell$ upper bounded by $\minval\lesssim1$ (uniformly in $d$),
\[\grigid[\vec\ell](n,nd\target)
\asymp
\exp\{
	O( nd^{-2y}(\log d)^{2\minval-1} )\}
\prod_{i=1}^n \P_{\target}(X_i\ge\ell_i).
\]

\bpf
Let $X_1,X_1,\ldots,X_n$ be i.i.d.\ $\Bin(d,\origin)$ random variables, with joint law $\P_\origin$: then
\[ \f{\grigid[\vec\ell](n,nd\target)}
	{\prod_{i=1}^n \P_{\target}(X_i\ge\ell_i)}
= \f{\P_\origin( \sum_{i=1}^n X_i=nd\target \giv X_i\ge\ell_i \ \forall i)}
	{ \P_\origin( \sum_{i=1}^n X_i=nd\target) }.\]
For any $\ell$ the conditional mean $\E_\origin[X\giv X\ge\ell]$ is increasing in $\origin$ (the derivative is the variance of a certain random variable), thus there is a unique value $\origin=[1+O(d^{-y}(\log d)^{\minval-1})]\,\target$ such that $\E_\origin[ \sum_{i=1}^n X_i \giv X_i\ge\ell_i \ \forall i ]=nd\target$. For this value of $\origin$, the local {\sc clt} (see \cite{MR1340834}) combined with Stirling's approximation gives
\[ \f{\grigid[\vec\ell](n,nd\target)}
	{\prod_{i=1}^n \P_{\target}(X_i\ge\ell_i)}
	\asymp \exp\{ nd\,\relent{\target}{\origin} \}
	=\exp\{ O(n d^{-2y} (\log d)^{2\minval-1} ) \},\]
concluding the proof.
\epf
\elem

\blem\label{l:large.pf}
For $k\ge k_0$ and $\dlbd\le d\le \dubd$,
the contribution to $\EZZ$ (see~\eqref{e:truncated.frozen.model}) from 
all $\be\le\bemax$ with $|2^{k+1}\be-1|\ge 2^{-k/8}$ is exponentially small in $n$
compared with $\EZZ$. Further
\beq\label{e:frozen.lbd}
e^{ O(n/2^k) }
= \EZZ
\ge \EZZ_{1/2^{k+1}}
\ge \exp\{
	n[ \Phi-\tf{1}{2^{k+1}}+O(\tf{k^{O(1)}}{2^{4k/3}}) ]\},\quad
\text{with $\Phi$ as in \eqref{e:ez}.}
\eeq

\bpf
Recall that $\ZZ_{n\be}$ denotes the contribution to the frozen model partition function from configurations with $n\be$ free variables. 

\smallskip\noindent\emph{Upper bound ignoring forcing constraints.}\\
Let $\mathbf{Y}_{n\be}$ denote the partition function of $\zof$ frozen configurations with $n\be$ frees where we \emph{ignore} the requirement that rigid variables be forced, so clearly $\mathbf{Y}_{n\be}\ge\ZZ_{n\be}$. In a given frozen configuration let $m\nu_j$ ($0\le j\le k$) count the number of clauses incident to exactly $j$ free variables; and let $\mathbf{p}^\be_\nu$ denote the probability of empirical measure $\nu$ of clauses with respect to a uniformly random matching between clause half-edges and variable half-edges with density $\be$ of frees. Then
\[\TS
	\EZZ_{n\be} \le \E\mathbf{Y}_{n\be}
	=2^{n(1-\be)}\tbinom{n}{n\be}
	\sum_\nu\mathbf{p}^\be_\nu(1-\slf{2}{2^k})^{m\nu_0}
	(1-\slf{4}{2^k})^{m\nu_1}.\]
Similarly to the calculation in the proof of Propn.~\ref{p:sat.below}, let $D_1,\ldots,D_m\sim\Bin(k,\be)$, and calculate $\mathbf{p}^\be_\nu = \P( \sum_a\Ind{D_a=j}=m\nu_j\text{ for all }0\le j\le k \giv\sum_a D_a=mk\be )$: since the local {\sc clt} implies $\P(\sum_a D_a=mk\be)=n^{O(1)}$, we find
\[\begin{array}{rl}
\E\mathbf{Y}_{n\be}
	\hspace{-6pt}&=\TS
n^{O(1)}\,
2^{n(1-\be)}\tbinom{n}{n\be}
\sum_\nu
\Ind{\sum_j j\nu_j=k\be}
\tbinom{m}{m\nu} \prod_j p_j^{m\nu_j}\quad\text{where}\\
&
p_0\equiv(1-\slf{2}{2^k})\,\bin_{k,\be}(0), \
p_1\equiv(1-\slf{4}{2^k})\,\bin_{k,\be}(1), \
p_j\equiv\bin_{k,\be}(j) \text{ for } 2\le j\le k.
\end{array}\]
The above is optimized at $\nu_j=p_j u^j/c$ where $c\equiv\sum_j p_j u^j$ and $u$ is chosen such that $k\be$ matches $\smash{\sum_j j\nu_j = (\sum_j jp_j u^j)/(\sum_j p_j u^j)}$. The latter is increasing in $u$, and it is straightforward to check that it has a unique solution $u=1+\slf{2}{2^k}+O(\slf{k}{4^k})$. This implies $c = 1-\slf{2}{2^k}+O(\slf{k^2}{8^k})$, thus
\[\begin{array}{rl} \E\mathbf{Y}_{n\be}
	\hspace{-6pt}&= n^{O(1)}\,\exp\{ n\,\mathbf{y}(\be) \}
\text{ with } \\
\mathbf{y}(\be)
	\hspace{-6pt}&\equiv(1-\be)\log2+H(\be)
	+(\slf{d}{k}) [ \log c-k\be\log u ]
	\\
	&=
	-\log2 + (\slf{d}{k})\log(1-\slf{2}{2^k})
	-\be\log2
	+H(\be)
	-d\be (\slf{2}{2^k})
	+O(\slf{k^2}{4^k})\\
&=\Phi+\be[\log(e/\be)
	-\log 2^{k+1}]
	+O(\slf{k^2}{4^k}),
\end{array}
\]
with $\Phi$ as in \eqref{e:ez}
(not depending on $\be$).

\smallskip\noindent\emph{Bounds with forcing constraints.}\\
Suppose we condition on an assignment of edges such that every clause is satisfied, and no $\free$-variables are illegally forced. Each of the $m\nu_0$ fully rigid clauses is forcing with probability $\vth\equiv \slf{2k}{(2^k-2)}$, and $m\nu_0$ is clearly sandwiched between $m$ and $\acute{m}\equiv m(1-k\be)$, therefore
\beq\label{e:Z.beta.lbd.ubd}
\EZZ_{n\be}
\begin{cases}
\le
	\E\mathbf{Y}_{n\be}
	\sum_\al \bin_{m,\vth}(m\al)
	\,\grigid(n(1-\be),m\al);\\
\ge
\E\mathbf{Y}_{n\be}
\sum_\al
\bin_{\acute{m},\vth}(\acute{m}\al)\,\grigid( n(1-\be),\acute{m}\al).
\end{cases}\eeq
For $\al=[1+O(2^{-k/3})]\,\vth$ we have
\[
\smf{m\al}{n(1-\be)d}
	= [1+O(2^{-k/3})]\,\slf{2}{2^k}
	=\smf{y\log d}{d}
	\quad\text{with }y
	= 1-\smf{\log k}{k\log2}+O(\slf1k);
\]
the same estimate holds with
$\acute{m},\acute{y}$ in place of $m,y$. Applying Lem.~\ref{l:forcing} then gives
\[ \begin{array}{rl}
\grigid(n(1-\be),m\al)
	\hspace{-6pt}&=[ 1
		-(1- \slf{2}{2^k} )^d\exp\{ O(\slf{k}{2^{k/3}}) \}
		]^{n(1-\be)}
	\exp\{ O(\slf{nk^{O(1)}}{4^k}) \}\\
	&=\exp\{ \slf{n}{2^k}[-1 + O(\slf{k}{2^{k/3}})] \}
	= \grigid(n(1-\be),\acute{m}\al).
	\end{array}
\]
For $|1-\al/\vth|\ge 2^{-k/3}$ we have $\bin_{m,\vth}(m\al) \le \exp\{ -\slf{n}{2^{k/2}} \}$: consequently, in each of the two sums on the right-hand side of \eqref{e:Z.beta.lbd.ubd}, the total contribution from such $\al$ is an exponentially small fraction of the sum. We therefore conclude
\[
\begin{array}{rl}
\EZZ_{n\be}
\hspace{-6pt}&= (\E\mathbf{Y}_{n\be})\,
	\exp\{ \slf{n}{2^k}[-1 + O(\slf{k}{2^{k/3}})] \}\\
&=\exp\{ n[ \Phi+\be[\log(e/\be)
	-\log 2^{k+1}]
	-\slf1{2^k} + O( \slf{k}{2^{4k/3}} )
	] \}.
\end{array}
\]
This is clearly optimized with $2^{k+1}\be\approx 1$, and estimating the second derivative of the exponent with respect to $\be$ implies the result.
\epf
\elem

\bppn\label{p:boundary}
The maximum of the $\zof$ auxiliary model exponent $\bPhi$ on $\simplex$ is not attained on the boundary $\pd\simplex$.

\bpf
Lem.~\ref{l:large.pf} shows that the maximum cannot be obtained on the boundary $\be=\bemax$, so it remains to show that the maximizer must be a strictly positive measure on $\supp\phi$. For $\bde\equiv(\dbde,\hbde)$ such that $\bh+t\bde$ lies in $\simplex$ for $t\ge0$ small, consider
$$\Tlog\bPhi(\bh;\bde)
\equiv\lim_{t\downarrow0}
\f{\bPhi(\bh+t\bde)-\bPhi(h)}{t\log(1/t)}
=
\dbde[(\supp\dbh)^c]
+(\slf{d}{k})\hbde[(\supp\hbh)^c]
-d\,\vh[(\supp\vh)^c].
$$
To show that $\bh\in\simplex$ is not a maximizer it suffices to exhibit $\Tlog\bPhi(\bh;\bde)>0$ for some $\bde$. In particular, it follows by convexity that for any $\bh\in\simplex$,
$\bh+t(\bhstar-\bh)\in\simplexint$ for $t>0$ small
and $\bhstar$ as in the statement of Thm.~\ref{t:first.moment.exponent}. Therefore, if $\bh$ is a maximizer such that the edge marginal has full support $\supp\vh=\msg$, then necessarily $\supp\bh=\supp\phi$, since otherwise $\Tlog\bPhi(\bh;\bhstar-\bh)>0$.

Suppose $\bh$ is a maximizer for $\bPhi$ on $\simplex$; recall $\dbh,\hbh$ must be symmetric functions. By Lem.~\ref{l:large.pf}, almost all variables are rigid except for $\asymp n 2^{-k}$ free variables; 
so some but not all edges are forcing. It is also clear that
the rigid variables will be divided roughly evenly between $\zro$'s and $\one$'s, so we obtain $\set{\ff,\zf,\of}\subseteq\supp\vh$ as well as
$\vh(\set{\free x,xx})>0$ for $x=\zro,\one$.
\bnm[1.]
\item \emph{Case $\vh(\fz)>0=\vh(\zz)$.}\\
By symmetry of $\bh$, $(\fz,\zf^{d-1})\in\supp\dbh$ and $(\fz,\zf^{k-1})\in\supp\hbh$. \\
Further $(\zf^k)\in\supp\hbh$, else
$\Tlog\bPhi(\bh;\bh'-\bh)>0$ for $\bh'\equiv(\dbh',\hbh')$ defined by
\[\dbh'=\I_{(\fz,\zf^{d-1})},\quad
(\slf{d}{k})\hbh'=\I_{(\fz,\zf^{k-1})}
	+(\slf{d}{k}-1) \I_{(\zf^k)}.\]
If $\vh(\zz)=0$ then consider
$$\dbde
	=\I_{(\zz^2,\zf^{d-2})}-\I_{(\fz,\zf^{d-1})},\quad
(\slf{d}{k})\hbde
	=2\cdot\I_{(\zz,\zf^{k-1})}
	-\I_{(\fz,\zf^{k-1})}
	-\I_{(\zf^k)};$$
this has marginal $d\vde = 2\cdot \I_{\zz}-\I_\fz-\I_\zf$
so we find $\Tlog\bPhi(\bh;\bde)=1+2-2>0$.

\item \emph{Case $\vh(\zz)>0=\vh(\fz)$.}\\
By symmetry of $\bh$, $(\zz,\zf^{k-1})\in\supp\hbh$.\\
Further
$(\zz^2,\zf^{d-2})\in\supp\dbh$, else
$\Tlog\bPhi(\bh;\bh'-\bh)>0$ for $\bh'\equiv(\dbh',\hbh')$ defined by
\[
\dbh'=\I_{(\zz^2,\zf^{d-2})},\quad
(\slf{d}{k})\hbh'=2\cdot \I_{(\zz,\zf^{k-1})}
	+ (\slf{d}{k}-2) \I_{(\zf^k)}.\]
If $\vh(\fz)=0$ then for $\bde$ as above we find
$\Tlog\bPhi(\bh;-\bde)=1+1-1>0$.
\enm
Clearly the same argument applies replacing $\zro$ with $\one$. In each case the conclusion contradicts the assumption that $\bh$ is a maximizer, concluding the proof.\footnote{In our setting we have checked $\supp\vh=\msg$ in a rather ad hoc manner. A simpler argument applies generally to any specification $\phi$ which is \emph{everywhere} positive on $\msg^d$, $\msg^k$: if $\si\notin\supp\vh$ then take $\acute\si\in\supp\vh$, and observe that $\Tlog\bPhi(\bh;\bde)>0$ for $\bde$ defined by $\dbde=\I_{(\si,\acute\si^{d-1})}-\I_{(\acute\si^d)}$, $(\slf{d}{k})\hbde=\I_{(\si,\acute\si^{k-1})}-\I_{(\acute\si^k)}$.}
\epf
\eppn

\subsection{Bethe recursion symmetries}
\label{ss:one.reduc}

\begin{figure}[ht]
\includegraphics[height=2.5in,trim=.6in 1.2in .8in .7in,clip]{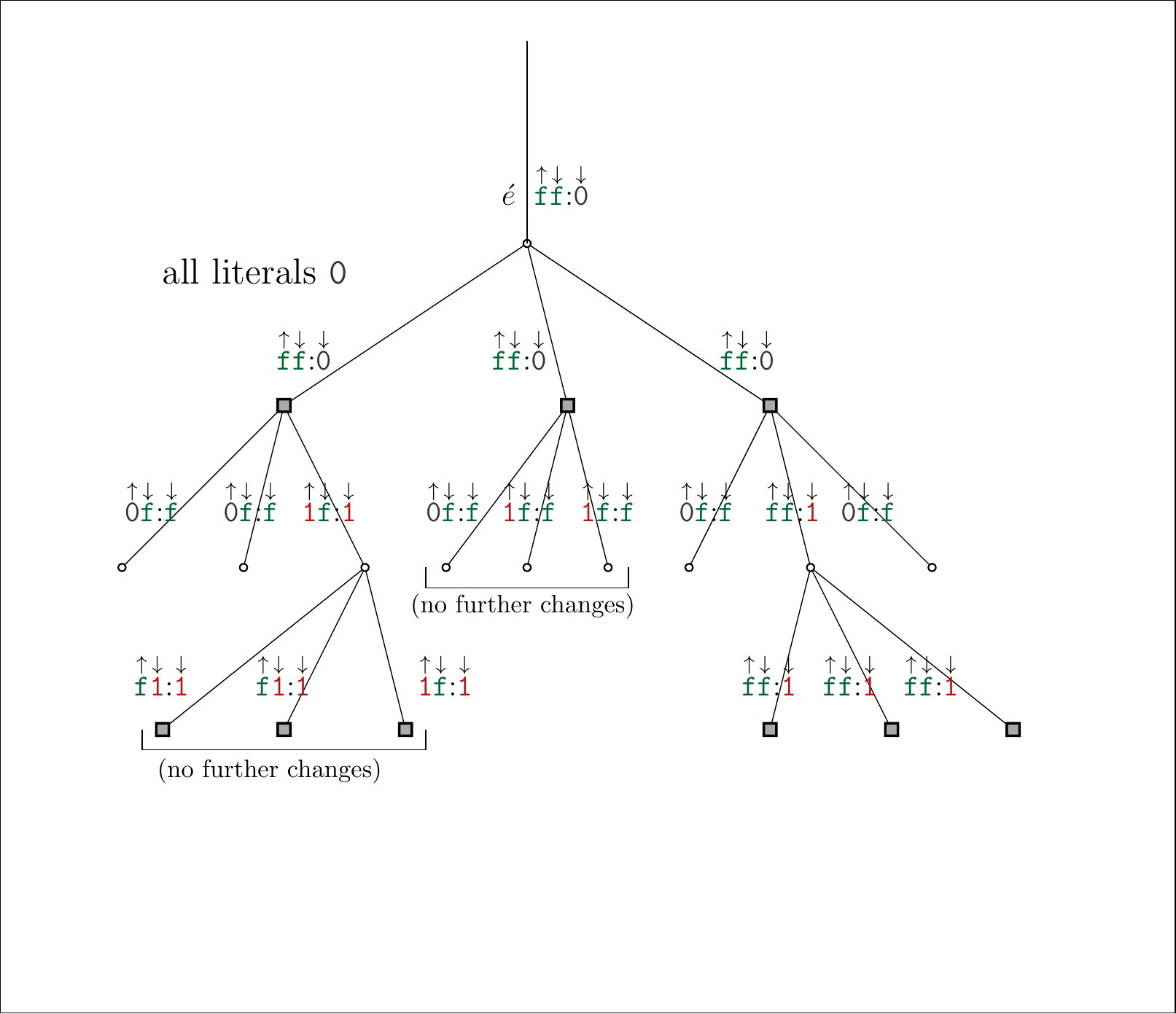}
\caption{Change of message incoming down to $\acute{e}$ is passed down
	$\etree$\\
	($\eta\acute\eta$:$\grave\eta$
	means message $\eta$ up, message $\acute\eta$ down in $\usi$,
	message $\grave\eta$ down in $\usi'$)
}
\label{f:bp.sym}
\end{figure}

Suppose $\bh$ is an interior maximizer for $\bPhi$ on $\simplex$, and so corresponds to a Bethe solution $h$. Let $\etree$ denote $\tree$ with a subtree incident to the root removed, leaving an unmatched half-edge $\acute{e}$ incident to $\rt$ (Fig.~\ref{f:bp.sym}). Consider defining a Gibbs measure on $\etree$ in the manner of \eqref{e:gibbs.aux}, with boundary law given by the Bethe solution $h$. Then the marginal law of $\si_{\acute{e}}$ will be $\hd$, and the marginal law of the tuple of spins incident to any given vertex will be $\dbh$ if the vertex is a variable, $\hbh$ if it is a clause. Further, the Gibbs measure on $\etree$ can be generated in Markovian fashion, starting with spin $\si_{\acute{e}}$ distributed according to $\hd$, generating the messages on the other $d-1$ edges incident to $\rt$ according to the conditional measure $\dbh(\dsi\,|\,\si_1=\si_{\acute e})$, and continuing iteratively down the tree.

Write $\si_{\acute{e}} \equiv \upmsg\downmsg$ where $\upmsg$ is the variable-to-clause message and $\downmsg$ the clause-to-variable message 
(in Fig.~\ref{f:bp.sym}, $\upmsg$ is directed upwards, $\downmsg$ downwards). Given any valid auxiliary configuration $\usi$ on the edges of $\etree$, changing $\downmsg$ and passing the changed message through the tree (via $\dmp_{d-1}$, $\hmp_{av}$) produces a new auxiliary configuration $\usi'$ (Fig.~\ref{f:bp.sym}). The symmetries \eqref{e:bp.sym} will follow by showing that for any fixed $\upmsg$, the effect of changing $\downmsg$ is \emph{measure-preserving} under the Gibbs measure $\nuaux$ corresponding to $\bh$. From our definition of the Gibbs measure via the boundary law, the measure-preserving property will follow by showing that the effect of changing $\downmsg$ almost surely does not percolate down the tree.

Indeed, recall that we already saw directly from the Bethe recursions that $\hh_\zz=\hh_\fz$: this came from the observation that $\hpsi$ does not distinguish between $\zz$ and $\fz$, which corresponds to the fact that changing the message incoming to a clause along a \emph{forcing} edge has no effect on the other $k-1$ edges. We also saw that $\hh_\zz=\hh_\fz$ implies $\hd_\zz=\hd_\zf$: this corresponds to the fact that if $\upmsg=\zro$, changing $\downmsg$  at most can change messages incoming to clauses in $\pd\rt$ along forcing edges, so the effect terminates before the second level of the tree.

\bpf[Proof of Thm.~\ref{t:first.moment.exponent}] By Propn.~\ref{p:boundary}, any maximizer $\bh$ for $\bPhi$ on $\simplex$ must lie in the interior $\simplexint$, and so corresponds to a solution $h$ of the Bethe recursions \eqref{e:bp}. From the above discussion it remains to show that $h$ satisfies $\hd_\ff=\hd_\fz$: meaning that in the Gibbs measure $\nuaux$ corresponding to $\bh$, changing $\downmsg$ with $\upmsg=\free$ fixed has a finite-range effect. Let $g$ correspond to $h$ via \eqref{e:bp.zof.rf}.

The effect of changing $\upmsg\downmsg$ from $\ff$ to $\fz$ can only propagate through clauses in which the parent variable and exactly one descendant variable send message $\free$, and the evaluation of the remaining $k-2$ messages under the clause literals is identically $\zro$ or $\one$.  
The vertex-preceding edges of $\etree$ whose spins will be affected by changing $\downmsg$ from $\free$ to $\zro$ form a branching process with mean
\beq\label{e:bp.sym.branching1}
(d-1)(k-1) (\slf{8}{2^k})
\hbh( \ksi\in(\ff^2,\rf^{k-2}) \giv \si_1=\ff )
\le
dk (\slf{8}{2^k} )
\f{\gd_\ff^{2} \gd_\rf^{k-2} }{ (1-\tf{4}{2^k}) \gd_\ff^{} \gd_\rf^{k-1} }
\lesssim k^2\gd_\ff/\gd_\rf,
\eeq
where the intermediate step follows from \eqref{e:bij}. Similarly, the effect of changing $\upmsg\downmsg$ from $\fz$ to $\ff$ can only propagate through clauses in which exactly one descendant variable sends message $\free$, and the evaluation of the remaining $k-1$ messages under the clause literals is identically $\zro$ or $\one$. This forms a branching process with mean
\beq\label{e:bp.sym.branching2}
(d-1)(k-1) (\slf{4}{2^k})
\hbh( \ksi\in(\rf^{k-1},\fr) \giv \si_1\in\rf )
\le dk (\slf{4}{2^k})
\f{
( 2^k-4 ) \gd_\ff^{} \gd_\rf^{k-1}
}
{(2^k-2-2k) \gd_\rf^k }
\lesssim k^2 \gd_\ff/\gd_\rf.
\eeq
To show that both processes are subcritical, we now estimate the ratio
$\DOT{u}\equiv\gd_\ff/\gd_\rf$. Recall from the proof of Lem.~\ref{l:large.pf} that the number $m\,\hbh(\rf^k)$ of fully rigid non-forcing clauses, is $m\nu_0\al = m(1-O(\tf{k}{2^k}))$ (otherwise the contribution to the partition function is an exponentially small fraction of the whole). Applying \eqref{e:bij} again we have
$$
\f{nd\be}{m\nu_0\al}
=
\hbh(\rf^k)^{-1}
\sum_{j\ge1} j\,\hbh(\perm(\ff^j,\rf^{k-j}))
=\sum_{j\ge1} j\smb{k}{j} \DOT{u}^j
= k
\DOT{u}
(1+\DOT{u})^{k-1},
$$
so we conclude $\DOT{u}=\gd_\ff/\gd_\rf=\be [1+O(\tf{k}{2^k})] \lesssim 2^{-k}$, which clearly shows that the effect of
changing $\downmsg$ given $\upmsg=\free$
does not percolate.
Therefore $h$ satisfies the symmetries \eqref{e:bp.sym},
and so corresponds to a solution $\dq$
of the frozen model recursions \eqref{e:q.rec}. Further
$1-\qrig=\gd_\ff/(\gd_\rf+\gd_\ff) \lesssim 2^{-k}$,
so Lem.~\ref{l:contract} implies $\dq=\dq^\star$ as claimed.
\epf

\subsection{Explicit form of first moment exponent}\label{ss:explicit}

We conclude this section by giving the explicit form of $\bPhistar\equiv\bPhistar_k(d)$.

\bppn\label{p:explicit}
For $k\ge k_0$, $\dlbd\le d\le \dubd$, $\bPhistar\equiv\bPhistar_k(d)$ is given by
\beq\label{e:explicit}
\bPhistar=
\log2-\log(2-\qrig)
	-d(1-k^{-1}-d^{-1})\log[ 1- 2 (q/2)^k]
	+(d-1)\log[1-(q/2)^{k-1}]
\eeq
where $\qrig$ is the unique solution of
\beq\label{e:d.of.q}
d=1+\Big(\log \smf{2(1-\qrig)}{2-\qrig}\Big)
	/ \Big( \log \smf{1-2(\qrig/2)^{k-1}}{1-(\qrig/2)^{k-1}} \Big)
\quad\text{with } 0\le 1-\qrig\le \tf{1}{2^k}.
\eeq
The function $\bPhistar$ is strictly decreasing in $d$ with $2^k[\Phi-\bPhistar]=\tf12+O(\tf{k^2}{2^k})$, and so has a unique zero $\dlbd< d_\star< \dubd$ satisfying
\beq\label{e:gap}
d_\star
= \Big(
	2^{k-1} - \smf12-\smf{1}{4\log2}  \Big) \,k\log2
	+O(\tf{k^3}{2^k})
= \dfm - \Big( \smf{1}{4\log2}-\smf{1}{6} \Big)\, k\log2
	+O(\tf{k^3}{2^k}),
\eeq
with $\dfm$ the first moment threshold of the original \acr{nae-sat} partition function \eqref{e:first.moment.threshold}.

\bpf
The equation \eqref{e:d.of.q} is a rewriting of the frozen model recursions \eqref{e:q.rec}, which by Lem.~\ref{l:contract} has a unique solution with $0\le 1-\qrig^\star \le 2^{-k}$. Throughout the following we write $\qrig\equiv\qrig^\star$,
$\qf\equiv 1-\qrig =\tf{1}{2^{k+1}} + O(\tf{k^2}{4^k})$,
and $v\equiv v(q)\equiv v_{k-1}(\qrig)$ as in \eqref{e:q.rec}.
We also abbreviate 
\beq\label{e:q.v}
\begin{array}{l}
Q\equiv Q(\qrig)\equiv(\qrig/2)^{k-1}=(1-v)/v
	\quad\text{and}\\
\DS \vr\equiv1-v = Q/(1-Q)
	= 2^{-k}[ 2-2(k-1)\qf
		+4/2^k
		+O(k^2/4^k) ].
\end{array}
\eeq
Clearly, Lem.~\ref{l:interior.bp} applies for both the $\zof$ and $\rig/\free$ auxiliary models. Substituting \eqref{e:bij} into \eqref{e:bethe} and rearranging gives
\beq\label{e:bethe.explicit}
\bPhistar=\bPhi(\bhstar)
	=\log\dbz_h + (d/k)\log\hbz_h - d\log\barz_h
	=\log\dbz_g + (d/k)\log\hbz_g - d\log\barz_g.
\eeq 
We use \eqref{e:grf.q} to calculate
\beq\label{e:explicit.bij.z}
\begin{array}{l}
\dbz_g = 2(\gh_\rr+\gh_\ff)^d-(\gh_\ff)^d
	=\smf{(2-v^d)}{2^d}
	= \smf{2(1+\qf\vr) }{2^d (1+\qf)},\\
\hbz_g = (\gd_\rf+\gd_\ff)^k-2(\gd_\rf/2)^k
	= \smf{1}{(2+\qf)^k} \smf{1+\qf\vr}{1+\vr},\\
\barz_g 
	= (\gd_\rr+\gd_\ff)(\gh_\rr+\gh_\ff)
		+ \gd_\ff\gh_\rr
	= \smf{1+\qf\vr}{2(2+\qf)}.
\end{array}
\eeq
(From the Bethe recursions \eqref{e:bp} we see that $\bpdzg\,\gd_\ff = (\gh_\ff)^{d-1}$ and $\bphzg\,\gh_\rr =(\gd_\rf/2)^{k-1}$, so we can use \eqref{e:grf.q} again to express $\bpdzg,\bphzg$ in terms of $\qf,\vr$ and confirm that the relations 
$\barz_g=\dbz_g/\bpdzg=\hbz_g/\bphzg$ of Lem.~\ref{l:interior.bp}
are indeed satisfied.) Then
\begin{align*}
\bPhistar
&=\log2-(d/k)\log(1+\vr)
	-d(1-k^{-1}-d^{-1})\log(1+\qf\vr)-\log(1+\qf)\\
&=\Phi-(d/k)\log[ (1-\tf{2}{2^k})(1+\vr) ]
	-d(1-k^{-1}-d^{-1})\log(1+\qf\vr)-\log(1+\qf)
\end{align*}
with $\Phi\equiv\log2+(d/k)\log(1-\tf{2}{2^k})$ the first-moment exponent for the original \acr{nae-sat} partition function \eqref{e:ez}.
From \eqref{e:q.v}
we have $ (1-\tf{2}{2^k})(1+\vr)= 1-\tf{2(k-1)}{2^k}\qf+O(\tf{k^2}{8^k})$, therefore
\beq\label{e:nae.frozen.gap}
\bPhistar-\Phi
=d (1-k^{-1}) \qf(\tf{2}{2^k}-\vr)
	-\log(1+\qf) + O(\tf{k^2}{4^k})
=-\qf + O(\tf{k^2}{4^k}).
\eeq
Let us now see that $\bPhistar$ is strictly decreasing in $d$. Recalling \eqref{e:q.v} that $v=Q/(1-Q)$, we find 
$$\bPhistar=
\log2+(d/k)\log(1-\qrig Q)
	+(d-1)\log\smf{1-Q}{1-\qrig Q}
	-\log(2-\qrig),$$
and rearranging gives \eqref{e:explicit}. This can be expressed
as a function of $\qrig$ alone by taking $Q=Q(\qrig)$ as in \eqref{e:q.v} and $d=d(\qrig)$ as in \eqref{e:d.of.q}.
With $\Dq$ denoting differentiation in $\qrig$, we calculate
\begin{align*}
&\Dq Q = \smf{(k-1)Q}{q} = \smf{2(k-1)}{2^k} [ 1 + O(\tf{k}{2^k})],\quad
	\Dq v = -\smf{\Dq Q}{(1-Q)^2}
	= -\smf{2(k-1)}{2^k} [ 1 + O(\tf{k}{2^k})],\\
&\Dq d=
	-[\qf (1+\qf)\log v]^{-1}
	+\smf{(d-1) \Dq v}{-v\log v}
	=(\qf\vr)^{-1}
		+O(k^2 2^k)
	= 4^k[1+O(\tf{k^2}{2^k})].
\end{align*}
The total derivative 
of $\bPhistar\equiv \bPhistar(d(\qrig))$ with respect to $\qrig$
is then straightforward to calculate: the main contribution comes from
\[\begin{array}{l}
\Dq[ (d/k) \log(1-\qrig Q) ]
	=
	-(dkQ)/(1-\qrig Q)
	+k^{-1} (\Dq d) \log(1-\qrig Q)  \\
\qquad = -k^{-1} Q (\Dq d) + O(k)
	= -(2/k)
	2^k [1+O(k^3/2^k)],\\
\text{while }
	\Dq[(d-1)\log[ (1-Q)/(1-\qrig Q) ]
	-\log(2-\qrig)]
	=O(k^2)
\end{array}
\]
Thus $\bPhistar\equiv\bPhistar_k$ is strictly decreasing on the interval $\dlbd\le d\le \dubd$ with derivative
\[\bPhistar'(d)=\smf{\Dq \bPhistar}{\Dq d}
= -\smf{2}{2^k k}[1+O(\tf{k^3}{2^k})],\]
so it must have a unique zero $\dlbd\le d_\star\le \dubd$.
We further estimate from \eqref{e:regime}, \eqref{e:q.v}, and \eqref{e:nae.frozen.gap} that $d_\star$ must satisfy
\eqref{e:gap}, concluding the proof.
\epf
\eppn

\section{Second moment of auxiliary model}\label{s:two}

In this section we compute the exponential growth rate $\bPhitstar=\lim_{n\to\infty} n^{-1}\log\E[\ZZ^2]$ of the second moment of the (truncated) frozen model partition function \eqref{e:truncated.frozen.model}. This is done in the same framework as in introduced in \S\ref{s:first}, regarding the second moment as the first moment of the partition function of \emph{pair} frozen model configurations $\uom\equiv(\ueta^1,\ueta^2)$ on the same underlying graph. The corresponding model of pair auxiliary configurations $\uta\equiv(\usi^1,\usi^2)$ has factors
$\dpsit\equiv\dpsi\otimes\dpsi$ and 
$\hpsit^a(\kta)
\equiv
\hpsi^\circ(\kta\oplus\ulit)
\equiv
\hpsi^\circ(\ksi^1\oplus\ulit)\,
\hpsi^\circ(\ksi^2\oplus\ulit)$, and rate function $\bPhit$ on the space $\simplext$ of empirical measures $\bht$ with both marginals in $\simplex$.
We can again average over literals to define the $\zof$ auxiliary model on $G$; note however that the pair auxiliary model does \emph{not} have a simple $\rig/\free$ projection as was found in \eqref{e:rf}.
In this section we prove
\bthm\label{t:second.moment}
The rate function $\bPhit$ on $\simplext$ attains its maximum only at the product measure $\bhtstar\equiv\bhstar\otimes\bhstar$, or at the measures $\bhtx$ ($x=\zro,\one$) with marginals $\bhstar$ supported on pair configurations $\uta=(\usi,x\oplus\usi)$.
\ethm

We begin with an {\it a priori} estimate in the frozen model. As before $\rig$ denotes $\set{\zro,\one}$. We partition $\rr\equiv\set{\zro,\one}^2$ into $\rreq\equiv\set{\zz,\oo}$ and $\rrne\equiv\set{\zo,\oz}$, and we decompose $\ZZ^2=\sum_\pi\ZZ^2_\pi$ where $\ZZ^2_\pi$ denotes the partition function of pair frozen configurations with associated empirical measure $\pi$ on $\pairs\equiv\set{\rreq,\rrne,\rf,\fr,\ff}$. We write $\pi_\fd\equiv\pi_\fr+\pi_\ff$, etc.; in view of \eqref{e:truncated.frozen.model} we always assume $\pi_\fd,\pi_\dotf\le\bemax$. Throughout the following we write $\al\equiv \pi_{\rreq}/\pi_\rr$ for the fraction of $\rr$-vertices taking the same spin in both coordinates.

\blem\label{l:rule.out.intermediate}
With $\al\equiv \pi_{\rreq}/\pi_\rr$, the function $\bPhit$ can only attain its global maximum on $\simplex$ either in the near-independent regime $\simplexind\subseteq\simplex$ of measures with $|2\al-1|\le k/2^{k/2}$, or in the near-identical regime $\simplexid\subseteq\simplex$ of measures with $\al\wedge(1-\al)\le 2^{-3k/4}$.

\bpf
Given empirical measure $\pi$ on $\pairs$, let $\mathbf{p}^\pi_\nu$ denote the probability, with respect to a uniformly random matching between variable and clause half-edges, that there are exactly $m\nu$ clauses which are incident to only $\rreq$ or only $\rrne$ variables: such clauses have two invalid literal assignments. Of the remaining $m(1-\nu)$ clauses, all but $O( m k^2/4^k )$ must have fewer than two frees in at least one of the two coordinates, hence must have at least \emph{four} invalid literal assignments. Therefore
$$
\TS\E[\ZZ^2_\pi]
\le
2^{n(1-\pi_\ff)}
	\tbinom{n}{n\pi}
	\sum_\nu
	\mathbf{p}^\pi_\nu
	(1-\slf{2}{2^k})^{m\nu}
	(1-\slf{4}{2^k})^{m(1-\nu)}
	\exp\{ O(n \tf{k^2}{4^k}) \}.
$$
The typical value of $\nu$ given $\pi$ is $\nubar\equiv \pi_\rr^k[\al^k+(1-\al)^k]$; conditioning and applying the local {\sc clt} (see \eqref{e:P.alpha.gamma} or the proof of Lem.~\ref{l:large.pf}) gives $\mathbf{p}^\pi_\nu\le n^{O(1)}\, \exp\{ -m\,\relent{\nu}{\nubar}\}$. The optimal contribution to the summation above comes from
$$\smf{\nu}{1-\nu}=\smf{\nubar}{(1-\nubar)(1-\vth)}
\quad\text{where }\vth\equiv \smf{2}{2^k-2}.$$
Since $\pi_\fd,\pi_\dotf\le\bemax$ we find $2^{n(1-\pi_\ff)}\tbinom{n}{n\pi} =\exp\{ n[\log 2 +H(\al) +O(k/2^k)]\}$. Combining and recalling
\eqref{e:ez}, \eqref{e:Phi2} gives
$$\TS\E[\ZZ^2_\pi]
\le
\exp\{ n[\Phi+\overline{\mathbf{a}}(\al) + O(\slf{k}{2^k})] \}.
$$
Recall from \eqref{e:frozen.lbd} that $\EZZ=e^{O(n/2^k)}$;
it therefore follows from the estimates
done in the proof of Propn.~\ref{p:sat.below}
that $\E[\ZZ^2_\pi]/(\EZZ)^2$ is exponentially small in $n$ throughout the regime $2^{-3k/4}\le\al\wedge(1-\al)\le \tf12 (1-k/2^{k/2})$.
\epf
\elem

\subsection{Near-independence regime}

In this subsection we complete our analysis of the near-independent regime $\simplexind$ to prove
\bppn\label{p:indep.maximizer}
The unique global maximizer of the restriction of $\bPhit$ to $\simplexind$ is $\bhtstar$.
\eppn

\blem\label{l:two.large.pf}
The contribution to $\E[\ZZ^2]$ from frozen configurations with $|2\al-1|\le k/2^{k/2}$ and $(2/3)\,\bemax\le \pi_\fd \vee \pi_{\dotf} \le \bemax$ is exponentially small in $n$ compared with $(\EZZ)^2$.

\bpf
Write $\pi_\free\equiv1-\pi_\rr$. Let $m\nu_j$ count the number of clauses with exactly $2j$ invalid literal assignments, and let $\nubar_j$ denote the typical value of $\nu_j$ given $\pi$:
\[
\begin{array}{rl}
\nubar_1 \hspace{-6pt}&\ge
	(\pi_\rr)^k[\al^k+(1-\al)^k]
	= \slf{2}{2^k} + O(\slf{k^2}{2^{3k/2}}),\\
\nubar_2 \hspace{-6pt}&\ge
	(\pi_\rr)^k[1-\al^k-(1-\al)^k]
	=1-\slf{2}{2^k}-k\pi_\free + O(\slf{k^2}{2^{3k/2}}),\\
\nubar_3 \hspace{-6pt}&\ge
	k(\pi_\rf+\pi_\fr) (\pi_\rr)^{k-1}[1-\al^{k-1}-(1-\al)^{k-1}]
	= k(\pi_\rf+\pi_\fr)+ O(\slf{k}{4^k}),\\
\nubar_4 \hspace{-6pt}&\ge
	k\pi_\ff (\pi_\rr)^{k-1}[1-\al^{k-1}-(1-\al)^{k-1}]
	= k\pi_\ff + O(\slf{k}{4^k}).
\end{array}
\]
By the argument of Lem.~\ref{l:rule.out.intermediate},
\[
\begin{array}{rl}
\E[\ZZ^2_\pi]
\hspace{-6pt}&\le
	2^{n(1-\pi_\ff)}\tbinom{n}{n\pi}
	\sum_\nu\mathbf{p}^\pi_\nu
	\prod_{j\ge0} (1-\slf{2j}{2^k})^{m\nu_j}\\
&\le\exp\{
	n[ (1-\pi_\ff)\log 2+H(\pi)
		+(\slf{d}{k})
		\log(1- \tf{2}{2^k}
			\sum_{j=1}^3  j\nubar_j) ]\}.
\end{array}
\]
Note that $H(\pi)$ is maximized at $\al=\slf12$, therefore
\[\begin{array}{rl}
H(\pi)
\hspace{-6pt}&\le 
-\sum_{\om\in\pairs\setminus\rr}\pi_\om\log\pi_\om
	-\pi_\rr\log\pi_\rr
	+\pi_\rr\log2\\
&\le\log2
-\sum_{\om\in\pairs\setminus\rr}\pi_\om\log\pi_\om
+\pi_\free\log(\slf{e}{2}).
\end{array}\]
From the above estimates on the $\nubar_j$ we find
\[\TS
1- \tf{2}{2^k} \sum_{j=1}^3  j\nubar_j
= (1-\slf{2}{2^k})^2
		-\tf{2k}{2^k} ( \pi_\rf+\pi_\fr+2\pi_\ff )
+O(\slf{k^2}{2^{5k/2}}).\]
Combining these estimates and recalling $\Phi=\log 2+(\slf{d}{k})\log(1-\slf{2}{2^k})$ from \eqref{e:ez} gives
\[ 
n^{-1}\log \smf{\E[\ZZ^2_\pi]}{e^{2n\Phi}}
\le-
	\underbrace{
	\pi_\rf\log \Big( \smf{\pi_\rf 2^{k+1}}{e}\Big)
	}_{\ge 1/2^{k+1}}
	-\underbrace{
		\pi_\fr\log \Big(\smf{\pi_\fr 2^{k+1}}{e}\Big)
	}_{\ge 1/2^{k+1}}
	-\underbrace{
		\pi_\ff\log \Big(\smf{\pi_\ff 4^{k+1}}{e}\Big)
	}_{\ge 1/4^{k+1}}
	+O\Big(\smf{k^2}{2^{3k/2}}\Big).
\]
Recalling \eqref{e:frozen.lbd} gives the two upper bounds
\[
n^{-1}\log \smf{\E[\ZZ^2_\pi]}{  (\EZZ)^2}
\le 
O\Big(\smf{k^{O(1)}}{2^{4k/3}}\Big)+
\begin{cases}
	-\pi_\rf\log( \pi_\rf 2^{k+1}/e )
	+\slf{3}{2^{k+1}};
		& \text{\sc (a)}\\
	-\pi_\ff\log(\pi_\ff 4^{k+1}/e)+\slf{2}{2^k}.
		& \text{\sc (b)}
\end{cases}
\]
Recall \eqref{e:truncated.frozen.model} that $\bemax=\cc/2^k$; ({\sc a}) implies that $\E[\ZZ^2_\pi]/(\EZZ)^2$ is exponentially small in $n$ for $2\pi_\rf\ge \bemax$, or symmetrically for $2\pi_\fr\ge \bemax$. However ({\sc b}) implies that $\E[\ZZ^2_\pi]/(\EZZ)^2$ is exponentially small in $n$ for $\pi_\ff \ge 4/(k2^k)$, and combining gives the result.
\epf
\elem

\bppn\label{p:boundary.two}
Any global maximizer of $\bPhit$ on $\simplexind$ must be an interior stationary point.

\bpf
Lem.~\ref{l:two.large.pf} shows that the maximum cannot be obtained on the boundary where density of frees in either coordinate is $\bemax$, so it remains to show that the maximizer must be a strictly positive measure on $\supp\phi$. For this we argue similarly as in the proof of Propn.~\ref{p:boundary}. If $\bht\equiv(\dbht,\hbht)\in\simplext$ is defined by
$\dbht(\dta)=\dbh(\dsi^1)\Ind{\dsi^1=x\oplus\dsi^2}$
and $\hbh(\ksi^1)\Ind{\ksi^1=x\oplus\ksi^2}$
for $\bh\in\simplex$ and $x\in\set{\zro,\one}$,
then clearly $\bPhit(\bht)=\bPhi(\bh)$,
so Propn.~\ref{p:boundary} implies
$$\set{(\dsi,x\oplus\dsi):\dsi\in\supp\dpsi}
\subseteq\supp\dbh,\quad
\set{(\ksi,x\oplus\ksi):\ksi\in\supp\hpsi}\subseteq\supp\hbh
\quad\text{for $x=\zro,\one$.}$$
In the following we write $r,s,x,y$ for elements of $\set{\zro,\one}$.

\bnm[1.]
\item If $\supp\vh$ does not contain $\spint{rr}{\free s}$ or $\spint{rr}{s\free}$
	then $\Tlog(\bPhit)(\bh;\bde)>0$ for
$$\begin{array}{rl}
\dbde\hspace{-6pt}
	&=\TS\sum_{x,y\ne\free}
	[
	\Ind{\dspint{\xx&\xx&\xf^{d-2}}{\fy&\yf&\yf^{d-2}}}
	-\Ind{\dspint{ \xx^2 & \xf^{d-2} }{ \yy^2 & \yf^{d-2} }]},
	\vspace{3pt}\\
(d/k)\hbde\hspace{-6pt}
	&=\TS\sum_{x,y\ne\free}
	[
	\Ind{\dspint{ \xx&\xf^{k-1} }{ \fy & \yf^{k-1} },
		\dspint{ \xx & \xf^j & \xf^{k-1-j} }
			{ \yf & \yf^j & \neg y\free^{k-1-j} }}
	-2\cdot\Ind{\dspint{ \xx & \xf^{k-1} }{ \yy & \yf^{k-1} }}
	]
	\quad \text{with }j= \flr{\tf{k-1}{2}},\vspace{3pt}\\
d\vde\hspace{-6pt}
	&=\TS
	\sum_{x,y\ne\free}[
	\Ind{\spint{\xx}{\fy},\spint{\xx}{\yf}}
	-2\cdot\Ind{\spint{\xx}{\yy}}
	].
\end{array}$$

\item If $\supp\vh$ does not contain
$\spint{r\free}{\free s}$ or $\spint{\free r}{s\free}$ then $\Tlog(\bPhit)(\bh;\bde)>0$ for
$$\begin{array}{rl}
\dbde\hspace{-6pt}
	&=\TS\sum_{x,y\ne\free}
		[\Ind{\dspint{ \xf&\fx&\xf^{d-2} }{ \fy&\yf&\yf^{d-2} }}
		-\Ind{\dspint{ \fx&\xf^{d-1} }{ \fy&\yf^{d-1} }}
		],\vspace{3pt}\\
(d/k)\hbde\hspace{-6pt}
	&=\TS\sum_{x,y\ne\free}
		[\Ind{\dspint{
			\neg\xf & \fx & \neg\xf^{k-2} }
			{ \fy & \neg\yf & \neg\yf^{k-2} } }
		-\Ind{\dspint{ \fx&\xf^{k-1}}
			{ \fy&\yf^{k-1}}
			}
		],\vspace{3pt}\\
d\vde\hspace{-6pt}
	&=\TS\sum_{x,y\ne\free}
	[ \Ind{ \spint{\xf}{\fy},\spint{\fx}{\yf} }
	-\Ind{ \spint{\fx}{\fy},\spint{\xf}{\yf} }
	].
\end{array}$$
\item If $\supp\vh$ does not contain any of
	$\spint{\xf}{\ff}$,
	$\spint{\xx}{\ff}$,
	$\spint{\fx}{\ff}$
then $\Tlog(\bPhit)(\bh;\bde)>0$ for
$$\begin{array}{rl}
\dbde\hspace{-6pt}
	&=\Ind{ \dspint{\fx & \xf^{d-1} }{\ff & \ff^{d-1} },
		\dspint{ \xx^2 & \xf^{d-2}}{\ff^2 & \ff^{d-2} } }
	-2\cdot\Ind{ \dspint{\ff^d}{\ff^d} },\vspace{3pt}\\
(d/k)\hbde\hspace{-6pt}
	&=\Ind{ \dspint{\fx&\xf^{k-1} }{\ff&\ff^{k-1}} }
	+2\cdot\Ind{\dspint{ \xx&\xf^{k-1} }{\ff&\ff^{k-1}} }
	+(2 \tf{d}{k}-3)
		\Ind{ \dspint{\xf^k}{\ff^k} }
	-2\tf{d}{k}\cdot
	\Ind{ \dspint{\ff^k}{\ff^k} },\vspace{3pt}\\
d\vde\hspace{-6pt}
	&=\Ind{ \spint{\fx}{\ff} }
	+2\cdot\Ind{ \spint{\xx}{\ff} }
	+(2d-3) \Ind{ \spint{\xf}{\ff} }
	-2d\cdot\Ind{ \spint{\ff}{\ff} }.
\end{array}$$
\enm
It follows by symmetry considerations that $\supp\vh=\msg^2$, hence any maximizer $\bh$ of $\bPhit$ must be positive on $\supp\phi$ since otherwise $\Tlog\bPhit(\bh;\bhtstar-\bh)$ would be positive.
\epf
\eppn

\blem\label{l:contract.two}
The pair frozen model tree recursions
on measures $\tq\equiv(\dtq,\htq)$
have the unique solution
$\tq=q^\star\otimes q^\star$
in the regime
$\dtq(\set{\ff,\rf,\fr})\lesssim 2^{-k}$,
$(\dqzz+\dqoo)/(\dqzo+\dqoz)=1+O(\slf{k}{2^{k/2}})$.

\bpf
The pair frozen model tree recursions are as follows. Write $\dqeq\equiv\dqzz+\dqoo$ and $\dqne\equiv\dqzo+\dqoz$.
By assumption, $\dqeq=(1+\ep)\dqne$ with $|\ep|\lesssim\tf{k}{2^{k/2}}$.
The clause recursions are
\[\begin{array}{l}
\hqeq\equiv\hqzz=\hqoo
	= (\slf{2}{2^k})\,(\dqeq)^{k-1},\quad
\hqne\equiv\hqzo=\hqoz
	= (\slf{2}{2^k})\, (\dqne)^{k-1},\\
\hqrf\equiv\hqzf=\hqof
	=(\slf{2}{2^k})\,[ (\dqrd)^{k-1}
		- (\dqeq)^{k-1}-(\dqne)^{k-1}],\\
\hqfr\equiv\hqfz=\hqfo
	=(\slf{2}{2^k})\,[ (\dqdr)^{k-1}
		- (\dqeq)^{k-1}-(\dqne)^{k-1}],\\
\hqff
	=1-2[\hqeq+\hqne+\hqrf+\hqfr]
	=1-2[\hqzz+\hqzo+\hqzf+\hqfz]
	=1-O(2^{-k}).
\end{array}\]

\medskip\noindent The variable recursions are
(with $\DOT{c}$ the normalizing constant)

\medskip\noindent\hspace{5pt}
$\begin{array}{l}
\DOT{c}\dqff
= \hqff^{d-1},\\
\DOT{c}\dqfz
=\DOT{c}\dqfo
	=(\hqff+\hqfr)^{d-1}-\hqff^{d-1},\quad
\DOT{c}\dqzf
=\DOT{c}\dqof
	=(\hqff+\hqrf)^{d-1}-\hqff^{d-1},\\
\DOT{c}\dqzz
=\DOT{c}\dqoo
	=(\hqeq+\hqrf+\hqfr+\hqff)^{d-1}
	-(\hqrf+\hqff)^{d-1}
	-(\hqfr+\hqff)^{d-1}
	+\hqff^{d-1}\\
\DOT{c}\dqzo
=\DOT{c}\dqoz
	=(\hqne+\hqrf+\hqfr+\hqff)^{d-1}
	-(\hqrf+\hqff)^{d-1}
	-(\hqfr+\hqff)^{d-1}
	+\hqff^{d-1}.
\end{array}$

\smallskip\noindent
By the assumption that
$|\ep|\lesssim\tf{k}{2^{k/2}}$,
the clause recursions give
$\hqrf=\tf{2}{2^k}[1+O(\tf{k}{2^k})]=\hqfr$,
and therefore from the variable recursions
we must have $\dqrf=(1+\de)\dqfr$ with $|\de|\lesssim 2^{-k}$.
But then the clause recursions give
$|\hqrf-\hqfr| \lesssim \de \tf{k}{4^k}$,
consequently
$$\DSf{\dqfr}{\dqrf}
=\DSf{
	1+O( \de \tf{dk}{4^k} )- \hqff^{d-1} / (\hqff+\hqrf)^{d-1}
	}{1- \hqff^{d-1} / (\hqff+\hqrf)^{d-1}}
=1+O(\de \tf{k^2}{2^k}),$$
proving that the recursion contracts to $\de=0$, i.e.\
$\dqrf=\dqfr$
and $\hqrf=\hqfr$. Similarly, the clause recursions give $\hqeq=\hqne = O( \ep \tf{k}{4^k})$, and substituting this into the variable recursions gives $\dqeq/\dqne=1+O(\ep \tf{k^2}{2^k})$,
so we also have contraction to $\ep=0$,
$\dqeq=\dqne\equiv \tf12 \dqrr$.

It remains to show $(\dqrr,\dqrf,\dqff)=(\dq^2,\dq(1-\dq),(1-\dq)^2)$ with $\dq=\dq^\star$. To this end write
$$q\equiv \dqrr+\dqrf \equiv 4\dqzz+2\dqzf,\quad
\qp\equiv \dqrr/q \equiv
\dqzz/(2\dqzz+\dqzf).$$
Writing $Q\equiv(q/2)^{k-1}$, $\Qp\equiv(\qp/2)^{k-1}$,
and $Q^\star\equiv (\dq^\star/2)^{k-1}$, we have
$$
\text{\footnotesize
$\DS
1-q=\f{\ip{(0,2,-1)}{W^{d-1}}}{\ip{(4,-4,1)}{W^{d-1}}}, \ \
1-\qp=\f{\ip{(0,1,-1)}{W^{d-1}}}{\ip{(2,-3,1)}{W^{d-1}}}, \
\begin{array}{r}
\hqff+2\hqzf+\hqzz=1-2Q+Q\Qp\equiv W_1,\\
\hqff+\hqzf=1-3Q+2Q\Qp\equiv W_2,\\
\hqff = 1-4Q+4Q\Qp\equiv W_3.
\end{array}$}
$$
By assumption, $x\equiv 2^k(|q-\dq^\star|+|\qp-\dq^\star|) \lesssim1$,
so $Q  + O(x \tf{k}{4^k}) = Q^\star  = \Qp + O(x \tf{k}{4^k})$
and consequently
$W^{d-1}[1+O( x \tf{dk}{4^k} )] = (W^\star)^{d-1}$.
It follows that
$$\text{\small$\DS\f{1-q}{1+O(x\tf{k^2}{2^k})}
=\f{1-\qp}{1+O(x\tf{k^2}{2^k})}
=\f{\ip{(0,2,-1)}{(W^\star)^{d-1}}}{\ip{(4,-4,1)}{(W^\star)^{d-1}}}
=1-\dq^\star$}
$$
implying that the recursion contracts to $x=0$, $q=\qp=\dq^\star$ as claimed.
\epf
\elem

\bpf[Proof of Propn.~\ref{p:indep.maximizer}]
By Propn.~\ref{p:boundary.two} and Lem.~\ref{l:interior.bp},
any maximizer $\bh$ of $\bPhit$ on $\simplexind$ corresponds to some solution $h$ of the Bethe recursions
for the pair $\zof$ auxiliary model. We now show that $h$ must satisfy the symmetries
\eqref{e:bp.sym}:
that is, $\hd(\bm{oi})=\hd(\bm{oi}')$, where $\bm{o}$ now indicates the outgoing \emph{pair} of variable-to-clause messages, and $\bm{i}$ or $\bm{i}'$ indicates the incoming pair of clause-to-variable messages.
Let
$\mathscr{R}_=
\equiv\set{\spint{\zf}{\zf},\spint{\of}{\of}}$
and $\mathscr{R}_{\ne}
\equiv\set{\spint{\zf}{\of},\spint{\of}{\zf}}$,
and write
$\mathscr{D}_j\equiv
\mathscr{R}_=^j
\cup
\mathscr{R}_{\ne}^j$
and
$\mathscr{R}\equiv\mathscr{D}_1$.
Analogously to the first-moment symmetries seen directly from the Bethe recursions, in the pair model it is easily seen that
$\hh(\spint{\si}{xx})=\hh(\spint{\si}{\free x})$ and $\hd(\spint{\si}{xx})=\hd(\spint{\si}{x\free})$ for any $\si\in\msg$ and $x\in\set{\zro,\one}$.
Further, the $\zro/\one$ symmetry in the clause factors implies
$\hh(\tau)=\hh(\neg\tau)$
and $\hd(\tau)=\hd(\neg\tau)$
for any $\tau\in\msg^2$.
It remains to show that
$\hd(\spint{\si}{\free x})=\hd(\spint{\si}{\ff})$
for $x\in\set{\zro,\one}$.

\smallskip\noindent\emph{Estimates on messages.}
The number of clauses incident to any variable which are free in either coordinate is $\lesssim m\tf{k}{2^k}$, while an~easy {\it a priori} estimate implies that the number of fully-rigid clauses which are non-forcing is $\asymp m$. Recalling \eqref{e:bij} then gives
\beq\label{e:free.message.ubd}
\slf{k}{2^k}
\gtrsim
\DSf{k\,\hbh(
	\set{
	\spint{\ff}{\rf},\spint{\rf}{\ff},\spint{\ff}{\ff}
	} , \mathscr{R}^{k-1} )
	}{ \hbh( \mathscr{R}^k )} 
\ge [1-O(\tf{k}{2^k})]
\,
\DSf{k\,
[
\hd(\spint{\ff}{\rf})
+\hd(\spint{\rf}{\ff})
+\hd(\spint{\ff}{\ff})
]}{\hd(\spint{\rf}{\rf})},
\eeq
where the last inequality follows because all the $\hpsit$ factor weights involved in the application of \eqref{e:bij} are $1-O(\tf{k}{2^k})$.

We now estimate the ratio $\acute{\gm}\equiv\hd(\spint{\zf}{\zf})/\hd(\spint{\zf}{\of})$. Another application of \eqref{e:bij} gives
$$
\gm \equiv\smf{\al}{1-\al}=
\smf{\pi_\rreq}{\pi_\rrne}
=\DSf{\hh( \spint{\set{\zz,\zf}}{\set{\zz,\zf}} )^d
-\hh( \spint{\set{\zz,\zf}}{\zf} )^d
-\hh( \spint{\zf}{\set{\zz,\zf}} )^d
+\hh(\spint{\zf}{\zf})^d}
{\hh( \spint{\set{\zz,\zf}}{\set{\oo,\of}} )^d
-\hh( \spint{\set{\zz,\zf}}{\of} )^d
-\hh( \spint{\zf}{\set{\oo,\of}} )^d
+\hh(\spint{\zf}{\of})^d}
$$
where we have used the symmetry
$\hh(\spint{\si}{xx})=\hh(\spint{\si}{\free x})$ noted above.
The ratio
$\acute{\gm}$
is given by the same expression with $d-1$ in place of $d$.
Writing $\hd^\otimes$ for the product measure with marginals $\hd$,
the Bethe recursions give
$$\bphzh\,\hh(\spint{\zz}{\zz})
= \tf{2}{2^k} \hd^\otimes(\mathscr{D}_{k-1}),\quad
\bphzh\,\hh(\spint{\zz}{\zf})
= \tf{2}{2^k} \hd^\otimes(\mathscr{R}^{k-1}\setminus\mathscr{D}_{k-1}),
\quad
\bphzh\,\hh(\spint{\zf}{\zf})
= (1-O(\tf{k}{2^k}))
\hd^\otimes(\mathscr{R}^{k-1}),
$$
where the last estimate uses \eqref{e:free.message.ubd}. 
Thus
$\hh( \spint{\set{\zz,\zf}}{\set{\zz,\zf}} )
= 
[1+\tf{2}{2^k}+O(\tf{k}{4^k})] \hh( \spint{\set{\zz,\zf}}{\zf} )$,
and so
\beq\label{e:message.bias}
\acute{\gm}
	\equiv 
	\DSf{\hd(\spint{\zf}{\zf})}{\hd(\spint{\zf}{\of})}
	=e^{O(1/2^k)}
	\DSf{\hh( \spint{\set{\zz,\zf}}{\set{\zz,\zf}} )^{d-1}}
	{\hh( \spint{\set{\zz,\zf}}{\set{\zz,\zf}} )^{d-1}}
	=e^{O(1/2^k)}
	( \gm e^{O(1/2^k)})^{(d-1)/d}
	=1+O(\tf{k}{2^{k/2}}),
\eeq
where the last step uses the assumption that 
$\bh$ lies in $\simplexind$.

\smallskip\noindent\emph{Finite-range effect of changed incoming message.}
We now show
$\hd(\spint{\si}{\free x})=\hd(\spint{\si}{\ff})$
for $x\in\set{\zro,\one}$.
The effect propagates through clauses which in the second copy are as described in the proof of Thm.~\ref{t:first.moment.exponent}:
that is, in the second copy,
exactly one descendant variable sends message $\free$,
and the evaluation of all the incoming $\zro/\one$ messages
(of which there are $k-2$ or $k-1$ depending on whether the parent variable sends $\free$ or not) under the clause literals is identically $\zro$ or $\one$. The mean of the branching process is bounded as in \eqref{e:bp.sym.branching1}~and~\eqref{e:bp.sym.branching2}
except that we must now condition on the \emph{pair} spin
$\tau_1$ on the edge preceding the clause.

We now explain the rather delicate case where the clause is forcing to its parent variable in the \emph{first} coordinate. Conditioned on spin $\tau_1=\spint{\zz}{\ff}$ on the preceding edge, the probability of having a clause as described above is (using \eqref{e:bij}~and~\eqref{e:message.bias})
$$\le
\DSf{
\tf{2}{2^k}
\hd(\spint{\zz}{\ff})
(k-1)
\hd(\spint{\rf}{\ff})
\hd^\otimes( \mathscr{D}_{k-2} )}
{\tf{2}{2^k}
\hd(\spint{\zz}{\ff})
\hd^\otimes(\mathscr{R}^{k-1} \setminus \mathscr{D}_{k-1}) }
\lesssim
\tf{k}{2^k} \hd(\spint{\rf}{\ff}) /\hd(\spint{\rf}{\rf})
\lesssim \tf{k}{4^k}
$$
and this is $\ll d^{-1}$ so the propagation through clauses started from $\tau_1=\spint{\zz}{\ff}$ is subcritical. The calculations for the remaining cases of $\tau_1$ are similar but easier, and so are left to the reader. We therefore see that $h$
satisfies the symmetries \eqref{e:bp.sym},
and so corresponds to a solution $q\equiv(\dq,\hq)$ of the pair frozen model recursions.
By \eqref{e:free.message.ubd} and \eqref{e:message.bias}
this solution falls in the regime
of Lem.~\ref{l:contract.two},
which uniquely identifies $\bh$ as $\bhtstar$.
\epf

\subsection{\emph{A priori} rigidity estimate}\label{ss:rigid}

Recalling Lem.~\ref{l:rule.out.intermediate}, let $\idZZ$ denote the contribution to $\ZZ^2$ from the near-identical regime $\bh\in\simplexid$. In this subsection we prove

\bppn\label{p:second.moment.id}
For $k\ge k_0$, $\dlbd\le d\le \dubd$, and $n\ge n_0(k)$, $\E[\idZZ]\le n^{O(1)}\,\EZZ$.
\eppn

\blem\label{l:forced.once}
Given a frozen configuration $\ueta$, for $1\le j\le k$ let $m\nu_j$ count the number of clauses incident to exactly $j$ $\ueta$-free variables, and write $\smash{\nu_{\ge2}\equiv 1-\nu_0-\nu_1}$. Let $\mf$ count the number of $\ueta$-forcing clauses, and let $\gm$ denote the fraction of rigid variables which are $\ueta$-forced only once. Then for $k\ge k_0$, $n\ge n_0(k)$ it holds that
\[
\begin{array}{ll}
\begin{array}{l}
	\E[\ZZ_{n\be};(\Om_A)^c]
	\le(\EZZ)\,\exp\{ -10 n k^2/2^k \}
	\end{array}
	&\begin{array}{l}
	\text{for }\Om_A\equiv\set{\nu_{\ge2}
	\le k^3 \be^2};
	\end{array}\\
\begin{array}{l}
	\E[\ZZ_{n\be};(\Om_B)^c]
	\le (\EZZ)\,\exp\{ -10 n k^2/2^{k/2} \} \\ {}
	\end{array}
	&\begin{array}{l}
	\text{for }\Om_B\equiv \\ {}
	\end{array}
	\hspace{-2pt}
	\bigg\{\hspace{-3pt}
	\begin{array}{rl}
		| 1-\mf/(m \cdot \slf{2k}{2^k}) |
			\hspace{-6pt}&\le 2^{-k/8}\\
		\text{and }\gm
			\hspace{-6pt}&\le k^2/2^{k/2}\end{array}
	\hspace{-3pt}\bigg\}
\end{array}
\]

\bpf
As in the proof of Lem.~\ref{l:large.pf}, let $\mathbf{p}^\be_{\nu_0,\nu_1}$ denote the probability of $\nu_0,\nu_1$
with respect to a uniformly random matching between
clause half-edges and variable half-edges with density $\be$ of frees.
Conditioned on all fully-rigid clauses being satisfied, the number $\mf$ of forcing clauses is distributed $\Bin(m\nu_0,\vth)$ with $\vth\equiv \slf{2k}{(2^k-2)}$. Conditioned on $\mf$, the probability of having $\gm$-fraction of the rigid variables forced only once is
\[\TS\mathbf{a}^\be_{\mf}(\gm)
	= \P( \sum_{i=1}^{np}\Ind{F_i=1}=np\gm
	\giv \sum_{i=1}^{np} F_i = \mf ),\quad
	F_i\sim\Bin(d,\al)\]
(where $0<\al<1$ may be arbitrarily chosen).
We therefore bound
\[
\E[\ZZ_{n\be};\nu_{\ge2},\mf,\gm]
\le 2^{n(1-\be)} \tbinom{n}{n\be}
	\sum_{\nu_0+\nu_1=1-\nu_{\ge2}}
		\mathbf{p}^\be_{\nu_0,\nu_1}
		(1-\slf{2}{2^k})^{m\nu_0}
		\,\bin_{m\nu_0,\vth}(\mf)\,\mathbf{a}^\be_{\mf}(\gm).
\]
From the trivial bound $\nu_0\ge1-k\be$, together with our estimate \eqref{e:frozen.lbd} that $\EZZ=e^{O(n/2^k)}$,
\[
\E[\ZZ_{n\be};\nu_{\ge2},\mf,\gm]
\le  (\EZZ)\,\exp\{ O(\slf{nk}{2^k})\}\,
		\sum_{\nu_0+\nu_1=1-\nu_{\ge2}}
		\mathbf{p}^\be_{\nu_0,\nu_1}
		\,\bin_{m\nu_0,\vth}(\mf)
		\,\mathbf{a}^\be_{\mf}(\gm).
\]
Summing over $\mf,\gm$ and simply upper bounding $\smash{\sum_{\mf,\gm}\bin_{m\nu_0,\vth}(\mf)\,\mathbf{a}^\be_{\mf}(\gm)\le1}$
yields 
the bound on $\E[\ZZ_{n\be};(\Om_A)^c]$, recalling that the typical value of $m\nu_{\ge2}$ is $\lesssim mk^2\be^2$. To bound
$\E[\ZZ_{n\be};(\Om_B)^c]$, we first estimate
\[
\bin_{m\nu_0,\vth}(\mf) \le \exp\{ -n 2^{-k/3} \}
\quad\text{on the event }
	| 1-\mf/(m \cdot\slf{2k}{2^k} ) | \ge 2^{-k/8}.
\]
On the complementary event
$\smash{| 1-\mf/(m \cdot\slf{2k}{2^k} ) | < 2^{-k/8}}$,
in the above expression for $\smash{\mathbf{a}^\be_{\mf}(\gm)}$
we can set $\al=\mf/(ndp)=(\slf{2}{2^k})[1+O(2^{-k/8})]$,
and apply the local {\sc clt} to bound
\[\begin{array}{rll}
\mathbf{a}^\be_{\mf}(\gm)
\hspace{-6pt}&\le \exp\{ -np\,\relent{\gm}{\overline{\gm}} \}
	&\text{with }
	\overline{\gm}
	\equiv
	d\al(1-\al)^{d-1}
	= 2^{-k} k\log2 [ 1+O(\slf{k^2}{2^{k/8}}) ]\\
&\le
	\exp\{ -n  k^{5/2}/2^{k/2} \}
	&\text{for }
	\gm\ge k^2/2^{k/2}.
\end{array}
\]
Combining these estimates gives the bound on $\E[\ZZ_{n\be};(\Om_B)^c]$.
\epf
\elem

We now decompose $\ZZ^2=\sum_\pi\ZZ^2[\pi]$ where $\ZZ^2[\pi]$ denotes the contribution from empirical measure $\pi$ on $\set{\zro,\one,\free}^2$. For $j=1,2$ we write $\pi^j$ for the projection of $\pi$ onto the $j$-th coordinate,
e.g.\ $\pi^1\equiv(\pi_\zd,\pi_\od,\pi_\fd)$,
and we decompose $\ZZ=\sum_{\pi^1}\ZZ[\pi^1]$.

\blem\label{l:ap.near.id}
For any empirical measure $\pi$ on $\set{\zro,\one,\free}^2$ with $\pi^1_\free\vee\pi^2_\free\le\bemax$ and with $\De\equiv n\pi(\eta^1\ne\eta^2) \le \slf{n}{2^{k/2}}$, it holds for $k\ge k_0$, $n\ge n_0(k)$ that
$$\E[\ZZ^2[\pi]]
\le
e^{-n/2^{k/2}}\EZZ
+n^{O(1)}\,2^{ -\De k/10} 
	(\E[\ZZ[\pi^1]]+\E[\ZZ[\pi^2]]).$$

\bpf
Write $p\equiv\pi_\rd$ and $\be\equiv\pi_\fd\equiv1-p$.
Given any $\ueta^1\in\set{\zro,\one,\free}^V$, the number of choices for $\ueta^2$ for which $|\set{v:\eta^1\ne\eta^2}| \le n/2^{k/2}$ is (crudely) upper bounded by $\exp\{ O(nk/2^{k/2}) \}$ even in absence of satisfiability constraints. Combining with Lem.~\ref{l:forced.once} gives $\E[\ZZ^2[\pi];(\Om_B)^c]\le (\EZZ)\,\exp\{ -5nk^2/2^{k/2} \}$, so we hereafter restrict consideration to the event $\Om_B$.

For the remainder of the proof let $\uom\equiv(\ueta^1,\ueta^2)$ be a \emph{fixed} spin configuration with empirical measure $\pi$, and for $\om\in\set{\zro,\one,\free}^2$ write $V_\om\equiv\set{v\in V:\om_v=\om}$. Decompose $\Om_B$ as the disjoint union of events $\Om_{B,\bm{x}}$ where $\bm{x}\equiv(\nu_0,\nu_1,\mf,\gm)$ is defined as in the statement of Lem.~\ref{l:forced.once} \emph{with respect to} $\ueta^1$. Let $\Rne\equiv V_\rd\setminus V_\rreq$ and $|\Rne|\equiv np\ep = n\pi_{\rrne}+n\pi_\rf$, and write $F_\de$ for the event that exactly $np\ep\de$ in $\Rne$ are $\ueta^1$-forced only once. We then bound
\begin{align*}
\E[\ZZ^2[\pi];\Om_B]
&\le
\sum_{\bm{x}}
\E[ \ZZ[\pi^1];\Om_{B,\bm{x}}]
\sum_\de
\mathbf{c}^{\pi,\bm{x},\de}_{\rig}\,
\mathbf{c}^\pi_\free\,
\P( \ueta^2\text{ valid} \giv
	\ueta^1\text{ valid},
	\Om_{B,\bm{x}},
	F_\de)\\
&\text{where }
	\mathbf{c}^{\pi,\bm{x},\de}_{\rig}
	\equiv
	2^{np\ep}
	\smb{np\gm}{np\ep\de}
	\smb{np(1-\gm)}{np\ep(1-\de)}
	\text{ and }
	\mathbf{c}^\pi_\free
	\equiv\smb{n\be}{n\pi_\fr}
	\le 2^{n\be}.
\end{align*}

\smallskip\noindent
\emph{Constraints on clauses incident to $\set{\rrne,\rf}$-variables.}\\
Let $Q_\rig\supseteq\set{\ueta^2\text{ valid}}$ denote the event that the variables in $\Rne$ do not violate any clauses. On the event $F_\de$ there must be at least $np\ep\de+2np\ep(1-\de)=np\ep(2-\de)$ clauses $\ueta^1$-forcing to $\Rne$, and for $Q_\rig$ to occur, each such clause must be incident to at least one other $\Rne$-variable. The density of edges from $\Rne$ among the non-$\ueta^1$-forcing edges is $\le ndp\ep/(mk\nu_0-\mf)\le 2\ep$, so
\[\begin{array}{rl}
	\mathbf{p}^{\pi,\bm{x},\de}_\rig
	\hspace{-6pt}&\equiv \P( Q_\rig
		\giv \ueta^1\text{ valid},\Om_{B,\bm{x}},F_\de)\\
	&\le
		\P( D_a>0 \ \forall a\le np\ep(2-\de)
		\giv
		\sum_{a=1}^{m\nu_0} D_a= (mk\nu_0-\mf)2\ep
		),\\
	&\text{with }D_a \text{ independent random variables distributed as}\\
		&{}\qquad\text{$\mathrm{Bin}(k-1,2\ep)$ for $a\le \mf$,
		$\mathrm{Bin}(k,2\ep)$ for $a>\mf$.}
	\end{array} \]
Lem.~\ref{l:forcing} gives (crudely) that $\mathbf{p}^{\pi,\bm{x},\de}_\rig
\le e^{O(n\ep)}\,[1-(1-\ep)^{k-1}]^{np\ep(2-\de)}
\le e^{O(n\ep)}\,(2k\ep)^{np\ep(2-\de)}$.
Combining with $\mathbf{c}^{\pi,\bm{x},\de}_\rig$ and rearranging gives
\[ \begin{array}{rl}
	(np)^{-1}
	\log (\mathbf{c}^{\pi,\bm{x},\de}_\rig
		\mathbf{p}^{\pi,\bm{x},\de}_\rig)
	\hspace{-6pt}
	&\le \gm\,[H(\slf{\ep\de}{\gm})
		+ (\slf{\ep\de}{\gm})\log(\slf{\ep}{\gm}) ]
		+(1-\gm)\,[ H(\tf{\ep(1-\de)}{1-\gm} )
			+ \tf{\ep(1-\de)}{1-\gm}\log\ep ]\\
	&\qquad {}+\ep\de\log\gm
		+\ep(1-\de)\log\ep\\
	&\le 2\ep+\ep\log(\ep+\gm),
	\end{array} \]
where we have used the trivial inequality $H(x)+x\log c\le \log(1+c)\le c$. Recall $\ep\lesssim 2^{-k/2}$ by assumption, and $\gm\le k^2/2^{k/2}$ by the restriction to $\Om_B$, therefore
\beq\label{e:ap.id.ubd}
\mathbf{c}^{\pi,\bm{x},\de}_\rig\mathbf{p}^{\pi,\bm{x},\de}_\rig
\le e^{O(n\ep)}\,
	\exp\{ -np\ep (k/2) \log2[1- O( \tf{\log k}{k} ) ]\}.
\eeq
Recalling $\mathbf{c}^\pi_\free\le 2^{n\be}$ 
we see that 
\[ \text{if }10\ep\ge\be\text{ then }
	\E[\ZZ^2[\pi];\Om_B]
	\le \f{\E[\ZZ[\pi^1];\Om_B]}{ 2^{n\ep k/3-n\be} }
	\le  \f{\E[\ZZ[\pi^1];\Om_B]}{
		2^{n(\ep+\be)k/5} }
	\le\f{\E[\ZZ[\pi^1];\Om_B]}{2^{\De k/5}}.\]

\smallskip\noindent\emph{Forcing of $\fr$-variables.}\\
Now suppose $\ep\lesssim\be$: the number of choices for $\ueta^2$ is then $\le\exp\{ O(nk/2^k) \}$, so combining with Lem.~\ref{l:forced.once} gives in this case $\E[\ZZ^2[\pi];(\Om_A)^c] \le (\EZZ)\,\exp\{ -5nk^2/2^k\}$. Therefore we restrict consideration hereafter to the event $\Om_A$.

On $\Om_A$, consider the event $Q_\free\supseteq\set{\ueta^2\text{ valid}}$ that every $\fr$-variable is $\ueta^2$-forced, conditioned on the preceding events $\set{\ueta^1\text{ valid}}$, $\Om_{B,\bm{x}}$, $F_\de$, and $Q_\rig$. A clause can be $\ueta^2$-forcing to an $\fr$-variable in only one of two ways:
\bnm[1.]
\item For $v\in V_\fd$ let $\Avp$ denote the number of clauses $a\in\pd v$ which are incident to no $\ueta^1$-free variables besides $v$; note that $\smash{\sum_{v\in V_\fd} \Avp=m\nu_1\le mk\be}$. A clause $a\in\pd v$ of this type will be $\ueta^2$-forcing to $v$ for certain arrangements of literals and of spins $\rr,\rrne$ among the neighbors $u\in\pd a\setminus v$. Since $v\in V_\fd$ the clause $a$ is conditioned not to be $\ueta^1$-forcing, so $\ueta^2$-forcing arrangements of $a$ will occur with conditional probability $\lesssim \ep k/2^k$.

\item In the $m\nu_{\ge2}$ clauses incident to more than one $\ueta^1$-free variable, the conditioning so far gives no information about the arrangement of the literals. Therefore, in each such clause distinguish a uniformly random edge to be \emph{potentially} $\ueta^2$-forcing. For $v\in V_\fd$ let $\Av$ denote the number of such edges incident to $v$, and write $\smash{\mA\equiv\sum_{v\in V_\fd}\Av \le m\nu_{\ge2}}$. A clause $a\in\pd v$ of this type is $\ueta^2$-forcing to $v$ for certain arrangements of literals and of spins $\rr,\fr$ among the neighbors $u\in\pd a\setminus v$. In particular, by definition at least one neighbor $u\in\pd a\setminus v$ is $\ueta^1$-free, thus $\ueta^2$-forcing arrangements of $a$ will occur with conditional probability $\lesssim k\pi_\fr/(2^k\be)$.
\enm
Crudely bounding $\Avp\le d$, there exists a uniform constant $C$ such that
\[
\mathbf{p}^{\pi,\bm{x},\de}_\free
	\equiv \P(Q_\free\giv
		\ueta^1\text{ valid},\Om_{B,\bm{x}},F_\de,Q_\rig,\Om_A)
	\le \E\Big[ \prod_{v\in V_\fr}
		\Big\{
		1-\Big(1-  \smf{C\ep k}{2^k} \Big)^d
		\Big( 1-\smf{Ck\pi_\fr}{2^k\be}  \Big)^{\Av}
		\Big\}\Big].
\]
In the above we have used the restriction to $\Om_A$ to see that the total number $\mA$ of {\sc a}-edges is $\le m\nu_{\ge2} \le mk^3\be^2 =(k^2\be)(nd\be)$: since this is much smaller than the total number $nd\be$ of half-edges leaving $V_\fd$, our $\lesssim k\pi_\fr/(2^k\be)$ estimate on the $\ueta^2$-forcing probability of each {\sc a}-edge remains valid even after revealing the states of some of the other {\sc a}-edges. 
Now apply the local {\sc clt} to bound
\[\mathbf{p}^{\pi,\bm{x},\de}_\free\le n^{O(1)}\,
	\Big[ 1-\Big(1-  \smf{C\ep k}{2^k} \Big)^d
	\E\Big\{
	\Big( 1-\smf{Ck\pi_\fr}{2^k\be}
	\Big)^B\Big\}\Big]^{n\pi_\fr},\quad
	B\sim\mathrm{Bin}\Big( d, \smf{\mA}{nd\be} \Big).\]
Recalling $\E[(1-x)^{\mathrm{Bin}(d,p)}]=( 1 - px )^d$ we see that
\[\mathbf{p}^{\pi,\bm{x},\de}_\free
	\le n^{O(1)}\,
	\Big(
		\smf{5C d\ep k}{2^k}
		+\smf{5C k\pi_\fr}{2^k\be} \smf{m\nu_{\ge2}}{nd\be}
		\Big)^{n\pi_\fr}
	\le n^{O(1)}\,
		[5Ck^2(\ep + y)]^{n\pi_\fr}
		\quad\text{with }y\equiv \smf{k\pi_\fr}{2^k}.\]
Again recalling $H(x)+x\log c\le\log(1+c)\le c$ we bound
\[\mathbf{c}^{\pi,\bm{x},\de}_\free\mathbf{p}^{\pi,\bm{x},\de}_\free
	\le\begin{cases}
	n^{O(1)}\,\exp\{ n\be[ H(y) + y\log(k^{O(1)}\ep ) ] \}
		\le n^{O(1)}\,e^{O(n\ep)}
		& \text{if }y\le\ep\le10\be;\\
	n^{O(1)}\, k^{ O(n\pi_\fr) }\exp\{ n\be[ H(y)+y\log y ] \}
		\le n^{O(1)}\, k^{ O(n\pi_\fr) }
		& \text{if }0\le\ep\le y.
	\end{cases}\]
Assume by symmetry that $\pi_\fr\le\pi_\rf$: then $2n\ep\ge np\ep+ n\pi_\fr=\De$, and in both the above cases we obtain $\mathbf{c}_\free\mathbf{p}_\free\le k^{O(n\ep)}$. Combining with \eqref{e:ap.id.ubd} then gives that
\[ \text{if }10\ep\le\be\text{ then } \E[\ZZ^2[\pi];\Om_B]
	\le \f{n^{O(1)}\,\E[\ZZ[\pi^1];\Om_B]}{2^{n\ep k/3 }}
	\le \f{n^{O(1)}\,\E[\ZZ[\pi^1];\Om_B]}{2^{\De k/6 }},\]
concluding the proof.
\epf
\elem

\bpf[Proof of Propn.~\ref{p:second.moment.id}]
Follows from Lem.~\ref{l:ap.near.id}.
\epf

\bpf[Proof of Thm.~\ref{t:second.moment}] 
Follows by combining Propns.~\ref{p:indep.maximizer}~and~\ref{p:second.moment.id}.
\epf

\section{Negative-definiteness of free energy Hessians}\label{s:nd}

In this section we prove Thm.~\ref{t:moments}.

\bppn\label{p:nd}
The Hessians $\hess\bPhi(\bhstar)$ and $\hess\bPhit(\bhtstar)$ are negative-definite.
\eppn

\subsection{Derivatives of the Bethe functional}
\label{ss:neg.def.derivs}

Let $\bh\in\simplexint$ with $\dbh$ and $\hbh$ both symmetric, and let $\bde$ be any signed measure on $\supp\vph$ (not necessarily symmetric) with $\bh+s\bde\in\simplexint$ for sufficiently small $\abs{s}$. Then
$$k\,\pd_s^2\bPhi(\bh+s\bde)|_{s=0}
=-k\,\angl{(\dbde/\dbh)^2}_{\dbh}
-d\,\angl{(\hbde/\hbh)^2}_{\hbh}
+dk\,\angl{(\vde/\vh)^2}_{\vh}
$$
where $a/b$ denotes the vector given by coordinate-wise division of $a$ by $b$, and $\angl{\cdot}_h$ denotes integration with respect to measure $h$, e.g.\ $\angl{(\vde/\vh)^2}_{\vh}=\sum_\si\bde(\si)^2/\vh(\si)$.

Given fixed marginals $\vde$, $\angl{(\dbde/\dbh)^2}_{\dbh}$ is minimized by $\dbde(\dsi)=\dbh(\dsi) \tf1d \sum_{i=1}^d \dchi_{\si_i}$ with $\dchi$ chosen to satisfy the margin constraint --- which, after a little algebra, becomes the vector equation
 $\vH^{-1}\vde = d^{-1}[I+(d-1)\dM]\dchi$ where $\vH\equiv\diag(\vh)$ and $\dM$ denotes the stochastic matrix
\beq\label{e:markov}
\dM_{\si,\si'}
\equiv
	\vh(\si)^{-1}
	\sum_{\dsi} \dbh(\dsi)
		\Ind{(\si_1,\si_2)=(\si,\si')}.
\eeq
If such $\dchi$ exists, then the minimal value of $\angl{(\dbde/\dbh)^2}_{\dbh}$ subject to marginals $\vde$ is $\ip{\vde}{\dchi}$. Define analogously the stochastic matrix $\hM$ corresponding to $\hbh$: if both
$\dL\equiv I+(d-1)\dM$
and
$\hL\equiv I+(k-1)\hM$
are
non-singular, then the maximum of $k\,\pd_\eta^2\bPhi(\bh+\eta\bde)|_{\eta=0}$ over all $\bde$ with marginal $\vde$ is given by
$$- dk \, \vde^t [ (\vH\dL)^{-1}+(\vH\hL)^{-1}-\vH^{-1} ]\vde
= -dk \, (\vH^{-1/2}\vde)^t F (\vH^{-1/2}\vde)$$
where
$F
\equiv ( \vH^{1/2}\dL\vH^{-1/2})^{-1}
	+( \vH^{1/2}\hL\vH^{-1/2})^{-1}-I$.
It is clear from \eqref{e:markov} that $\dM$ and $\hM$ are $\vh$-reversible, therefore $F$ is symmetric. Since $\sum_\si\vde(\si)=0$ we consider only the action $F'$ of $F$ on the space of vectors orthogonal to $\vh^{1/2}$. At a global maximizer we know $F'$ to be negative-semidefinite, so if $\det F\ne0$ then it is in fact negative-definite.
Thus let $\dM$, $\hM$, $\dM_2$, $\hM_2$ denote the Markov transition matrices corresponding (via \eqref{e:markov}) to $\dbhstar$, $\hbhstar$, $\dbhtstar$, $\hbhtstar$ respectively. In \S\ref{ss:neg.def.matrices} we will prove that the matrices
\beq\label{e:L.matrices}\begin{array}{rl}
\dL= I+(d-1)\dM,& \dL_2\equiv I+(d-1)\dM_2,\\
\hL= I+(k-1)\hM,&\hL_2\equiv I+(k-1)\hM_2,\\
L\equiv I-(d-1)(k-1)\dM\hM,&
L_2\equiv I-(d-1)(k-1)\dM_2\hM_2
\end{array}\eeq
are all non-singular.
Propn.~\ref{p:nd} then follows by noting that
$F = \vH^{1/2} \dL^{-1} L \hL^{-1}\vH^{-1/2}$.

\subsection{Calculation of transition matrices}
\label{ss:neg.def.matrices}

Recall the notation $\dq\equiv\dq^\star$, $\qf\equiv 1-\dq \asymp 2^{-k}$, and $\vr\equiv1-v \asymp 2^{-k}$. Recalling \eqref{e:bij} that $\vh^\star(\si)$ is proportional to $\hd^\star_\si\hh^\star_\si$, we record here that \vspace{-5pt}
\beq\label{e:vh}
\vh=
[2(1+\qf\vr)]^{-1}\times
\hspace{-8pt}
\kbordermatrix{
& \zf & \zz & \fz & \of & \oo & \fo & \ff\\
& \dq v & \dq \vr & 2 \qf \vr
	& \dq v & \dq \vr & 2 \qf \vr
	& \qf v
}
\eeq

\blem\label{l:dM}
The eigenvalues of $\dM$ counted with geometric multiplicity are
$$\eig(\dM)=
(1,1,1,\lm,\lm,-\lm,-\lm)
\quad\text{with }\lm\asymp \tf{1}{2^{3k/2}}.$$
The matrix $\dM_2$ is given by $\dM\otimes\dM$; consequently both
$\dL$ and $\dL_2$ are non-singular.

\bpf
The transition matrix $\dM\in\R^{7\times7}$ is block diagonal with blocks $\DOT{\mathfrak{m}}_\free$,
$\DOT{\mathfrak{m}}_\zro$,
$\DOT{\mathfrak{m}}_\one$
where 
$\DOT{\mathfrak{m}}_\free$ is the one-dimensional identity matrix (the action of $\dM$ on $\set{\ff}$),
and for $x=\zro,\one$ the matrix
$\DOT{\mathfrak{m}}_x\in\R^{3\times3}$
gives action of $\dM$ on
$\set{x\free,xx,\free x}$.
Recalling the definition \eqref{e:markov}, the entries of $\dM$ are straightforwardly calculated from \eqref{e:bp.zof.rf}, \eqref{e:grf.q}, and \eqref{e:bij}: for example,
$$
\dM_{\zz,\zz}
=\DSf{\gh_\rr^{}(\gh_\rr+\gh_\rf)^{d-2}}
	{ (\gh_\rr+\gh_\rf)^{d-1}-(\gh_\rf)^{d-1} }
= \smf{\vr}{1-v^{d-1}}
= \vr \smf{1+\qf}{1-\qf}$$
where the last step uses \eqref{e:q.rec}.
We therefore find $\DOT{\mathfrak{m}}_\zro=\DOT{\mathfrak{m}}_\one=\DOT{\mathfrak{m}}$ where\vspace{-5pt}
$$
\DOT{\mathfrak{m}}
\equiv
\text{\small$\kbordermatrix{
&\rf &\rr&\fr\\
\rf & 1-a
	& a(1-b)
	& ab\\
\rr & 1-a
	& a
	& 0\\
\fr & 1   & 0 & 0
}$}
\quad\text{with }
\left\{\hspace{-4pt}
\begin{array}{l}
a\equiv \vr(1+\qf)/(1-\qf)\asymp 2^{-k},
	\vspace{3pt}\\
b \equiv(2\vr\qf)/(v\dq)
	\asymp 4^{-k},
\end{array}\right.
$$
which has $\eig(\DOT{\mathfrak{m}})=(1,a b^{1/2},-a b^{1/2})$.
Thus $\eig(\dM)=(1,\eig(\DOT{\mathfrak{m}}),\eig(\DOT{\mathfrak{m}}))$
is as stated above. Since $\bhtstar=\bhstar\otimes\bhstar$, clearly $\dM_2=\dM\otimes\dM$, so the lemma is proved.
\epf
\elem

\blem
There exist (explicit) irreducible $\vh$-reversible transition matrices $\hM^\zro$, $\hM^\one$ such that
$\hM=\tf12[\hM^\zro+\hM^\one]$,
$\hM_2=\tf12[\hM^\zro\otimes\hM^\zro
	+\hM^\one\otimes\hM^\one]$.
The eigenvalues of $\hM^x$ counted with geometric multiplicity are
$\eig(\hM^x)=(1,\lm_2,\ldots,\lm_7)$
with $2^{k/2} |\lm_i|\lesssim 1$ for all $i>2$;
consequently both $\hL$ and $\hL_2$ are non-singular.

\bpf
Clearly $\hM_{\rr,\si}=\hM_{\fr,\si}=\Ind{\si=\rf}$, and it is straightforward to calculate that
\[\hM_{\ff,\ff}
=1-\hM_{\ff,\rf}=\de\equiv\qf\smf{1+\vr}{1-\vr}
	\asymp 2^{-k}.\]
The remaining entries of $\hM$ are easily determined by $\vh$-reversibility (see~\eqref{e:vh}): writing $\gm\equiv \vr/v \asymp 2^{-k}$ and $\ep\equiv 2\qf/\dq \asymp 2^{-k}$, we calculate
$$
\hM
=\kbordermatrix{
& \msg_\zro & \msg_\one & \ff \\
\msg_\zro
	& \hat{\mathfrak{m}}/2
	& \hat{\mathfrak{m}}/2 & \hat{\mathfrak{b}}\\
\msg_\one& \hat{\mathfrak{m}}/2 & \hat{\mathfrak{m}}/2 & \hat{\mathfrak{b}}\\
\ff&
	\hat{\mathfrak{a}}^t/2 & 
	\hat{\mathfrak{a}}^t/2 & \de
},\quad\text{where }
\bpm
\hat{\mathfrak{m}}&\hat{\mathfrak{b}}\\
\hat{\mathfrak{a}}^t&\de
\epm
\equiv\text{\small$\kbordermatrix{
	& \rf&\rr&\fr&\ff\\
\rf & 1-\de-\gm & \gm & \ep\gm & \de-\ep\gm\\
\rr & 1 & 0 & 0 & 0\\
\fr & 1 & 0 & 0 & 0\\
\ff & 1-\de & 0 & 0 & \de
}$}
$$
($\hat{\mathfrak{m}}\in\R^{3\times3}$ while $\hat{\mathfrak{a}},\hat{\mathfrak{b}}\in\R^{3\times1}$). Now consider decomposing $\hM=\tf{1}{2^k}\sum_{\ulit}\hM(\ulit)$ where $\hM(\ulit)$ is the transition matrix corresponding to $\hbh(\cdot\giv\ulit)$: there are only two possibilities $\hM^\zro,\hM^\one$ for $\hM(\ulit)$, depending on whether $\ulit_1\oplus\ulit_2=\zro,\one$. Since $\bhtstar(\kta) =\tf{1}{2^k}\sum_{\ulit}
	\hbhstar(\ksi^1\giv\ulit)
	\hbhstar(\ksi^2\giv\ulit)$, we conclude $\hM_2=\tf12[\hM^\zro\otimes\hM^\zro+\hM^\one\otimes\hM^\one]$.

The entries of each matrix
$\hM^x$ are easily read from $\hM$ except the ones giving the transition probabilities within $\set{\zf,\of}$. We calculate these from \eqref{e:zof} to find that
$$\hM^\zro
=\text{\small$
\kbordermatrix{
& \zf & \zz & \fz & \of & \oo & \fo & \ff\\
\zf & \tf{1+B}{2}\hat{\mathfrak{m}}_{\rf,\rf}
	& 0 & 0 & \tf{1-B}{2}\hat{\mathfrak{m}}_{\rf,\rf}
	& \gm
	& \ep\gm
	& \de-\ep\gm \\
\zz & 0 & 0 & 0 & 1 & 0 & 0 & 0 \\
\fz & 0 & 0 & 0 & 1 & 0 & 0 & 0 \\
\of & \tf{1-B}{2}\hat{\mathfrak{m}}_{\rf,\rf}
	& \gm
	& \ep\gm
	& \tf{1+B}{2}\hat{\mathfrak{m}}_{\rf,\rf}
	& 0 & 0 &
	\de-\ep\gm \\
\oo & 1 & 0 & 0 & 0 & 0 & 0 & 0 \\
\fo & 1 & 0 & 0 & 0 & 0 & 0 & 0 \\
\ff & \tf12 (1-\de) & 0 & 0 &
	\tf12 (1-\de)
	& 0 & 0 &
	\de \\
}$}$$
where $B= \vr/(1-\qf-2\vr-\qf\vr)$.
Then $\hM^\one$ is 
defined by exchanging the roles of $\zro$ and $\one$.

We see in particular that $\hM^\zro,\hM^\one$ is $\vh$-reversible,\footnote{To see this without explicit calculation of $\hM^x$, simply observe that (by symmetry) each $\hbh(\cdot\giv\ulit)$ has marginal $\vh$, and so from \eqref{e:markov} we have
$\vh(\si)\hM^x_{\si,\si'}
=\sum_{\ksi} \hbh(\ksi\giv L_1\oplus L_2=x)
=\vh(\si')\hM^x_{\si',\si}$.}
so, writing $\vH\equiv\diag(\vh)$, the matrix $\vH^{1/2}\hM^x \vH^{-1/2}$ is symmetric and hence orthogonally diagonalizable. Let $u$ be a left eigenvector of $\hM^x$ with eigenvalue $\lm$, such that $u$ has norm $1$ and is orthogonal to the constant vector $\vec 1$. 
Suppose $2^k|\lm|\gg 1$: the eigenvalue equations
\[
\begin{array}{l}
\lm u_{xx} = \gm \, u_{\neg x\free},\quad
	\lm u_{\free x} = \ep\gm \, u_{\neg x\free},\\
\lm u_\ff = (\de-\ep\gm)(u_\zf+u_\of) + \de u_\ff\\
\end{array}
\]
imply $|u_\si|\lesssim (2^k|\lm|)^{-1} \ll 1$ for all $\si\notin\rf$, so by the conditions $\|u\|=1$, $\ip{u}{\vec 1}=0$ we must have $\abs{u_\zf}\asymp\abs{u_\of}\asymp1$ with $|u_\zf+u_\of|\lesssim (2^k|\lm|)^{-1}$. The eigenvalue equation for $u_\zf$ then gives $|\lm|\asymp|\lm u_\zf| \lesssim  B + (2^k|\lm|)^{-1}$; rearranging and recalling $B\lesssim 2^{-k}$ then gives $2^k|\lm|^2 \lesssim |\lm|+1 \le2$ which proves the lemma.
\epf
\elem

\blem
The matrices $L$ and $L_2$ are non-singular.

\bpf
For $A\subseteq\msg$ write $\vh_A\equiv (\vh\I_A)/\vh(A)$,
the stationary distribution $\vh$ conditioned on $A$. The vectors
\[u_1\equiv\vh,\quad
u_2\equiv
(\vh_{\msg_\zro}+\vh_{\msg_\one})/2-\vh_\ff,\quad
u_3\equiv \vh_{\msg_\zro}-\vh_{\msg_\one}\]
are left eigenvectors of $\dM$ with eigenvalue $1$ such that $(w_i\equiv u_i/\vh^{1/2})_{i=1}^3$ forms an orthogonal basis for the $1$-eigenspace of the symmetrized matrix $\dS\equiv\vH^{1/2}\dM\vH^{-1/2}$. Define likewise $\hS\equiv\vH^{1/2}\hM\vH^{-1/2}$; if $w$ is orthogonal to this $1$-eigenspace, then Lem.~\ref{l:dM} implies $\nrm{w^t \dS\hS} =O( \nrm{w^t\dS} ) =O(2^{-3k/2}\nrm{w})$, so it remains to consider the action of $\dS\hS$ on the $1$-eigenspace of $\dS$. Clearly $u_1^t\hM=u_1^t$, and $u_3^t\hM=0$ by symmetry. Since $\vh_\ff$ is simply the indicator of $\ff$, clearly
$\nrm{w_2} \asymp \| \vh_\ff/\vh^{1/2} \|
\asymp \vh(\ff)^{-1/2} \asymp 2^{k/2}$;
and we calculate
$$u_2^t\dM\hM
= u_2^t\hM
= \bpm \slf14 & 0 & 0 & \slf14 & 0 & 0 & 0 \epm + O(2^{-k}),$$
therefore $\nrm{w_2^t\dS\hS}/\nrm{w_2}=\nrm{(u_2^t\dM\hM)/\vh^{1/2}}/\nrm{w_2}\asymp 2^{-k/2}$. It follows that $\dM\hM$ (equivalently $\dS\hS$) can have no eigenvalue with absolute value $\asymp 1/(dk)$, hence $L$ is non-singular.

For $1\le i,j\le 3$ let $w_{ij}\equiv w_i\otimes w_j$, and note that if $w$ is orthogonal to the span of $(w_{11},w_{12},w_{21})$ then $\nrm{w^t \dS_2\hS_2} = O(2^{-3k/2})\nrm{w}$. Next note that
$$w_{12}^t\dS_2\hS_2
=\f{(u_1^t \otimes u_2^t)
	[\hM^\zro\otimes\hM^\zro
	+\hM^\one\otimes\hM^\one]}{2\vh^{1/2}}
= \f{u_1^t \otimes (u_2^t\hM)}{\vh^{1/2}}
= w_1^t \otimes (w_2^t\hS),$$
so $\nrm{w_{12}^t\dS_2\hS_2}/\nrm{w_{12}}
=\nrm{w_2^t\hS}/\nrm{w_2} \asymp 2^{-k/2}$. Since $w_{12}\dS_2\hS_2$ and $w_{21}\dS_2\hS_2$ are orthogonal,
$$\f{\nrm{(aw_{12}+bw_{21})^t
	\dS_2\hS_2}^2}
	{\nrm{aw_{12}+bw_{21}}^2}
= \f{
	\abs{a}^2\nrm{w_{12}^t\dS_2\hS_2}^2
	+\abs{b}^2\nrm{w_{21}^t\dS_2\hS_2}^2}
	{(\abs{a}^2+\abs{b}^2) \nrm{w_{12}}^2}
\asymp 2^{-k/2}$$
for any $a,b\in\C$ (not both zero). It follows that $\dM_2\hM_2$ (equivalently $\dS_2\hS_2$) can have no eigenvalue with absolute value $\asymp 1/(dk)$, proving that $L_2$ is also non-singular.
\epf
\elem

\bpf[Proof of Propn.~\ref{p:nd}]
As shown in \S\ref{ss:neg.def.derivs} the result follows by verifying that the matrices defined in \eqref{e:L.matrices} are non-singular, which is done by the lemmas of \S\ref{ss:neg.def.matrices}.
\epf

\bpf[Proof of Thm.~\ref{t:moments}]
Recall (Defn.~\ref{d:simplex}) that $\ZZ$ is the sum of $\ZZ(\bh)$ over probability measures $\bh\equiv(\dbh,\hbh)$ on $\supp\phi$ such that $\bg\equiv(\dbg,\hbg)\equiv(n\dbh,m\hbh)$ is integer-valued, and lies in the kernel of matrix $\smash{H_{\simplex}\equiv\bpm \dH & -\hH \epm}$. Let $\starZZ$ denote the contribution to $\ZZ$ from (non-normalized) measures $\bg$ within euclidean distance $\smash{n^{1/2}\log n}$ of $\bgstar$. Thm.~\ref{t:first.moment.exponent} and Propn.~\ref{p:nd} together imply $\EZZ=[1+o(n^{-1})]\,\E[\starZZ]$. By Lem.~\ref{l:lin.bij}, the integer matrix $H_{\simplex}$ defines a surjection
$$L'\equiv\set{\bde\in\R^{\supp\phi}:\ip{\dbde}{1}=\ip{\hbde}{1}=0}
\quad\text{to}\quad
\set{\vde\in\R^\msg:\ip{\vde}{1}=0},$$
so $L\equiv L'\cap (\ker H_{\simplex})\cap\Z^{\supp\phi}$ is an $(\dbs+\hbs-\bar{s}-1)$-dimensional lattice with spacings $\asymp_k1$. The measures $\bg$ contributing to $\starZZ$ are given by the intersection of the euclidean ball $\smash{\set{\|\bg-\bgstar\|\le n^{1/2}\log n}}$  with an affine translation of $L$. The expansion \eqref{e:poly.correction} then shows that $\E[\starZZ]$ defines a convergent Riemann sum, therefore $\E[\starZZ]\asymp_k \exp\{ n\,\bPhi(\bhstar) \}$ as claimed.

In the pair partition function $\ZZ^2$, let $\indZZ$ denote the contribution from (non-normalized) measures $\bgt$ within euclidean distance $\smash{n^{1/2}\log n}$ of the independent-copies local maximizer $\bgtstar=\bgstar\otimes\bgstar$. Recall from the statement of Lem.~\ref{l:rule.out.intermediate}
that $\idZZ$
denotes the contribution to $\ZZ^2$
from the near-identical measures $\simplexid$. Decompose
\beq\label{e:pair.decompose}
\ZZ^2
=
\Big(\begin{array}{c}
	\text{near-independent}\\
	\text{contribution $\indZZ$}
	\end{array}\Big)
+ \Big(\begin{array}{c}
	\text{near-identical}\\
	\text{contribution $\idZZ$}
	\end{array}\Big)
+ (\text{remainder}),
\eeq
and note that for $\dlbd\le d\le \dubd$, Thm.~\ref{t:second.moment} and Propn.~\ref{p:nd} together imply that the expectation of the remainder is a negligible fraction of $\E[\ZZ^2]$:
\[\E[\ZZ^2]
=[1+o(n^{-1})]\, (\E[\indZZ] + \E[\idZZ]).\]
Repeating the argument above gives
$\E[\indZZ]\asymp_k
\exp\{ n\,(\bPhit(\bhtstar)) \}
\asymp_k (\EZZ)^2$, and combining with
Propn.~\ref{p:second.moment.id}
gives the conclusion
$\E[\ZZ^2]\lesssim_k(\EZZ)^2 + n^{O(1)}\,(\EZZ)$.
\epf

\section{From constant to high probability}\label{s:whp}

As in the proof of Thm.~\ref{t:moments}, $\starZZ$ denotes the contribution to the auxiliary model partition function on $(G,\ulit)$ from configurations whose non-normalized empirical measure $\bg\equiv(n\dbh,m\hbh)$ lies within euclidean distance $n^{1/2}\log n$ of $\bgstar$. The main result of this section is the following
\bppn\label{p:var}
For $k\ge k_0$ and $\dlbd\le d<d_\star$, $\Var\log(\starZZ+\ep\,\EZZ)\lesssim_k 1+o_\ep(1)$.
\eppn
The proposition easily implies the following strengthened statement of Thm.~\ref{t:whp}\ref{t:whp.frozen}:

\bthm\label{t:whp.ball}
For $k\ge k_0$ and $\dlbd\le d< d_\star$, $\lim_{\ep\downarrow0}\liminf_{n\to\infty}\P(\starZZ > \ep\,\EZZ)=1$.

\bpf
Write $\Lep\equiv\log(\starZZ+\ep\,\EZZ)$. For $\dlbd\le d<d_\star$, Thm.~\ref{t:moments} gives $\E[(\starZZ)^2]\asymp_{k,d}(\EZZ)^2$. Therefore there exist a constant $\de\equiv\de(k,d)>0$ for which $\P(\ZZ\ge \de\,\EZZ)\ge\de$, therefore $\E\Lep\ge \de\log\de+(1-\de)\log\ep+\log\EZZ$. On the other hand, $\starZZ\le\ep\,\EZZ$ if and only if $\Lep\le\log(2\ep\,\EZZ)$, which for small $\ep>0$ is much less than the lower bound on $\E\Lep$. Applying Chebychev's inequality and Propn.~\ref{p:var} therefore gives
\[ \P(\starZZ\le\ep\,\EZZ)
	\le\f{\Var\Lep}{( \E\Lep-\log(2\ep\,\EZZ) )^2}
	\lesssim_k \f{1+o_\ep(1)}{( \de\log(\de/\ep)-\log2 )^2}.\]
Taking $n\to\infty$ followed by $\ep\downarrow0$ proves the theorem.
\epf
\ethm

We prove Propn.~\ref{p:var} by controlling the increments of the Doob martingale of the random variable $\Lep\equiv\log(\starZZ+\ep\,\EZZ)$ with respect to the edge-revealing filtration $(\filt_i)_{1\le i\le m}$ for the graph $G_n\sim\Gndk$. We will show that the variance of $\Lep$ has two dominant components: the first is an ``independent-copies contribution'' coming from pair configurations with empirical measure near $\bhtstar\equiv\bhstar\otimes\bhstar$, which we will show in this section to be $\lesssim_k1+o_\ep(1)$. The other component is an ``identical-copies contribution'' coming from pair configurations with empirical measure near $\bhtzero$ or $\bhtone$, which we will easily see to be exponentially small in $n$ simply by the assumption $d<d_\star$. In \cite{dss-is} we demonstrate how to control the identical-copies contribution assuming only that the first moment is bounded below by a large constant.

\subsection{Doob martingale}

Consider forming the graph $G$ by beginning with $n$ vertices each incident to $d$ half-edges, and choosing, for $1\le i\le m$, a random set of $k$ unmatched half-edges to be joined into the $i$-th clause $a_i$ (every clause comes with literals). Let $\filt\equiv(\filt_i)_{0\le i\le m}$ denote the associated filtration; to prove Propn.~\ref{p:var} we will control the increments of the Doob martingale of $\Lep$ with respect to $\filt$:
\[\begin{array}{l}
\Var\Lep=\sum_{i=1}^m \Var_i\Lep\quad\text{where}\\
\qquad\qquad\Var_i\Lep
	\equiv \E[ (
	\E[\Lep\giv\filt_i]
	-\E[\Lep\giv\filt_{i-1}]
	)^2].\end{array}\]
Note the term $i=m$ is zero, since there is no randomness left when only two unmatched half-edges remain. Since a maximum of $k$ clauses will use any subset of the half-edges of size $k$, the random graph $G\giv\filt_i$ can be coupled with $G\giv\filt_{i-1}$ such that the graphs differ only in the placement of clauses on $k^2$ half-edges. Therefore
\beq\label{e:z.acute}
\begin{array}{l}
\Var_i\Lep \le\TS (\ep\,\EZZ)^{-2}
	\max_{\Ak,\Aprime}
	\E[|\Vep_i|^2]\quad\text{with} \\
\qquad\qquad
	\Vep_i \equiv(\ep\,\EZZ)\,
	\E_{i-1}[\log(\starZZ+\ep\,\EZZ)-\log(\starZZprime+\ep\,\EZZ)],
\end{array}
\eeq
where $\E_{i-1}$ is expectation conditioned on the graph $G^\circ$ with clauses $\smash{a_1,\ldots,a_{i-1}}$ (hence with $nd-k(i-1)$ unmatched variable-incident half-edges); and $\smash{\starZZ\equiv \starZZ(\Ak)}$, $\smash{\starZZprime\equiv \starZZ(\Aprime)}$ are the partition functions for coupled completions $G,\Gprime$ of $G^\circ$ to $(d,k)$-regular graphs: choose a random subset $\KK$ of $k^2 \wedge [(m-i+1)k]$ of the unmatched half-edges in $G^\circ$, and place on $\KK$ the clauses $\Ak$ for $G$, $\Aprime$ for $\Gprime$. Complete the graphs by placing random clauses $W$ on the set $\WW$ of remaining unmatched half-edges, using the same $W$ for both $G$ and $\Gprime$.\footnote{Though we suppress it from the notation, $\Vep_i$ should be regarded a function of $\Ak,\Aprime$, with $\E_{i-1}$ averaging over the possibilities of $W$.} By Jensen's inequality, the right-hand side of \eqref{e:z.acute} is increasing in $i$ for $i\le m-k$, therefore we may fix any threshold $m'\le m-k$ and bound
\beq
\label{e:jensen}\TS
(\ep\,\EZZ)^2 \Var\Lep
\le m' \max_{A,\Aprime}\E[|\Vep_{m'}|^2]
	+\sum_{i=m'+1}^{m-1} \max_{A,\acute A}\E[|\Vep_i|^2].
\eeq
The bound with $m'=m-k$ will suffice for the proof of Thm.~\ref{t:whp.ball}, but we will treat more generally $\E[|\Vep_i|^2]$ for $m-i\lesssim_k\log n$, for use in \cite{dss-is}.

\begin{figure}[ht]
\includegraphics[trim=.5in 1.1in .7in 1in,clip]{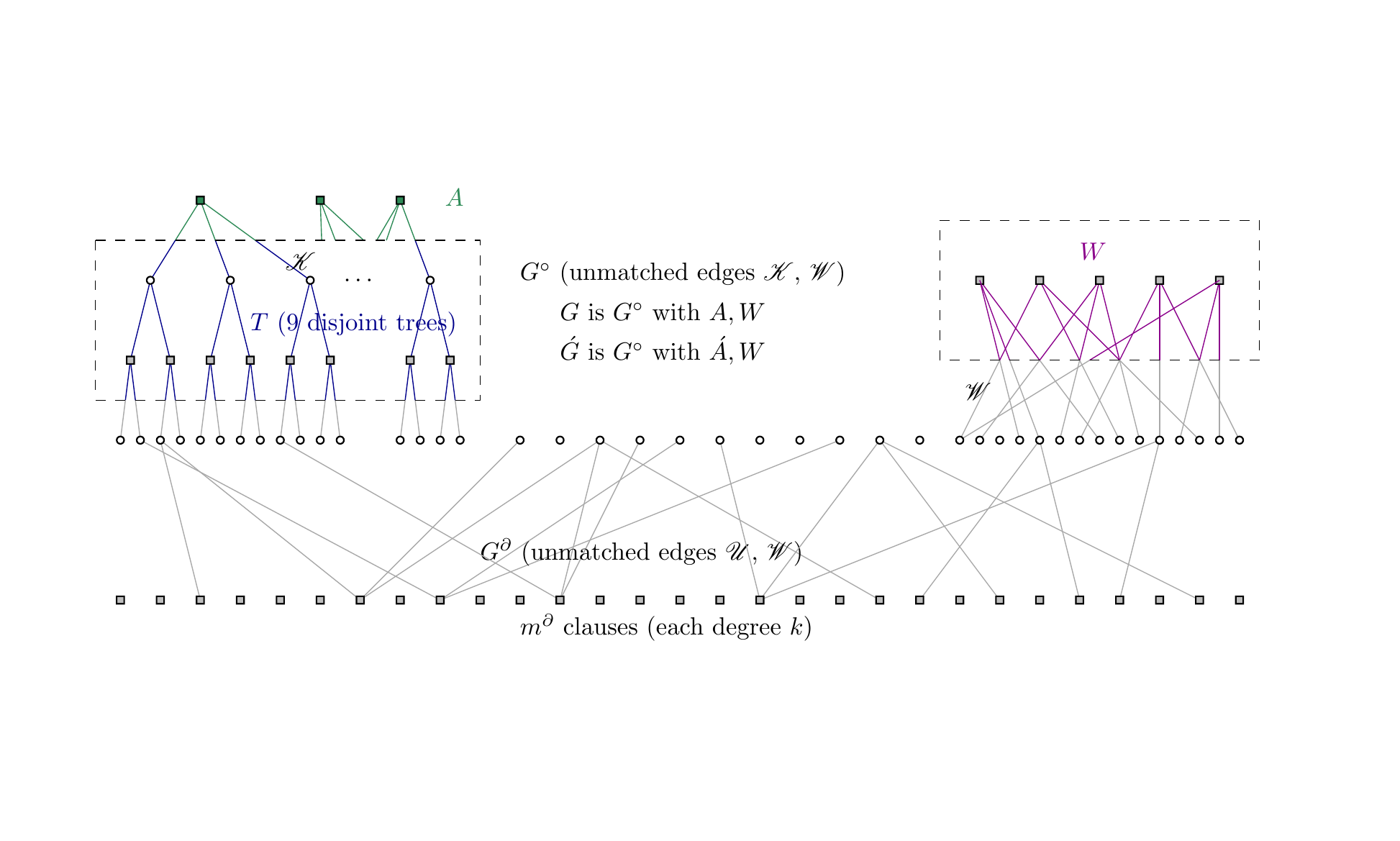}
\caption{$G^\circ\equiv$ graph with clauses $a_1,\ldots,a_{i-1}$;
$\Gmin \equiv G^\circ\setminus T$ (gray)}
\label{f:graph.fourier}
\end{figure}

Let $i$ be fixed, and note that since $(a-b)/a\le \log(a/b) \le (a-b)/b$ for any $a,b>0$,
$$
|\Vep_i|
\le |\Dep_i|\quad\text{for }
\Dep_i\equiv \E_{i-1}[\starZZ-\starZZprime].$$
Consider the graph $G^\circ$, with its unmatched half-edges partitioned into disjoint subsets $\KK$ and $\WW$. We shall define a certain local neighborhood $\TT$ of the half-edges in $\KK$, such that $\Gmin\equiv G^\circ\setminus\TT$ is a graph with unmatched half-edges in disjoint sets 
$\WW$ (as before)
and $\UU$ (leaves of $\TT$ without $\WW$); see Fig.~\ref{f:graph.fourier}.\footnote{Each unmatched half-edge is incident to a variable, and does not include a clause.} Then,
writing $\YY\equiv\UU\cup\WW$, we decompose
\[\begin{array}{rl}
\starZZ \hspace{-6pt}&=\sum_{\usi_\WW}\Psi_W(\usi_\WW)
	\sum_{\usi_\UU}
	\ka(\usi_\UU)\,
	\starZZpd[\usi_\UU,\usi_\WW],\\
\starZZprime \hspace{-6pt}&=\sum_{\usi_\WW}\Psi_W(\usi_\WW)
	\sum_{\usi_\UU}
	\kaprime(\usi_\UU)\,
	\starZZpd[\usi_\UU,\usi_\WW]
\end{array}
\]
where $\Psi_W$ is the partition function on $W$,
$\starZZpd[\usi_\YY]$ is the partition function on $\Gmin$ given boundary conditions $\usi_\YY$, $\ka(\usi_\UU)\equiv\ka_\TT(\usi_\UU|\usi_\WW)$ is the partition function on $\TT\cup\Ak$ given boundary conditions $\usi_\YY$, and $\smash{\kaprime(\usi_\UU)\equiv\kaprime_\TT(\usi_\UU|\usi_\WW)}$ is the partition function on $\TT\cup\Aprime$ given boundary conditions $\usi_\YY$. Averaging over $W$ and squaring gives
\[
|\Dep_i|^2
\le \sum_{\uta_\YY \equiv (\usi^1_\YY,\usi^2_\YY)}
	\psiw(\usi^1_\WW)
	\psiw(\usi^2_\WW)
	\underbrace{
	[\ka(\usi^1_\UU)-\kaprime(\usi^1_\UU)]
	}_{ \vpi(\usi^1_\UU) }
	\underbrace{
	[\ka(\usi^2_\UU)-\kaprime(\usi^2_\UU)]
	}_{ \vpi(\usi^2_\UU) }
	\tZZpd[\uta_\YY]
\]
where $\psiw$ denotes the average of $\Psi_W$ over the possibilities of $W$, and $\tZZpd$ denotes the \emph{pair} partition function on $\Gmin$. In the manner of \eqref{e:pair.decompose} we decompose $\tZZpd$ into a near-independent contribution $\indZZpd$, a near-identical contribution $\idZZpd$, and a remainder term which has expectation $o(n^{-1})\,\E[\ZZ^2]$. Substituting into the above gives the corresponding decomposition
\beq\label{e:Dep.decompose}
|\Dep_i|^2
\le\Big(\begin{array}{c}
	\text{near-independent}\\
	\text{contribution $\Deptind_i$}
	\end{array}\Big)
+\Big(\begin{array}{c}
	\text{near-identical}\\
	\text{contribution $\Deptid_i$}
	\end{array}\Big)
+(\text{remainder}).
\eeq
The remainder term has expectation $o(n^{-1})\,\E[\ZZ^2]$ uniformly over $\dlbd\le d\le \dubd$, and so can be ignored. The near-identical contribution $\Deptid_i$ has expectation $o(n^{-1})\,\E[\ZZ^2]$ uniformly over the integers $\dlbd\le d < \dubd$, and so can be ignored for the main result of this paper; however we will keep track of it for use in \cite{dss-is}.

Let $B^\circ_\ell(\KK)$ denote the ball of graph distance $\ell$ about $\KK$ in the graph $G^\circ$; the leaves of $B^\circ_\ell(\KK)$ are half-edges incident to variables ($\ell$ odd) or clauses ($\ell$ even). We shall fix a constant maximum depth $2t$ and set
\beq\label{e:TT}
\TT=B^\circ_{2t'}(\KK),\quad
t' = t\wedge \min\set{\ell: B^\circ_{2\ell} \text{ has fewer than $|\KK|$ connected components} }.
\eeq
$\WW$ can only intersect $\TT$ in its leaves; and we shall let $\UU$ denote the leaves of $\TT$ without $\WW$.
\beq\label{e:Dep.ind}
\TS
\E[\Deptind_i;\TT]
=\P(\TT)
\sum_{\uta_\WW}
\psiw(\usi^1_\WW)
\psiw(\usi^2_\WW)
\sum_{\uta_\UU}
\vpi(\usi^1_\UU)\vpi(\usi^1_\UU)
\ETT[ \indZZpd[\uta_\YY] ]
\eeq
where $\ETT$ denotes expectation conditioned on $\TT$. We shall soon see (Lem.~\ref{l:dangling} below) that in the graph $\Gmin$ the distribution of the boundary spins $\uta_\UU$ is very close to
$$\TS
\bp(\uta_\UU)\equiv
\bp(\usi^1_\UU)\bp(\usi^2_\UU)
\equiv
\prod_{u\in\UU}
	\hd_{\si^1_u}
	\hd_{\si^2_u},
$$
so we shall bound \eqref{e:Dep.ind} by projection onto a Fourier basis for $L^2(\msg^{2\UU},\bp)$: take $(\fourier_1,\ldots,\fourier_{|\msg|})$ to be an orthonormal basis for $L^2(\msg,\hd)$ with $\fourier_1\equiv1$. Then the functions
$\fourier_{\uss}(\usi_\UU)
\equiv\prod_{u\in\UU} \fourier_{s(u)}(\si_u)$
($\uss\in[|\msg|]^\UU$)
form an orthonormal basis for $L^2(\msg^\UU,\bp)$, and the functions $\fourier_{\uss^1,\uss^2}(\uta_\UU)\equiv\fourier_{\uss^1}(\usi^1_\UU)\fourier_{\uss^2}(\usi^1_\UU)$ form an orthonormal basis for $L^2(\msg^{2\UU},\bp)$. By Plancherel's identity,
\beq\label{e:plancherel}
\E[\Deptind_i;\TT]
		=\P(\TT)
		\sum_{\uta_\WW}
		\psiw(\usi^1_\WW)\,
		\psiw(\usi^2_\WW)
		\sum_{\uss^1,\uss^2}
		\vpi^\wedge_{\uss^1}\,
		\vpi^\wedge_{\uss^2}\,
		\F^\wedge_{\uss^1,\uss^2}
\eeq
where $^\wedge$ indicates the Fourier transform
with respect to the basis $\fourier$,
and
\beq\label{e:F.fourier}
\F(\uta_\UU)\equiv\F_\TT(\uta_\UU|\uta_\WW)
	\equiv
	\bp(\uta_\UU)^{-1}
	\,\ETT[ \indZZpd[\uta_\UU,\uta_\WW] ].
\eeq

\subsection{Expansion of partition function}\label{ss:z.expand}

On the graph $\Gmin$ we now analyze the partition function $\starZZpd[\usi_\YY]$ and its marginals
\beq\label{e:z.dangling.margin}
\begin{array}{rll}
\starZZpd[\cdot,\usi_\WW]
	\hspace{-6pt}&\equiv
	\sum_{\usi_\UU}\starZZpd[\usi_\UU,\usi_\WW]
	& \text{(marginal over $\UU$);}\\
\starZZpd
	\hspace{-6pt}&\equiv
	\sum_{\usi_\YY}\starZZpd[\usi_\YY]
	=\sum_{\usi_\WW} \starZZpd[\cdot,\usi_\WW]
	& \text{(marginal over $\YY=\UU\cup\WW$).}
\end{array}
\eeq
and likewise the pair partition function
$\indZZpd[\uta_\YY]$ and its marginals. Write $n^\pd$ and $m^\pd$ for the number of variables and clauses in $\Gmin$. We assume throughout that $|\WW|\lesssim_k\log n$ while $|T|\lesssim_{k,\ep} 1$. The main result of this subsection is Cor.~\ref{c:z.fourier} bounding the Fourier coefficients of the function $\F$ in \eqref{e:plancherel}.

It is more convenient here to work with the non-normalized tuple empirical measures $\bg\equiv(\dbg,\hbg)\equiv(n\dbh,m\hbh)$. The associated non-normalized marginal edge counts are given by $\dH\dbg$ and $\hH\hbg$ where $\dH$, $\hH$ are the marginalization matrices corresponding to $\dpsi$, $\hpsi$ as defined in \S\ref{ss:bethe}. The pair $\bg\equiv(\dbg,\hbg)$ can contribute to $\starZZpd[\usi_\YY]$ only if 
$$
\ip{\dbg}{1}=n^\pd,\quad
\ip{\hbg}{1}=m^\pd,\quad\text{and}\quad
\dH\dbg-\hH\hbg=\vH_{\usi_\YY}$$
where $\ip{\dbg}{1}$ indicates the total mass of $\dbg$, and $\vH_{\usi_\YY}$ denotes the non-normalized edge empirical measure associated to $\usi_\YY$. The contribution from such $\bg$ is given by
\[
\Xi(\bg)
\equiv
\smb{\ip{\dbg}{1}}{\dbg}
\smb{\ip{\hbg}{1}}{\hbg}
\f{(\dH\dbg)!}{ (nd)! }
\dpsi^{\dbg}
\hpsi^{\hbg}
\]
where we adopt the shorthand
$\dpsi^{\dbde}\equiv \prod_{\dsi}\dpsi(\dsi)^{\dbde(\dsi)}$, $\dbg!\equiv \prod_{\dsi} \dbg(\dsi)!$, etc.

The following lemma estimates the contribution in expectation from $\starZZpd[\usi_\YY]$ to $\starZZpd$.

\blem\label{l:dangling}
Suppose there exist measures $\bde^{\usi_\YY}$ (non-zero only on the support of $\phi$) satisfying $\hH\hbde^{\usi_\YY}=\dH\dbde^{\usi_\YY}+\vH_{\usi_\YY}$ for all $\usi_\YY$, with $(\ip{\dbde^{\usi_\YY}}{1},\ip{\hbde^{\usi_\YY}}{1})=(\nu,\mu)$ constant in $\usi_\YY$. Then
\beq
\label{e:dangling}
\f{\ETT[\starZZpd[\usi_\YY]]}
	{ \bp(\usi_\YY)\,\ETT[\starZZpd] }
= 1 + O_k\Big( \smf{\|\bde\|_1 \log n}{ n^{1/2} } \Big)
\eeq
where $\|\bde\|_1\equiv \sum_{\dsi} |\dbde(\dsi)|+\sum_{\ksi} |\hbde(\ksi)|$.

\bpf
$\Gmin$ has $n^\pd$ variables, $m^\pd$ clauses, and $|\YY|$ unmatched edges incident to variables, so $n^\pd d=m^\pd k+|\YY|$. Fix $\usi_\YY$ and abbreviate  $\bde\equiv\bde^{\usi_\YY}$; note from the assumptions that $\nu d=\mu k-|\YY|$. Thus $n'd=m'k$ for $n'\equiv n^\pd+\nu$ and $m'\equiv m^\pd+\mu$, so we may compare $\starZZpd$ with the partition function $\starZZ'$ on a full bipartite $(d,k)$-regular graph with $n'$ variables and $m'$ clauses. Away from the simplex boundary, empirical measures contributing to $\starZZpd[\usi_\YY]$ can be parametrized as $\bg-\bde$ where $\bg$ runs over the empirical measures contributing to $\starZZ'$. Writing $(a)_b$ for the falling factorial $a!/(a-b)!$, we have
\begin{align}\nonumber
\f{\Xi(\bg-\bde)}{\Xi(\bg)}
&= \DSf{ (\dbg)_{\dbde} (\hbg)_{\hbde} / (\dH\dbg)_{\dH\dbde} }
	{[(n')_\nu (m')_\mu / (n'd)_{\nu d}]
	\times 
	(\dpsi)^{\dbde} (\hpsi)^{\hbde} }\\ \nonumber
&= \underbrace{
	\f{1}{\bm{c}}
	\f{ (n')^\nu (m')^\mu / (n'd)^{\nu d} }
	{ (n')_\nu (m')_\mu / (n'd)_{\nu d} }
	}_{ 1/\wt{\bm{c}} }
\times
\underbrace{
	\f{ (\dbg)_{\dbde} (\hbg)_{\hbde} / (\dH\dbg)_{\dH\dbde} }
	{ (\dbgstar)^{\dbde} (\hbgstar)^{\hbde} / (\dH\dbgstar)^{\dH\dbde} }
	}_{ \bm{e}_{\bg,\bde} }
\times
\underbrace{
	\DSf{
	\bm{c}\cdot
	(\dbhstar)^{\dbde} (\hbhstar)^{\hbde} }
	{
	(\vh^\star)^{\dH\dbde} (\dpsi)^{\dbde} (\hpsi)^{\hbde}   }
	}_{\bm{r}_{\bg,\bde}}\\
\text{where }
	\bm{c}&
	\equiv
	\f{\dbz^\nu \hbz^\mu}{\barz^{\nu d}}
	=\Big(
		\f{ \dbz \hbz^{d/k} }
		{\barz^d} \Big)^\nu
	\hbz^{|\YY|/k} 
	= \exp\{ \nu\bPhistar\}
		\cdot \hbz^{|\YY|/k} 
		\quad\text{(using \eqref{e:bethe.explicit}).}
		\label{e:dangling.proportionality.constant}
\end{align}
The factor $\wt{\bm{c}}=\bm{c}[1+O_k(\|\bde\|^2_1/n)]$ is a proportionality constant not depending on $(\bg,\bde)$. For $\|\bg-\bgstar\|\le n^{1/2}\log n$, we find $\bm{e}_{\bg,\bde}\sim1$ while $\bm{r}_{\bg,\bde}$ gives the main dependence on $\usi_\YY$:
\begin{align*}
&
\bm{e}_{\bg,\bde}=1+ O_k\Big(\f{\|\bde\|_1\log n}{n^{1/2}}\Big),\\
\nonumber
&\bm{r}_{\bg,\bde}
=\f{(\dbz\dbhstar/\dpsi)^{\dbde}
	(\hbz\hbhstar/\hpsi)^{\hbde}}
	{ (\barz\vh^\star)^{\dH\dbde} }
= \prod_\si \hd_\si^{(\hH\hbde-\dH\dbde)(\si)}
= \bp(\usi_\YY)
\end{align*}
(using \eqref{e:bij} to calculate $\bm{r}$).
Recalling Propn.~\ref{p:nd} we have
(with $\bde\equiv\bde^{\usi_\YY}$)
$$
\f{\ETT[\starZZpd[\usi_\YY]]}
{\ETT[\starZZ']}
=[1+o(n^{-2})]
\f{
\sum_{\bg} \Ind{\|\bg-\bgstar\| \le n^{1/2}\log n}\,
	\Xi(\bg-\bde)
}
{\sum_{\bg} \Ind{\|\bg-\bgstar\| \le n^{1/2}\log n}\,
	\Xi(\bg)},$$
from which \eqref{e:dangling} easily follows.
\epf
\elem

Lem.~\ref{l:dangling} applies for any factor model with free energy attaining a local maximum at $\bgstar$ with negative-definite Hessian. In particular, it applies to both the first- and second-moment versions of the auxiliary model. We now show how to construct the required measures $\bde^{\uta_\YY}$ for the pair auxiliary model (the construction for the first-moment version being similar but simpler). The construction is based on the following
\blem\label{l:lin.bij}
For any $\tau,\tau'\in\msg^2$ there exists a signed integer measure $\bde=\bde_{\tau-\tau'}=(\dbde,\hbde)$ with $\supp\bde\subseteq\supp\psit$ such that
$$\ip{\dbde}{1}=0=\ip{\hbde}{1}\quad\text{and}\quad
\dH\dbde-\hH\hbde=\I_\tau-\I_{\tau'}.$$

\bpf
We shall show that in fact one can always find $\bde_{\tau-\tau'}$ of form $(0,\hbde)$. Define a graph on $\msg^2$ by placing an edge between $\tau,\tau'$ if and only if there exist $\kta,\kta'\in\supp\hpsit$ such that $\hH(\I_{\kta}-\I_{\kta'})=\I_\tau-\I_{\tau'}$; it suffices to show this graph is connected. It is easy to check that in the first-moment $\zof$ auxiliary model, $\msg$ is connected via differences $\hH(\I_{\ksi}-\I_{\ksi'})$ with $\ksi,\ksi'\in\supp\hpsi$. With $\neg\free\equiv\free$,
 if $\ksi\in\supp\hpsi$ then both $(\ksi,\ksi)$ and $(\ksi,\neg\ksi)$ belong to $\supp\hpsi$, and it follows that
$\De\equiv\set{ (\si,\si),(\si,\neg\si) }_{\si\in\msg}$ forms a connected subset of $\msg^2$. We will conclude by showing that any element of $\msg^2\setminus\De$ is connected to $\De$.

In a non-forcing clause we may freely set any spin $\tau$ to $\spint{\ff}{\ff}$, so it remains to consider those spins $\tau\in\msg^2\setminus\De$ which can only appear in forcing clauses:
\bnm[1.]
\item \emph{Clause forcing in both coordinates.} By considering a clause with literals identically $\zro$ we see that if $\ksi^1=(\oo,\zf^{k-1})$ and $\ksi^2=(\fo,\zf^{k-1})$ then $(\ksi^1,\ksi^2)\in\supp\hpsit$. Comparing with $(\ksi^1,\ksi^1)$ proves that $\spint{\oo}{\fo}$ is connected to $\spint{\oo}{\oo}\in\De$.
\item \emph{Clause forcing in one coordinate only.} For $\eta\in\set{\one,\free}$ and $\eta'\in\set{\zro,\one,\free}$ we have
$$
\dspint{\eta\one&\zf^{k-4}&\zf^3}{\eta'\free&\zf^{k-4}&\of^3},
\dspint{\ff&\zf^{k-4}&\zf^3}{\ff&\zf^{k-4}&\of^3}
\in\supp\hpsi
$$
(consider a clause with literals identically $\zro$ in the first case, and $(\zro^{k-1},\one)$ in the second). Taking the difference proves that $\spint{\eta\one}{\eta'\free}$ is connected to $\spint{\ff}{\ff}\in\De$.
\enm
The remaining cases follow by symmetry.
\epf
\elem

There exist $n_\TT,m_\TT\ge0$ bounded by the numbers of variables and clauses in $\TT$ such that $n'd=m'k$ for $n'=n^\pd+n_\TT$, $m'=m^\pd+m_\TT+k^{-1}|\WW|$. (In particular, we may certainly take $n'=n$ and $m'=m$.) Fix a reference spin $\tau^\mathrm{p}$, say $\tau^\mathrm{p}=\spint{\ff}{\ff}$. It is clear from Lem.~\ref{l:lin.bij} that we can find signed integer measures $\bde^\mathrm{v}\equiv(\dbde^\mathrm{v},\hbde^\mathrm{v})$ and
$\bde^\mathrm{c}\equiv(\dbde^\mathrm{c},\hbde^\mathrm{c})$
with $\ip{\hbde^\mathrm{v}}{1}=0=\ip{\dbde^\mathrm{c}}{1}$ and
$d^{-1}\dH\dbde^\mathrm{v}=\I_{\tau^\mathrm{p}}=k^{-1}\hH\hbde^\mathrm{c}$.\footnote{In the pair auxiliary model for {\sc nae-sat} we may trivially take $\bde^\mathrm{v}=(\I\dspint{\ff^d}{\ff^d},0)$ and $\bde^\mathrm{c}=(0,\I\dspint{\ff^k}{\ff^k})$; but note that for a general model the measures $\bde^\mathrm{v},\bde^\mathrm{c}$ can be constructed merely from the statement of Lem.~\ref{l:lin.bij}.}
Then set
\beq\label{e:delta.YY}
\TS
\bde^{\uta_\YY}
= n_T\bde^\mathrm{v}
+ (m_T+k^{-1}|\WW|)\bde^\mathrm{c}
+ \sum_{y\in\YY} \bde_{\tau_y-\tau^\mathrm{p}};\eeq
clearly this satisfies the conditions of Lem.~\ref{l:dangling}. Applying Lem.~\ref{l:dangling} with this choice of $\bde$ for $n'=n$ and $m'=m$ gives

\bcor\label{c:dangling.indep}
Let $\TT$ be as in \eqref{e:TT}.
\bnm[(a)]
\item\label{c:dangling.indep.a}
Recalling the notation of \eqref{e:z.dangling.margin}, we have
\[\begin{array}{rl}
\ETT[\indZZpd[\uta_\YY]]
\hspace{-6pt}&= \DS \bp(\uta_\YY)\,
\ETT[\indZZpd]\,
\Big( 1 + O_k(\smf{|\WW|\log n}{n^{1/2}} ) \Big)
	\vspace{2pt}\\
&\asymp_k\ETT[\starZZpd[\usi^1_\YY]] \cdot \ETT[\starZZpd[\usi^2_\YY]].
\end{array}
\]

\item\label{c:dangling.indep.b}
If $\bm{\Om}$ is any $T$-measurable event with $\P(\TT\in\bm{\Om})=o(n^{-1})$ then $\E[\Deptind_i;\TT\in\bm{\Om}]$
is at most $o(n^{-1})(\EZZ)^2$ where $\ZZ$ is the partition function of a full $(d,k)$-regular bipartite graph with $n$ variables and $m$ clauses.
\enm

\bpf
\eqref{c:dangling.indep.a}
Write $\nu\equiv n-n^\pd$, $\mu\equiv m-m^\pd$. The calculation of Lem.~\ref{l:dangling} applied to the auxiliary model gives
\beq\label{e:compare.full.graph}
\ETT[\starZZpd[\usi_\YY]]
=\bp(\uta_\YY)\,
\ETT[\indZZpd]\,
\Big( 1 + O_k\Big( \smf{|\YY|\log n}{n^{1/2}} \Big) \Big)
	\sim \bm{c}^{-1}\,\bp(\usi_\YY)\,\EZZ
\eeq
with $\ZZ$ the partition function on a full bipartite $(d,k)$-regular graph on $n$ variables, and with proportionality constant $\bm{c}\equiv \exp\{\nu\bPhistar\} \cdot \hbz^{|\YY|/k}$
where $\hbz$ is the normalizing constant for $\hbhstar$ as defined by \eqref{e:bij}.
The corresponding normalizing constant for $\hbhtstar$ is $\hbz^2$,
so we see that the proportionality constant corresponding to the pair auxiliary model is simply $\bm{c}^2$. Thus, applying \eqref{e:compare.full.graph} in the pair auxiliary model gives
\[\begin{array}{rl}
\E_T[\indZZpd[\uta_\YY]]
\hspace{-6pt}&\sim
\bm{c}^{-2}
	\,\bp(\uta_\YY)\,\E[\indZZ]
\asymp_k
(\bm{c}^{-1}\,\bp(\usi^1_\YY)\,\EZZ)
(\bm{c}^{-1}\,\bp(\usi^2_\YY)\,\EZZ)\\
&\sim
\E_T[\starZZpd[\usi^1_\YY]\,\E_T[\starZZpd[\usi^2_\YY].
\end{array}
\]

\smallskip\noindent
\eqref{c:dangling.indep.b} Substituting the result of \eqref{c:dangling.indep.a} into \eqref{e:Dep.ind} gives
\[\begin{array}{rl}
\ETT[\Deptind_i]
\hspace{-6pt}&\lesssim_k
	\big(\sum_{\usi_\YY}
		\psiw(\usi_\WW)
		(\ka+\kaprime)(\usi_\UU)
		\ETT[\starZZpd[\usi_\YY]] \big)^2\\
&\lesssim_{k,t}
	\big(
	\sum_{\usi_\UU}
	\sum_{\usi_\WW}
	\psiw(\usi_\WW)
	\ETT[\starZZpd[\usi_\UU,\usi_\WW]] \big)^2
\end{array}\]
The inner sum over $\usi_\WW$ is simply the expected partition function on the graph $G^\circ\setminus\TT$ with boundary condition $\usi_\UU$; applying \eqref{e:compare.full.graph} with $\WW=\emptyset$ we see that the expression above must be $\lesssim_{k,t}(\EZZ)^2$. In particular, if $\bm{\Om}$ is any event in the $\si$-algebra of $\TT$ with $\P(\TT\in\bm{\Om})=o(n^{-1})$ then
$\E[\Deptind_i;\TT\in\bm{\Om}]
\lesssim_{k,t}\P(\bm{\Om})(\EZZ)^2
= o(n^{-1})(\EZZ)^2$.
\epf
\ecor

The method of Lem.~\ref{l:dangling} can also be applied to estimate the dependence on $\uta_\UU$ fixing $\uta_\YY$: 

\blem\label{l:tree.expand}
Recalling the notation of \eqref{e:z.dangling.margin}, we have
\beq\label{e:fix.W}
\f{\ETT[\indZZpd[\uta_\YY]]}
	{ \bp(\uta_\UU)\,
\ETT[\indZZpd[\cdot,\uta_\WW]] }
=
 1 + O_k\Big( \smf{|T|\log n }{n^{1/2}} \Big)
\eeq
If $\TT$ is such that $k$ divides $|\UU|$, then there are coefficients $\xi$, $(\xi_j)_{j\le C}$ with $\|\xi\|_\infty\lesssim_k n^{-1/2}\log n$, $C\lesssim_k 1$ and $\|\xi_j\|_\infty\lesssim_k n^{-1/2}$ such that
\beq\label{e:part.fn.exp}
\f{\ETT[\indZZpd[\uta_\YY]]}
	{ \bp(\uta_\UU)\,
\ETT[\indZZpd[\cdot,\uta_\WW]] }
=1+\ip{\vH_{\uta_\UU}}{\xi}
	+\sum_{j=0}^{C}
		\ip{\vH_{\uta_\UU}}{\xi_j}^2
	+O_k\Big( \smf{(\log n)^{6}}{n^{3/2}}\Big)
\eeq

\bpf
Let $\bg$ denote a non-normalized empirical measure contributing to the partition function of a bipartite $(d,k)$-regular graph with $n'$ variables and $m'-k^{-1}|\WW|$ clauses, subject to boundary configuration $\uta_\WW$ on the unmatched half-edges $\WW$. Away from the simplex boundary, empirical measures contributing to $\tZZpd[\uta_\YY]$ can be parametrized as $\bg-\bde^{\uta_\UU}$ where
$
\bde^{\uta_\UU}
\equiv n_T \bde^\mathrm{v} + m_T\bde^\mathrm{c}
	+\sum_{u\in\UU}\bde_{\tau_u-\tau^\mathrm{p}}
$
(cf.~\eqref{e:delta.YY}). Comparing $\Xi(\bg)$ with $\Xi(\bg-\bde)$ gives \eqref{e:fix.W}.

In the special case that $|\UU|$ is divisible by $k$, we may simply take $n_\TT=0$ and $m_\TT=|\UU|/k$, so $\bde$ can be expressed as a linear function of the $\UU$-marginal:
$$\TS
\bde^{\uta_\UU}=\sum_\tau \vH_{\uta_\UU}(\tau)\,\bde^\tau\quad\text{where}\quad
\bde^\tau
\equiv(\dbde^\tau,\hbde^\tau)\equiv
\bde_{\tau-\tau^\mathrm{p}} + k^{-1}\bde^\mathrm{c}.$$
Assume for simplicity that $\dbde^\tau=0$ for all $\tau$:
writing $\hba\equiv(\hbg-\hbg^\star)/\hbg^\star$, we estimate
\[\begin{array}{rl}
\smf{(\hbg)_{\hbde}}{(\hbg^\star)^{\hbde}}
\hspace{-6pt}&=
	1
	+\ip{\hbde}{\hbA}
	+\smf{(\ip{\hbde}{\hbA})^2}{2}
	-\ip{\hbde^2}{\hbB}
	+ O_k\Big(
		\smf{(\log n)^6}{n^{3/2}} \Big)
	\vspace{4pt}\\
&\qquad
\text{for }
\hbA\equiv \hba+\tf12\hba^2+(2\hbg^\star)^{-1}, \
\hbB\equiv (2\hbg^\star)^{-1}.
\end{array}\]
We then simplify $\ip{\hbde^{\uta_\UU}}{\hbA}=\ip{\vH_{\uta_\UU}}{\xi}$ and $\hbde^{\uta_\UU}(\kta)^2 \hbB(\kta)=\ip{\vH_{\uta_\UU}}{\xi_{\kta}}^2$ where $\xi(\tau)\equiv \ip{\hbde^\tau}{\hbA}$ and $\xi_{\kta}(\tau)\equiv \hbB^{1/2}(\kta) \hbde^\tau(\kta)$. Relabelling gives \eqref{e:part.fn.exp} in the case that $\dbde^\tau=0$ for all $\tau$. The result in general follows by an easy generalization of the above calculation.
\epf
\elem

Lem.~\ref{l:tree.expand} easily implies bounds on the Fourier coefficients of the function $\F(\uta_\UU)$ of \eqref{e:plancherel}: for $\uss\in[|\msg|]^\UU$ write $|\uss|\equiv|\set{u\in\UU:s_u\ne1}|$, and write $\emptyset$ for the identically-$1$ vector in $[|\msg|]^\UU$. For $\uss^1,\uss^2\in[|\msg|]^{\UU}$ write $|\uss^1\uss^2|\equiv|\set{u\in\UU: s^1_us^2_u\ne1}|$. We then have the following

\bcor\label{c:z.fourier}
For any realization $\TT$ of \eqref{e:TT}
and $\F(\uta_\UU)\equiv\F_\TT(\uta_\UU|\uta_\WW)$ as in \eqref{e:plancherel}, 
$$
\mathbf{R}_{\uss^1,\uss^2}
\equiv
\mathbf{R}_\TT(\uss^1,\uss^2|\uta_\WW)\equiv
	(\ETT[\indZZpd[\cdot,\uta_\WW]])^{-1}
	\F^\wedge_{\uss^1,\uss^2} 
	$$
satisfies
$\mathbf{R}_{\emptyset,\emptyset}\asymp1$
and
$|\mathbf{R}_{\uss^1,\uss^2}|\lesssim_{k,t} n^{-1/2}\log n$ for $|\uss^1\uss^2|\ge1$, consequently
\beq\label{e:fourier.zero}\TS
(\EZZ)^2
\asymp_k
\sum_{\uta_\WW}
\psiw(\usi^1_\WW)\psiw(\usi^2_\WW)
(\ka^\wedge_\emptyset)^2
\F^\wedge_{\emptyset,\emptyset}.
\eeq
If further $\TT$ is disjoint from $\WW$ with $|\UU|$ divisible by $k$, then $|\mathbf{R}_{\uss^1,\uss^2}|\lesssim_k n^{-1}$ for $|\uss^1\uss^2|\ge2$, and $|\mathbf{R}_{\uss^1,\uss^2}|\lesssim_{k,t} n^{-3/2}(\log n)^6$ for $|\uss^1\uss^2|\ge3$. All estimates hold uniformly in $\uta_\WW$.

\bpf
The estimates on $\mathbf{R}$ follow straightforwardly from Lem.~\ref{l:tree.expand}. To see \eqref{e:fourier.zero}, we calculate
\begin{align*}
(\EZZ)^2
&\asymp_k
\E\Big[
\sum_{\uta_\WW}\psiw(\usi^1_\WW)\psiw(\usi^2_\WW)
\sum_{\uss^1,\uss^2}
\ka^\wedge_{\uss^1}
\ka^\wedge_{\uss^2}
\F^\wedge_{\uss^1,\uss^2}\Big]\\
&=
\Big[
	1 + O_{k,t}\Big( \smf{\log n}{n^{1/2}}\Big)
\Big]
\,\E
\Big[
(\ka^\wedge_\emptyset)^2
\sum_{\uta_\WW}\psiw(\usi^1_\WW)\psiw(\usi^2_\WW)
\F^\wedge_{\emptyset,\emptyset}\Big]
\end{align*}
where the first step uses 
Cor.~\ref{c:dangling.indep}\ref{c:dangling.indep.a} and the second step uses \eqref{e:fix.W} together with the trivial bound $|\ka^\wedge_{\uss}|\le \|\fourier_{\uss}\|_\infty \ka^\wedge_\emptyset \lesssim_{k,t} \ka^\wedge_\emptyset$, since the orthogonality relation $\|\fourier_s\|_{\dq}=1$ implies $\|\fourier_{\uss}\|_\infty \le (\min_\si \hd_\si)^{-|\uss|/2} \lesssim_{k,t} 1$.
\epf
\ecor

\subsection{Local neighborhood Fourier coefficients}

Let us recall again the definition \eqref{e:TT} of the local neighborhood $\TT$ of the unmatched $\KK$ in $G^\circ$ (equipped with random literals), with leaves joining $\TT$ to the graph $\Gmin$ considered in \S\ref{ss:z.expand}. Recall also that $\Ak$, $\Aprime$ are arbitrary choices of clauses (with literals) to place on $\KK$, and we write $\ka(\usi_\UU)$ (resp.\ $\kaprime(\usi_\UU)$) for the partition function on $\TT\cup\Ak$ (resp.\ $\TT\cup\Aprime$) given boundary configuration $\usi_\YY$. In this subsection we control the Fourier coefficients $\vpi^\wedge=(\ka-\kaprime)^\wedge$ appearing in \eqref{e:plancherel}.

Let $\mathbf{T}$ denote the event that $\TT$ consists of $|\KK|$ tree components with $\TT\cap\WW=\emptyset$. Let $\mathbf{C}^\circ$ denote the event that $\TT$ either contains a single cycle or has a single intersection with $\WW$ (but not both), but still consists of $|\KK|$ components.

\blem\label{l:tree.sym}
For $\TT\in\mathbf{T}$, $\vpi^\wedge_{\uss}=0$ for all $|\uss|\le1$; and $\ka^\wedge_\emptyset|_{\mathbf{T}}$ takes a constant value $\emptyka$, which does not depend on the literals on $\TT$ or on the clauses $\Ak$. For $\TT\in\mathbf{C}^\circ$, $\vpi^\wedge_\emptyset=0$.

\bpf
On the event $\mathbf{T}\cup\mathbf{C}^\circ$, the graphs $\TT\cup\Ak$ and $\TT\cup\Aprime$ are isomorphic ignoring the literals. If $|\uss|\le1$ then $\fourier_{\uss}$ depends at most on the spin of a single edge $e\in\UU$. For $\TT\in\mathbf{T}$, using the symmetry of \acr{nae-sat} one can produce an involution $\iota:\usi_\UU\mapsto\acute\usi_\UU$ on $\msg^\UU$ which keeps $\si_e$ fixed, is measure-preserving with respect to $\bp$, and satisfies $\ka(\usi_\UU)=\kaprime(\acute\usi_\UU)$: set $\acute\si_u$ to be $\si_u$ or $\neg\si_u$ depending on whether the sum of literals along the unique path joining $e$ to $u$ in $\TT\cup A$ differs in parity from the corresponding sum in $\TT\cup\Ak$. Then
$$\TS
\ka^\wedge_{\uss}
=\sum_{\usi_\UU} \bp(\iota\usi_\UU)
	\fourier_{\uss}(\iota\usi_\UU) \ka(\iota\usi_\UU)
=\sum_{\usi_\UU} \bp(\usi_\UU)
	\fourier_{\uss}(\usi_\UU) \kaprime(\usi_\UU)
=\kaprime^\wedge_{\uss},
$$
proving our claim on $\mathbf{T}$. A similar argument proves $\vpi^\wedge_\emptyset=0$ on $\mathbf{C}^\circ$.
\epf
\elem

\blem\label{l:tree.cov}
For $\TT\in\mathbf{T}$,
$\vpi^\wedge_{\uss}
\lesssim_k \emptyka/[ 4^{(k-4)t} ]$
for all $|\uss|=2$.

\bpf
Since $|\uss|=2$ we may write $\fourier_{\uss}(\usi_\UU)=f(\si_u)g(\si_w)$ for $f\equiv\fourier_{s(u)}$ and $g\equiv\fourier_{s(w)}$. If $u,w$ belong in the same connected component of $\TT$, arguing as in the proof of Lem.~\ref{l:tree.sym} gives $\vpi^\wedge_{\uss}=0$, so assume they belong to different components. Since $\TT\in\mathbf{T}$ we may define the random measure $\mu_\TT(\usi_\UU)\equiv\bp(\usi_\UU)\ka(\usi_\UU)/\emptyka$, and similarly $\muprime_\TT$. Then
$$\vpi^\wedge_{\uss}
=\emptyka \ip{ \mu_\TT-\muprime_\TT}{fg}
=\emptyka [ \Cov_\mu(f,g)-\Cov_{\muprime}(f,g) ],
$$
since Lem.~\ref{l:tree.sym} implies $f$ has the same expectation with respect to $\mu$ or $\muprime$ (and likewise $g$). Let $\gm_u$ denote the (unique) path joining $u$ to $\KK$ in $\TT$, and likewise $\gm_w$. Let $N$ denote the event that on the path $\gm\equiv\gm_u\cup\gm_w$ there exists a clause $a$ such that, among the $k-2$ variables in $(\pd a)\setminus\gm$, there exist two variables $v',v''$ with
$$
\si_{v'\to a}=\si_{v''\to a}=\free\quad\text{or}\quad
\lit_{a v'}\oplus\si_{v'\to a}
=\neg\lit_{a v''}\oplus\si_{v''\to a}.
$$
Note that $N$ is $(\TT,\usi_{\UU\setminus\set{u,w}})$-measurable, and on this event $\Cov_\mu(f,g)=\Cov_{\muprime}(f,g)=0$. For \emph{any} fixed realization of literals on $\TT$ the probability (with respect to the law $\bp$ of the spins on $\UU\setminus\set{u,w}$) that $N$ fails is $\lesssim 4^{-(k-4)t}$ simply by the $\zro/\one$ symmetry in $\bp$, so we find
$$\TS|\vpi^\wedge_{\uss}|
\lesssim 4^{-(k-4)t}\,\emptyka
\max_{s\in[|\msg|]}\|\fourier_s\|_\infty^2
\lesssim_k 4^{-(k-4)t}\,\emptyka,$$
since the relation $\|\fourier_s\|^2_{\dq}=1$ implies $\|\fourier_s\|_\infty^2\le (\min_{\si\in\msg} \dq_\si)^{-1} \lesssim_k1$.
\epf
\elem

Consider \eqref{e:plancherel} for $\TT\in\mathbf{T}$: by Lem.~\ref{l:tree.sym} there is no contribution from terms $|\uss^i|\le1$. The number of pairs $(\uss^1,\uss^2)$ with $|\uss^1|=|\uss^2|=|\uss^1\uss^2|=2$ is at most $[|\msg||\KK| (dk)^t]^2 \lesssim_k (k^5 4^k)^t$, and combining Cor.~\ref{c:z.fourier} and Lem.~\ref{l:tree.cov} we see that the dominant contribution to \eqref{e:plancherel} comes from these pairs:
for $\TT\in\mathbf{T}$, 
\beq\label{e:tree.bound}
\E[\Deptind_i;\TT] \lesssim_k
	\smf{1 + o_t(1)}{n}
	\smf{(k^5 4^k)^t}{16^{(k-4)t}}
	(\emptyka)^2
	\sum_{\uta_\WW}\psiw(\usi^1_\WW)\psiw(\usi^2_\WW)
	\ETT[\indZZpd[\cdot,\uta_\WW]]
\lesssim_k 
	\smf{(k^6 4^{-k})^t}{n} (\EZZ)^2
\eeq
where $o_t(1)$ indicates the contribution from pairs $(\uss^1,\uss^2)$ with $|\uss^1\uss^2|\ge3$, and the last bound uses \eqref{e:fourier.zero}.
Similarly, 
$\P(\mathbf{C}^\circ)\lesssim_k n^{-1}[\log n+O_t(1)]$,
and by Lem.~\ref{l:tree.cov} the terms $\uss^i=\emptyset$ in \eqref{e:plancherel} vanish on $\mathbf{C}^\circ$,
 so combining with Cor.~\ref{c:z.fourier} gives
\beq\label{e:bound.interior.cycle}
\E[\Deptind_i;\mathbf{C}^\circ]
\lesssim_{k,t}
(\EZZ)^2 \,n^{-3/2}(\log n)^2
\eeq
(applying \eqref{e:fourier.zero} as above).

Lastly we consider the event $\mathbf{C}_{t'}$ that $B^\circ_{2t'-2}(\KK)$ has $|\KK|$ components, but $\TT=B^\circ_{2t'}(\KK)$ has $|\KK|-1$ components (cf.~\eqref{e:TT}) and is disjoint from $\WW$. On this event we again decompose into contributions from the different local maxima, but differently than in \eqref{e:Dep.decompose}:
\begin{align}\nonumber
\smf{|\Vep_i|^2}{3\,(\ep\,\EZZ)^2}
&\le
\Big(\E_{i-1}
	\Big[\log\smf{\starZZ+\ep\EZZ}
	{\starZZ_\emptyset+\ep\EZZ}\Big]
	\Big)^2
+\Big(\E_{i-1}
	\Big[\log\smf{\starZZprime_\emptyset+\ep\EZZ}
	{\starZZprime+\ep\EZZ}\Big]\Big)^2
+\Big(\E_{i-1}
	\Big[\log\smf{\starZZ_\emptyset+\ep\EZZ}
	{\starZZprime_\emptyset+\ep\EZZ}\Big]\Big)^2\\
&\le
\Big( \smf{\E_{i-1}[\starZZ-\starZZ_\emptyset]}{\ep\,\EZZ} \Big)^2
+\Big(
	\smf{\E_{i-1}[\starZZprime-\starZZprime_\emptyset]}{\ep\,\EZZ}
	\Big)^2
+\Big(\E_{i-1}
	\Big[\smf{\starZZ_\emptyset-\starZZprime_\emptyset}
	{\starZZ_\emptyset\wedge \starZZprime_\emptyset}
	\Big]\Big)^2\label{e:Vep.triangle}
\end{align}
where $\starZZ_\emptyset$ is the leading term in the Fourier expansion of $\starZZ$:
\[\begin{array}{l}
\starZZ_\emptyset
\equiv\sum_{\usi_\WW}\Psi_W(\usi_\WW)
	\FF^\wedge_\emptyset
	\ka^\wedge_\emptyset
\quad\text{where }\\
\qquad\qquad\FF(\usi_\UU)
\equiv
\FF(\usi_\UU|\usi_\UU)
\equiv\bp(\usi_\UU)^{-1}\,
\starZZpd[\usi_\UU,\usi_\WW].
\end{array}\]
Then the last term in \eqref{e:Vep.triangle} simplifies to $|\vpi^\wedge_\emptyset|^2/[\ka^\wedge_\emptyset\wedge \ka^\wedge_\emptyset]^2$. We decompose the other two terms in the manner of \eqref{e:pair.decompose} and \eqref{e:Dep.decompose}: write $\Uep_i \equiv \E_{i-1}[\starZZ-\starZZ_\emptyset]$ and $\Uprimeep_i \equiv \E_{i-1}[\starZZprime-\starZZprime_\emptyset]$; then for example the near-independent contribution to $|\Uep_i|^2$ is
\[ \indUep_i
=\sum_{\uta_\YY}\psiw(\usi^1_\WW)\psiw(\usi^2_\WW)
	\indZZpd[\uta_\YY]
	(\ka(\usi^1_\UU)-\ka^\wedge_\emptyset)
	(\ka(\usi^2_\UU)-\ka^\wedge_\emptyset).
\]
Then in place of \eqref{e:Dep.decompose} we have
\beq\label{e:Vep.decomp}
\tf13|\Vep_i|^2
\le
\Big( \smf{\vpi^\wedge_\emptyset\,(\ep\,\EZZ)}
	{\ka^\wedge_\emptyset\wedge \ka^\wedge_\emptyset} \Big)^2
	+\Veptind_i+\Veptid_i
	\quad\text{with }
	\begin{cases}
	\Veptind_i\equiv\indUep_i+\indUprimeep_i
		&\text{(near-independent)},\\
	\Veptid_i\equiv\idUep_i+\idUprimeep_i
		&\text{(near-identical)}.
	\end{cases}
\eeq
It follows straightforwardly from Cor.~\ref{c:z.fourier} that
$\E[ \indUep_i ]\lesssim_{k,t} (\EZZ)^2 n^{-1/2} \log n$, and taking into account $\P(\mathbf{C}_{t'})\lesssim_{k,t} n^{-1}$ gives
\beq\label{e:fourier.diff.small.cycle}
\E[\Veptind_i;\mathbf{C}_{t'}]\lesssim_{k,t} 
	\ep^{-2} n^{-3/2}(\log n).
\eeq

\blem\label{l:small.cycle}
If $\TT\in\mathbf{C}_{t'}$ then
$\ka^\wedge_\emptyset\asymp\kaprime^\wedge_\emptyset$
and
$|\vpi^\wedge_\emptyset|\lesssim_k
	4^{-(k-4)t'}
	[\ka^\wedge_\emptyset
	\wedge
	\kaprime^\wedge_\emptyset]$.

\bpf
Let $\kacirc(\usi_\UU,\usi_\KK)$ denote the partition function on $\TT$ given boundary conditions $(\usi_\UU,\usi_\KK)$, and define the random measure
$$\TS
\mucirc(\usi_\KK)\equiv
	[(\kacirc)^\wedge_\emptyset]^{-1}
	\sum_{\usi_\UU} \bp(\usi_\UU)
		\kacirc(\usi_\UU,\usi_\KK).$$
If $\iota(\usi_\KK)$ is the indicator that $\usi_\KK$ is valid for clauses $\Ak$ (likewise $\ioprime$ for clauses $\Aprime$) then
$$\TS
\vpi^\wedge_\emptyset
=\sum_{\usi_\UU,\usi_\KK}
	\bp(\usi_\UU)
	\kacirc(\usi_\UU,\usi_\KK)
	[(\iota-\ioprime)(\usi_\KK)]
=(\kacirc)^\wedge_\emptyset
	\sum_{\usi_\KK}
	\mucirc(\usi_\KK)
	[(\iota-\ioprime)(\usi_\KK)].
$$
By definition of the event $\mathbf{C}_{t'}$,
$$\TS
\mucirc(\usi_\KK)
= \mucirc(\si_{e'},\si_{e''})
	\prod_{e\in\KK\setminus\set{e',e''}} \mucirc(\si_e)$$
where $\set{e',e''}$ is the unique pair of edges in $\KK$ such that $B^\circ_{2t'}(e')$ and $B^\circ_{2t'}(e'')$ intersect. The graph $\TT$ (without $\Ak$ or $\Ak'$) contains no cycles, so it follows from the symmetry argument of Lem.~\ref{l:tree.sym} that the marginal of $\mucirc$ on each $e\in\KK$ does not depend on the literals on $\TT$; further each marginal must simply be $\bp$ from the Bethe recursions. It remains to note (arguing as in the proof of Lem.~\ref{l:tree.cov}) that $|\mucirc(\si_{e'},\si_{e''})-\mucirc(\si_{e'})\mucirc(\si_{e''})| \lesssim 4^{-(k-4)t'}$, from which we conclude
$$\TS
[1+O_k(4^{-(k-4)t'})]\ka^\wedge_\emptyset
=(\kacirc)^\wedge_\emptyset
	\sum_{\usi_\KK}\bp(\usi_\KK)\iota(\usi_\KK)
=
[1+O_k(4^{-(k-4)t'})]\kaprime^\wedge_\emptyset$$
(where $\ip{\bp}{\iota}=\ip{\bp}{\ioprime}$ again by symmetry).
\epf
\elem

\bpf[Proof of Propn.~\ref{p:var}]
Writing $\mathbf{C}\equiv\bigcup_{t'\le t}\mathbf{C}_{t'}$, we have
\begin{align*}
\E[|\Vep_i|^2]
&\le
\E[|\Vep_i|^2;\mathbf{C}]
+
\E[|\Dep_i|^2;\mathbf{C}^c]\\
&\lesssim_k
		\underbrace{n^{-1}
		[1 +(k^64^{-k})^t/\ep^2 + o_t(1)]
		(\ep\,\EZZ)^2}_{\text{independent-copies contribution}}
+\underbrace{\E[\Veptid_i;\mathbf{C}]+\E[\Deptid_i;\mathbf{C}^c]}_{\text{identical-copies contribution}}
\end{align*}
where the leading (order-$n^{-1}$) contribution comes from the bound Lem.~\ref{l:small.cycle} on the first term in the decomposition \eqref{e:Vep.decomp} of $\E[|\Vep_i|^2;\mathbf{C}]$, together with the bound \eqref{e:tree.bound} on $\E[\Deptind_i;\mathbf{T}]$. The $o_{k,t}(n^{-1})$ error term combines the bounds \eqref{e:bound.interior.cycle} and \eqref{e:fourier.diff.small.cycle} together with the result of Cor.~\ref{c:dangling.indep}\ref{c:dangling.indep.b} which gives that the variance contribution from $\bm{\Om}=(\mathbf{T}\cup\mathbf{C}^\circ\cup\mathbf{C})^c$ is $o_{k,t}(n^{-1})$. Since we have assumed $\bPhistar>0$ the remaining contribution, from $\bhtzero$ and $\bhtone$, is exponentially small in $n$, so the proof is concluded by choosing $t= t(k,\ep)$ such that $(k^64^{-k})^t\le\ep^2$.
\epf

\subsection{Variance bound for general factor models}

We summarize the result of this section by abstracting a variance bound (Cor.~\ref{c:general.var} below) which applies to a general class of factor specifications $\phi$ on $(d,k)$-regular graphs.

Consider forming a random $(d,k)$-regular graph on $n$ variables by adding clauses randomly one at a time, and for $i\le m-k$ let $G^\circ$ denote the graph with the first $i-1$ clauses. Randomly partition the remaining unmatched half-edges into disjoint sets $\KK$ and $\WW$ with $|\KK|=k^2$. 
Let $B^\circ_\ell(\KK)$ denote the ball of graph distance $\ell$ about $\KK$ in the graph $G^\circ$; the leaves of $B^\circ_\ell(\KK)$ are half-edges incident to variables ($\ell$ odd) or clauses ($\ell$ even). As in \eqref{e:TT} set
$$
\TT=B^\circ_{2t'}(\KK),\quad
t' = t\wedge \min\set{\ell: B^\circ_{2\ell} \text{ has fewer than $|\KK|$ connected components} }.
$$
Let $\UU$ denote the leaves of $\TT$ without $\WW$, and $\YY$ the disjoint union of $\UU$ and $\WW$.

Let $\Ak,\Aprime$ be two arbitrary ways to form $k$ clauses on $\UU$.
Let $\ka(\usi_\UU)$ denote the partition function on $\TT\cup\Ak$ subject to boundary conditions $\usi_\UU$, and let $\ka^\wedge_\emptyset\equiv\sum_{\usi_\UU}\bp(\usi_\UU)\ka(\usi_\UU)$. Define similarly $\kaprime(\usi_\UU)$ and $\kaprime^\wedge_\emptyset$ with respect to $\Aprime$ in place of $\Ak$. Let $\psiw(\usi_\WW)$ denote the probability that $\usi_\WW$ is valid with respect to a random formation of clauses on $\WW$. 

\bcor\label{c:general.var}
Suppose $\phi\equiv(\dpsi,\hpsi)$ specifies a factor model on $(d,k)$-regular bipartite factor graphs such that which the following hold:
\bnm[(i)]
\item (Factor support) The space $\msg$ of spins
is connected by measures $\bde_{\si-\si'}$ on the support of $\phi$
 (in the sense of Lem.~\ref{l:lin.bij}); likewise the space $\msg^2$ of pair spins
is connected by measures 
$\bde_{\tau-\tau'}$ on the support of the second-moment factors
 $\psit$.
\item (Moment conditions) The first-moment rate function $\bPhi$ has negative-definite Hessian at its global maximizer $\bhstar$, and the second-moment rate function $\bPhit$ has a local maximum at $\bhtstar\equiv\bhstar\otimes\bhstar$ with negative-definite Hessian.
\item (Tree isomorphisms) On the event $\mathbf{T}$ that $\TT$ consists of $|\KK|$ tree components with $\TT\cap\WW=\emptyset$, or the event $\mathbf{C}^\circ$ that $\TT$ either contains a single cycle or has a single intersection with $\WW$ (but not both), $\ka^\wedge_\emptyset=\kaprime^\wedge_\emptyset$. Further $\ka^\wedge_\emptyset|_{\mathbf{T}}$ takes a constant value $\emptyka$ not depending on $\Ak$,
and $\Cov_{\bm{p}}(\ka-\kaprime,f)=0$ for any function $f$ depending only on a single spin $\si_e$, $e\in\UU$.\footnote{Note that this property is not immediate in the \acr{nae-sat} setting because $\Ak$ and $\Aprime$ may have different literals, but in a model with non-random factors (e.g.\ the hard-core model) it follows immediately from the isomorphism between $T\cup\Ak$ and $T\cup\Aprime$.}
\item (Correlation decay) The tree Gibbs measure $\nu$ corresponding to $\bhstar$ has correlation decay at rate faster than the square root of the tree's branching rate: for variables $u,v$ are separated by $t$ clauses, $|\Cov_\nu(f(\si_u),g(\si_v))|
	\lesssim_k\|fg\|_\infty [c (d-1)(k-1)]^{-t/2}$
for $c>1$.
\enm
Let $r>0$ such that $\bPhitstar(\bh)<\bPhitstar(\bhtstar)$ for all $\bh\ne\bhtstar$ within distance $2r$ of $\bhtstar$, and let $\idZZpd$ refer to the contribution to the pair partition function on $\Gmin\equiv G^\circ\setminus\TT$ from empirical measures at distance more than $r$ from $\bhtstar$. For $m-m'\lesssim_k\log n$ and $m'\le i\le m-k$ define
$$\TS
\vid^t_i
\equiv
\E\big[
\sum_{\uta_\YY}
	\psiw(\usi^1_\WW)\psiw(\usi^2_\WW)
	\wt\ka(\usi^1_\UU)
	\wt\ka(\usi^2_\UU)
	\idZZpd[\uta_\YY]\big],\quad
	\wt\ka
	\equiv
	\ka+\kaprime+(\ka+\kaprime)^\wedge_\emptyset.$$
If $\starZZ$ denotes the contribution to the partition function from empirical measures within distance $n^{-1/2}\log n$ of $\bhstar$, then
for
$t=t(k,\ep)= 4\log_c(1/\ep)$ we have
$$\TS
\Var\log(\starZZ+\ep\,\EZZ)
\lesssim_k 
	1 + o_\ep(1)
	 +(\ep\,\EZZ)^{-2}
	 	[m'\vid^t_{m'}+ \sum_{i=m'+1}^{m-k} \vid^t_i].$$
In particular, if $\bPhi(\bhstar)>0$ and
 $\bhtstar$ is the unique global maximizer of $\bPhit$ on $\simplext$, then $\Var\log(\starZZ+\ep\,\EZZ)
\lesssim_k 
	1 + o_\ep(1)$, hence $\lim_{n\to\infty}\P(\starZZ>0)=1$.
\ecor

\bpf[Proof of Thm.~\ref{t:whp}\ref{t:whp.below}]
The original \acr{nae-sat} model can also be regarded as a factor model in the sense of Cor.~\ref{c:general.var}, with factors
 $\dpsi(\dsi)$ and $\hpsi_a(\ksi_a)\equiv\hpsi^\circ(\ksi_a\oplus \ulit_a)$ where (compare~\eqref{e:zof})
$$
\dpsi(\dsi)\equiv
\text{\footnotesize$\begin{cases}
1,& \dsi=(\zro^d) \text{ or } (\one^d),\\
0,& \text{else;}
\end{cases}$},\quad
\hpsi^\circ(\ksi)
\equiv\text{\footnotesize$\begin{cases}
1, & \ksi\in\perm[ (\zro^j,\one^{d-j})_{1\le j\le d-1} ]\\
0, & \text{else.}
\end{cases}$}
$$
For $d\le \dlbd$, the moment condition of Cor.~\ref{c:general.var} holds by Propn.~\ref{p:sat.below}. The remaining conditions are easily verified so we indeed have $\lim_{n\to\infty}\P(Z>0)=1$ as claimed.
\epf

\section{From clusters to assignments}\label{s:sat}

\bpf[Proof of Thm.~\ref{t:main}]
Given an auxiliary model configuration $\usi$ on the edges of $(G,\ulit)$, our aim is to complete $\usi$ to an \acr{nae-sat} solution $\ux$ on $(G,\ulit)$ (meaning that whenever $x_v$ agrees with $\eta_v\equiv\dmp_d(\dsi_v)$ whenever $\eta_v\ne\free$). Clearly, the potential issue is that setting a free variable may cause a chain of forcings resulting in an invalid assignment. We therefore let $F^\sharp\equiv F^\sharp(G,\ulit,\usi)$ denote the subset of clauses $a\in F$ such that at least two variables in $\pd a$ are free, and all rigid variables $v\in\pd a$ have the same evaluation $\lit_{av}\oplus\eta_v\equiv\xi_a$.\footnote{If $\ueta_{\pd a}=(\free^k)$ we also include $a\in F^\sharp$, and arbitrarily define $\xi_a=\zro$.}

Let $G^\sharp\equiv G^\sharp(G,\ulit,\usi)$ denote the subgraph of $G$ induced by the free variables together with the clauses $F^\sharp$. We claim that \emph{$\usi$ has a valid completion to an \acr{nae-sat} solution provided each connected component of $G^\sharp$ contains at most one cycle}. Indeed, in a tree component of $G^\sharp$ one may choose an arbitrary root vertex and assign it an arbitrary value --- this may cause a chain of forcings, but no conflict results since there is no cycle. In a unicyclic component $C$ with cycle $v_0,a_0,v_1,\ldots,a_{n-1},v_n$ (with indices taken modulo $n$ so $v_0=v_n$), setting $x_{v_i}=\neg\lit_{a_i v_i}\oplus\xi_{a_i}$ ensures that all clauses along the cycle are satisfied. Then, by the preceding argument for tree components, there exists a valid completion of $\ux$ to the remainder of $C$, proving our claim.

By Thm.~\ref{t:whp}\ref{t:whp.below} and Thm.~\ref{t:whp.ball} it suffices to show that for $k\ge k_0$ and $\dlbd\le d< d_\star$, 
the limit
$\lim_{n\to\infty}\P(Z>0\giv\ZZ(\bh)>0)=1$
holds uniformly over empirical measures
$\bh\in\simplex$ with $\|(n\dbh,m\hbh)-\bgstar\|\le n^{1/2}\log n$. Conditioned on $\ZZ(\bh)>0$ we may generate $(G,\ulit,\usi)$
--- where $G$ has the law of $\Gndk$, $\ulit$ is uniformly random, and $\usi$ has empirical measure $\bh$ ---
 as follows: 
start with a set $V$ of $n$ variables each incident to $d$ half-edges
and a set $F$ of $m$ clauses each incident to $k$ half-edges, and place spins on half-edges according to $\dbh$ and $\hbh$. Then construct 
the graph by randomly matching clause and variable half-edges in breadth-first search manner started from an initial variable $\rt$, and respecting the given spins $\si$. It is clear from this construction that up to the time that the process has explored say $n^{1/3}$ vertices, the evolution of the spins $\usi$ on the leaves of the exploration tree is very close to the Markovian evolution of the Gibbs measure $\nuaux$ described in \S\ref{ss:one.reduc}. In particular, 
starting from any free variable $v$,
the exploration of its connected component $T_v$ in $G^\sharp$ is
dominated by a Galton--Watson branching process with offspring numbers distributed as a random variable $0\le Y\le dk$ with
$$\E Y \lesssim
dk \, \hbhstar(F^\sharp \giv \si_1=\ff)
= \f{dk\sum_{j=2}^k \tbinom{k-1}{j-1}
	\hd_\ff^j [\hd_\zf^{k-j}+\hd_\of^{k-j}]}
		{ \hbz \vh(\ff)  }
	\lesssim k^3/2^k.$$
By a standard argument the total size of the Galton--Watson tree has an exponential tail,\footnote{Suppose $Y$ is a non-negative integer random variable with $\Lm(\lm)\equiv\log\E e^{\lm Y}<\infty$ for some $\lm>0$, and $\Lm'(0)=\E Y<1$. Let $(Y_j)_{j\ge1}$ be a sequence of i.i.d.\ random variables distributed as $Y$, and $Z_n\equiv 1 + \sum_{j=1}^n(Y_j-1)$. Then the total size of a Galton--Watson tree with offspring distribution $Y$ has the same law as $\tau\equiv\inf\set{n:Z_n=0}$, and it is clear that the distribution of $\tau$ has exponential decay: $\P(\tau>j)\le\P(Z_j\ge1)\le e^{-t}\E e^{t Z_j} = e^{j[\Lm(t)-t]}$, and since $\Lm'(0)=\E Y<1$, by considering $t>0$ sufficiently small we can find a constant $c>0$ such that $\P(\tau>j)\le e^{-cj}$.} so we may take $C\equiv C(k)$ such that $\P(|T_v|\ge C\log n)\le n^{-10}$. The probability of seeing more than one cycle in $T_v$ is then crudely $\le n^{-3/2}$. Taking a union bound over all free variables shows that w.h.p.\ no component of $G^\sharp$ contains more than a single cycle, so $\usi$ corresponds to a true \acr{nae-sat} solution as claimed.
\epf


The above analysis completes the analysis of the \textsc{sat--unsat} transition in the case that the critical threshold $d_\star$ (see Propn.~\ref{p:explicit}) is non-integer. We conclude by briefly sketching a proof that if $d_\star\in\Z$, then at $d=d_\star$ the probability that a random \acr{nae-sat} instance $(\Gndk,\ulit)$ is solvable is \emph{asymptotically bounded away from zero and one}.

\bpf[Proof for case of $d_\star$ integer
	(sketch)]
That the probability of solvability
is bounded away from zero follows by carrying out a somewhat more careful second moment argument to remove the $n^{O(1)}$ factor appearing in Thm.~\ref{t:moments}. To see that the probability is bounded away from one, it suffices to show that $\E[\ZZ \giv \Om_C] < 1$ for an event $\Om_C$ of asymptotically positive probability. We shall take $\Om_C$ to be the event that there is a large (but constant) number of disjoint triangles in the graph.  We show below that each additional triangle decreases the expected partition function by a constant factor, so that $\E[\ZZ \giv \Om_C] < 1$ for a sufficiently large (but constant) number of cycles. It is well-known that the number of triangles is asymptotically a non-degenerate Poisson random variable, so $\Om_C$ has asymptotically positive probability as required.

\begin{figure}[h]
\centering
\begin{subfigure}[h]{.32\textwidth}
\includegraphics[width=\textwidth,trim=0in 0.2in 0in 0in,clip]{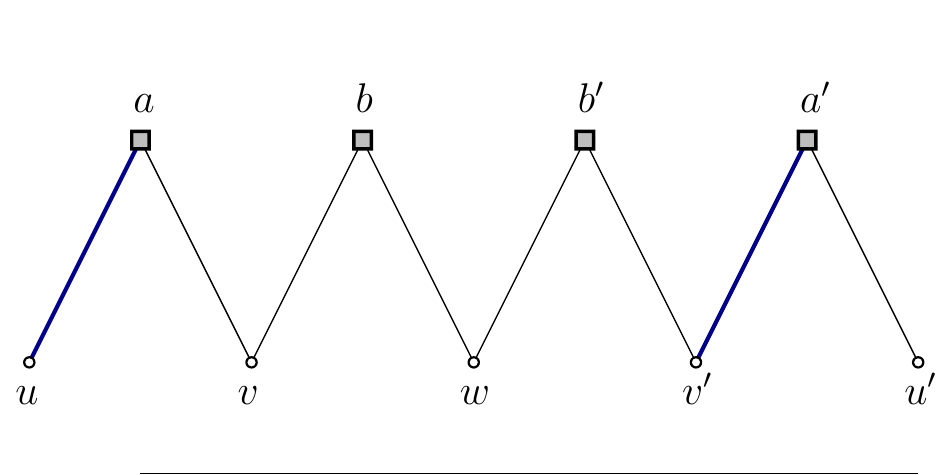}
\caption{Local neighborhood in $G_{(\ell-1)}$}
\label{f:switching.cycles.line}
\end{subfigure}
\qquad
\begin{subfigure}[h]{.32\textwidth}
\includegraphics[width=\textwidth,trim=0in 0.2in 0in 0in,clip]{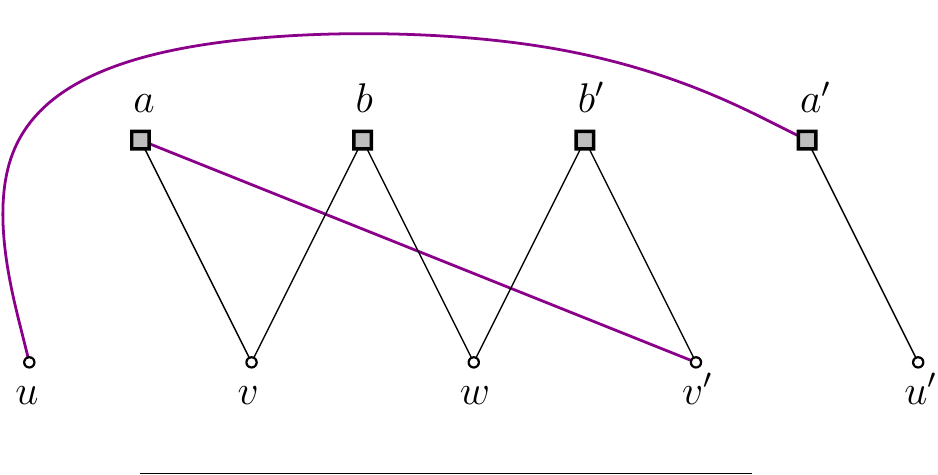}
\caption{Switched neighborhood in 
	$G_{(\ell)}$}
\label{f:switching.cycles.cycle}
\end{subfigure}
\caption{Switching argument for case of
	integer-valued $d_\star$}
\label{f:switching.cycles}
\end{figure}

We  define recursively a sequence of graphs 
$(G_{(\ell)})_{\ell\ge0}$ by the so-called ``switching method.'' Start from $G_{(0)}\equiv\Gndk$ ($d=d_\star$).
For $\ell\geq 1$, let $v,v'$ be a random pair of vertices at distance two in
the hypergraph $G_{(\ell-1)}$, with common neighbor $w$. 
Say $v$ is joined to $w$ by clause $b$,
and let $u\ne w$ be another neighbor of $v$
via a different clause $a\ne b$.
Likewise say $v'$ is joined to $w'$ by clause $b'$, and let $u'\ne w'$ be another neighbor of $v'$ via a different clause $a'\ne b'$
(Fig.~\textsc{\ref{f:switching.cycles.line}}).
Let $G_{(\ell)}$ be defined by making the switching shown in 
Fig.~\textsc{\ref{f:switching.cycles.cycle}}.
The result will follow by showing that
for $\ell$ bounded by a large constant,
this
switching decreases the expected partition function by a constant factor.

\begin{figure}[h]
\includegraphics[width=.5\textwidth]{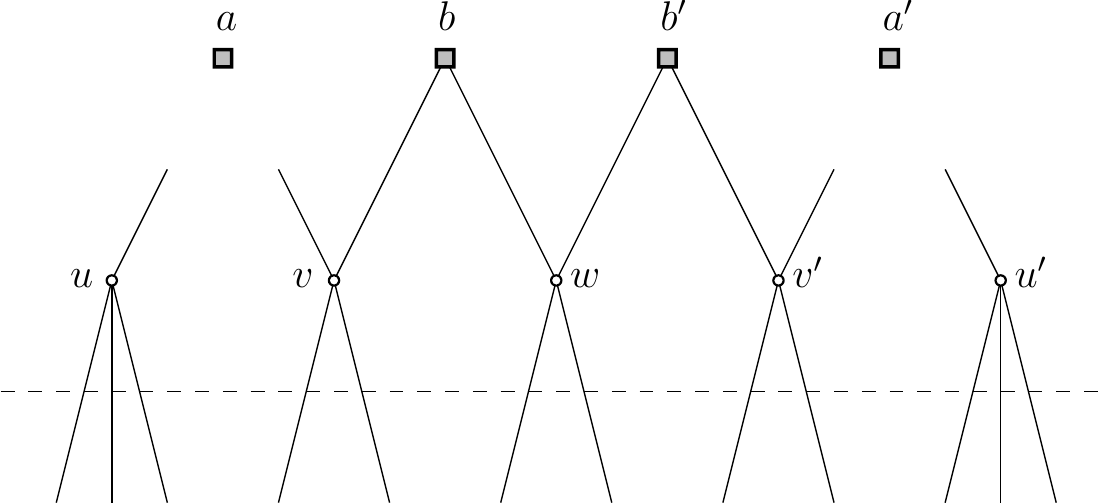}
\caption{Local neighborhood
	in graph $G_{(\ell)}^\star$
	with $a,a'$ removed}
\label{f:switching.gstar}
\end{figure}

Note that with high probability all previous switchings occur at distance at least say $(\log n)^{1/2}$ away, so it suffices to prove the claim with $\ell=1$. Consider the graph $G^\star$ with the clauses $a$ and $a'$ removed, leaving unmatched half-edges incident to variables (Fig.~\ref{f:switching.gstar}). Write $\P^\star$ for the marginal law, with respect to the auxiliary model on $G^\star$, for the spins $\si_u,\si_v,\si_{u'},\si_{v'}$ on the unmatched half-edges incident to $u,v,u',v'$; and write each $\si$ as $\bm{i}\bm{o}$ where $\bm{i}$ is the clause-to-variable message while $\bm{o}$ is the variable-to-clause message. (For example, $\bm{o}_u$ will correspond to $\si_{u\to a}$ in the original graph versus $\si_{u\to a'}$ in the switched graph.) We shall compare the probability for $a$ and $a'$ to be satisfied within the original graph versus the switched graph: with $\rig\equiv\set{\zro,\one}$ as above,
we claim
\begin{align}
\label{e:gstar.compare.a}
\P^\star( \bm{o}_v\in\rig,\bm{o}_{v'}\in\rig )
	\,[1+o(1)]
	&\ge 
	\P^\star( \bm{o}_v\in\rig,\bm{o}_u\in\rig )
	+ 2^{-5k},
	\\
\label{e:gstar.compare.aprime}
\P^\star( \bm{o}_u\in\rig,\bm{o}_{u'}\in\rig)
	\,[1+o(1)]
	&\ge
	\P^\star( \bm{o}_{v'}\in\rig,\bm{o}_{u'}\in\rig ).
\end{align}
(In the above, the first display concerns clause $a$ while the second concerns $a'$; in both displays the left-hand side is relevant to the switched graph while the right is relevant to the original graph.)
Recalling Lem.~\ref{l:dangling},
$\P^\star$
is the same up to $1+o(1)$ factors
as the measure $\mathbf{P}$ induced on the local neighborhood by taking boundary conditions given by $\bp$
(on the edges cut by the dashed line in Fig.~\ref{f:switching.gstar},
without regard to the structure of $G^\star$).
Under $\mathbf{P}$, clearly $\bm{o}_u$, $\bm{o}_{u'}$, and $\bm{o}_{v'}$ are mutually independent. Since $v'$ has only $d-2$ neighbors coming from the rest of the graph, it is slightly biased towards $\free$, which proves \eqref{e:gstar.compare.aprime}.

To prove \eqref{e:gstar.compare.a} we need to take two effects into account: first, $\bm{o}_v$ and $\bm{o}_{v'}$ are correlated while $\bm{o}_v$ and $\bm{o}_u$ are independent; and secondly, as noted in the proof of \eqref{e:gstar.compare.aprime}, marginally $\bm{o}_{v'}$ is slightly more likely than $\bm{o}_u$ to be $\free$ due to the different structure of the local neighborhood. The correlation goes in our favor while the marginal bias goes against, and we argue that the former dominates. Indeed, as we have seen in the proof of Lem.~\ref{l:tree.cov}, there is an event of probability $\gm\asymp 2^{-2k}$ such that on this event $\bm{o}_v$ and $\bm{o}_{v'}$ must be both rigid (with probability $z$) or both free (with probability $z$), but given the complementary event they are conditionally independent with probability $x$ for $\bm{o}_v$ to be rigid and probability $y$ for $\bm{o}_{v'}$ to be rigid. Thus $\mathbf{P}(\bm{o}_v\in\rig)=\gm z+(1-\gm)x$ which implies $x_\free\equiv 1-x \asymp 2^{-k}$; likewise and $\mathbf{P}(\bm{o}_{v'}\in\rig)=\gm z+(1-\gm)y$ which implies $y_\free\equiv 1-y \asymp 2^{-k}$. Combining, $\mathbf{P}(\bm{o}_v,\bm{o}_{v'}\in\rig)-\mathbf{P}(\bm{o}_v\in\rig) \mathbf{P}(\bm{o}_{v'}\in\rig)$ is quadratic in $z$, and it is straightforward to compute the derivative
and see that it is $\asymp\gm$ (hence positive) for $0\le z\le 1$.
Evaluating at $z=1$ gives
\[
\mathbf{P}(\bm{o}_v,\bm{o}_{v'}\in\rig)-\mathbf{P}(\bm{o}_v\in\rig) \mathbf{P}(\bm{o}_{v'}\in\rig)
\ge \gm(1-\gm)(1-x)(1-y) \asymp 2^{-4k}.
\]
As for the marginal bias, note that the increased chance for $\bm{o}_{v'}$ to be free compared with $\bm{o}_u$ comes from the fact that $v$ receives only $d-2$ incoming messages from the rest of the graph: thus $\si_{v\to b}$ is slightly biased towards $\free$, and this effect can percolate through the chain $\si_{b\to w}$, $\si_{w\to b'}$, $\si_{b'\to v'}$ to affect $\bm{o}_{v'}$. However the initial bias on $\si_{v\to b}$ is $\lesssim 4^{-k}$, and the effect decreases by a factor $2^k$ passing through each step in the chain, so the overall bias is $\lesssim 2^{-6k}$. Combining these estimates proves \eqref{e:gstar.compare.a}.

The result follows from \eqref{e:gstar.compare.a} and \eqref{e:gstar.compare.aprime} by noting that in a clause with $k$ random incoming messages which are mutually independent except for possible correlation among the first two, the probability for the clause to be satisfied decreases if the probability for the first two messages to be both rigid increases.
\epf

\bibliographystyle{abbrv}
\bibliography{refs}

\end{document}